%% file: nc-bv23.tex
\documentclass{birkjour}
\usepackage{amsmath, amssymb, epic, eepic}
\usepackage[all]{xypic}
\newtheorem{thm}[equation]{Theorem}
\newtheorem{prop}[equation]{Proposition}
\newtheorem{cor}[equation]{Corollary}
\newtheorem{lemma}[equation]{Lemma}

\newtheorem{lemdfn}[equation]{Lemma-Definition}

\theoremstyle{definition}
\newtheorem{defn}[equation]{Definition}
\newtheorem{caut}[equation]{CAUTION}

\theoremstyle{remark}
\newtheorem{exam}[equation]{Example}

\newtheorem{ntn}[equation]{Notation}
\newtheorem{rem}[equation]{Remark}
\makeatletter
\renewcommand{\subsection}{\@startsection{subsection}{2}{0pt}{-3ex
plus -1ex minus -0.2ex}{-2mm plus -0pt minus
-2pt}{\normalfont\bfseries}} 
\renewcommand{\subsubsection}{\@startsection{subsubsection}{2}{0pt}{-3ex
plus -1ex minus -0.2ex}{-2mm plus -0pt minus
-2pt}{\normalfont\bfseries}} \makeatother

\numberwithin{equation}{section}
\numberwithin{equation}{subsection}

\newdir^{ (}{{}*!/-5pt/@^{(}}
\def\hp{\hphantom{x}}
\newcommand{\mmod}{\text{-}\mathrm{mod}}
\newcommand{\erem}{\hfill$\lozenge$\end{rem}\vskip 3pt }
\newcommand{\Df}{\mathcal{D}}
\newcommand{\Ff}{\mathcal{F}}
\newcommand{\ev}{\operatorname{\mathsf{ev}}}
\renewcommand{\div}{\operatorname{div}}
\newcommand{\ctr}{\lrcorner\,\,}

\newcommand{\iso}{{\;\stackrel{_\sim}{\to}\;}}

\newcommand{\beq}{\begin{equation}\label}
\newcommand{\eeq}{\end{equation}}

\newcommand{\into}{\hookrightarrow}
\newcommand{\too}{\longrightarrow}
\newcommand{\onto}{\twoheadrightarrow}

\newcommand{\pr}{\operatorname{pr}}
\newcommand{\GL}{\operatorname{\mathsf{GL}}}

\newcommand{\ds}{\displaystyle}

\newcommand{\DDer}{\mathbb{D}\mathrm{er}}
\newcommand{\DDA}{\mathbb{D}\mathrm{er}(A)}
\newcommand{\WDDA}{\mathcal{W}\mathbb{D}\mathrm{er}(A)}
\newcommand{\cW}{\mathcal{W}}
\newcommand{\cM}{\mathcal{M}}
\newcommand{\cP}{\mathcal{P}}
\newcommand{\cQ}{\mathcal{Q}}
\newcommand{\cA}{\mathcal{A}}
\newcommand{\cV}{\mathcal{V}}
\newcommand{\Der}{\operatorname{Der}}

\newcommand{\op}{{\operatorname{op}}}
\newcommand{\inn}{{\operatorname{inn}}}
\newcommand{\out}{{\operatorname{out}}}

\newcommand{\Hom}{\operatorname{Hom}}
\newcommand{\tra}{\mathop{\rightarrow}}
\newcommand{\traa}{\mathop{\longrightarrow}}
\newcommand{\Spec}{\operatorname{Spec}}

\newcommand{\en}{\enspace }
\newcommand{\mto}{\mapsto}

\newcommand{\vi}{${\en\mathsf {(i)}}\;$}
\newcommand{\vii}{${\;\mathsf {(ii)}}\;$}

\renewcommand{\mto}{\mapsto}
\newcommand{\La}{{\mathsf{\Lambda}}}

\def\O{\mathcal{O}}

\def\k{\mathbf{k}}

\def\o{\otimes}
\def\ad{\mathrm{ad}}

\def\End{\mathrm{End}}
\def\dq{{\overline{Q}}}

\def\Id{\mathrm{Id}}
\def\cyc{\mathrm{cyc}}

\def\gr{\operatorname{gr}}

\def\rk{\mathrm{rk}}
\def\dim{\operatorname{\mathrm{dim}}}

\def\Rep{{\operatorname{Rep}}}
\def\RA{{\operatorname{RA}}}
\def\WRep{{\operatorname{WRep}}}
\def\WEnd{{\operatorname{WEnd}}}
\def\Ind{\operatorname{Ind}}

\def\Sym{{\text{Sym\ }}}
\def\Syns{{\operatorname{Sym}}}
\def\SuperSym{\operatorname{SuperSym}}
\def\Z{{\mathbb Z}}
\def\tr{{\text{tr}}}

\def\1{\mathbf{1}}

\def\ldb{\mathopen{\{\!\!\{}} \def\rdb{\mathclose{\}\!\!\}}}

\def\ldp{\mathopen{\langle\!\langle}} \def\rdp{\mathclose{\rangle\!\rangle}}

\def\SS{\mathbb{S}}
\begin{document}
\title{Differential operators and BV structures in noncommutative geometry} 
\author{Victor Ginzburg}
\address{Dept. of Math., U. Chicago \\ 5734 S. University Ave\\ Chicago, IL 60637}
\email{ginzburg@math.uchicago.edu}

\author{Travis Schedler} 
\address{Dept. of Math., MIT, Rm 2-236 \\ 77 Massachusetts Ave\\ Cambridge, MA 02139}
\email{trasched@gmail.com}

\begin{abstract} We introduce a new formalism of differential
  operators for a general associative algebra $A$. It replaces
  Grothendieck's notion of differential operators on a commutative
  algebra in such a way that derivations of the commutative algebra
  are replaced by $\DDA$, the bimodule of double derivations.  Our
  differential operators act not on the algebra $A$ itself but rather
  on $\Ff(A)$, a certain `Fock space' associated to any noncommutative
  algebra $A$ in a functorial way. The corresponding algebra
  $\Df(\Ff(A))$ of differential operators is filtered and
  $\gr \Df(\Ff(A))$, the associated graded algebra, is commutative in
  some `wheeled' sense.  The resulting `wheeled' Poisson structure on
  $\gr \Df(\Ff(A))$ is closely related to the double Poisson structure
  on $T_A \DDA$ introduced by Van den Bergh.  Specifically, we prove
  that $\gr \Df(\Ff(A))\cong \Ff(T_A(\DDA)),$ provided the algebra $A$
  is smooth.

  Our construction is based on replacing vector spaces by the new
  symmetric monoidal category of \emph{wheelspaces}.  The Fock space
  $\Ff(A)$ is a commutative algebra in this category (a ``commutative
  \emph{wheelgebra}'') which is a structure closely related to the
  notion of wheeled PROP.  Similarly, we have Lie, Poisson, etc.,
  wheelgebras. In this language, $\Df(\Ff(A))$ becomes the universal
  enveloping wheelgebra of a Lie wheelgebroid of double derivations.

  In the second part of the paper we show, extending a classical
  construction of Koszul to the noncommutative setting, that any
  Ricci-flat, torsion-free bimodule connection on $\DDA$ gives rise to
  a second order (wheeled) differential operator, a noncommutative
  analogue of the Batalin-Vilkovisky (BV) operator, that makes
  $\Ff(T_A(\DDA))$ a BV wheelgebra.

  In the final section, we explain how the wheeled differential
  operators $\Df(\Ff(A))$ produce ordinary differential operators on
  the varieties of $n$-dimensional representations of $A$ for all $n
  \geq 1$.
\end{abstract}

\maketitle

\bigskip

\centerline{Table of Contents}
\vskip -1mm
{\parbox[t]{115mm}{
\noindent
\hp1.{ $\;\,$} {Introduction}\newline
\hp2.{ $\;\,$} {Differential operators in noncommutative settings}\newline
\hp3.{ $\;\,$} {Differential operators on associative algebras and wheelgebras}\newline
\hp4.{ $\;\,$} {Torsion of bimodule connections on $\DDer(A)$ and $\Omega^1A$}\newline
\hp5.{ $\;\,$} {The BV operator $D_\nabla$}\newline
\hp6.{ $\;\,$} {The representation functor}
}}

\bigskip
\section{Introduction} 
\subsection{} The general definition of differential operator on an
abstract commutative algebra was introduced by Grothendieck.  For
noncommutative algebras, such as tensor algebras, Grothendieck's
definition still makes sense, but does not necessarily lead to a good
notion (see Remark \ref{groth}, \ref{groth2} for details).

In this paper, we introduce a new notion of differential operators on
associative algebras. Our approach is based on the observation that
the tensor algebra of a vector space may be viewed as a \emph{twisted}
commutative algebra.  The notion of twisted commutative algebra dates
back at least to the 1950's, appearing in algebraic topology, cf. also
\cite{Bar}.  One way to think about twisted commutative algebras is to
interpret them as {\em commutative algebra objects} in the monoidal
category of $\SS$-modules, where an $\SS$-module is a graded vector
space equipped with symmetric group actions on its homogeneous
components, cf. \S\ref{twist_sec}.

As a first step of our construction, we introduce
a rather general definition of  differential operators
for a  commutative  algebra object
in  an \emph{arbitrary} abelian symmetric monoidal
category. Applying our definition
in the special case of the monoidal
category  of $\SS$-modules, one obtains quite a reasonable theory of
differential operators on a twisted commutative algebra. Thus,
for any twisted commutative 
algebra $A$,  our construction produces
 a filtered algebra $\Df(A)$ of differential
operators.  We show that $\gr\Df(A)$, the associated graded algebra,
is twisted commutative and, moreover, it has a natural structure of a
twisted Poisson algebra (Theorem \ref{dot}).  Twisted Poisson
structures on tensor algebras are related to Van den Bergh's double
Poisson structures \cite{VdB}, as explained in \cite{Stw}.

The theory of differential operators on twisted commutative algebras
is not sufficient for the purposes of noncommutative geometry:
it is essential to be able to consider differential operators on
tensor algebras over a general \textit{noncommutative} base. This
amounts, formally, to replacing various tensor products over the
ground field $\k$ by tensor products over a noncommutative algebra
$A$. If $A$ is the path algebra of a quiver, for instance, tensoring
over $A$ may be interpreted as a gluing operation that joins heads and
tails of various paths. Notice that joining the head and tail of the
 {\em same} path creates
a ``wheel,'' i.e., an oriented cycle in the quiver. 
Thus, the next step of our construction is to introduce
the concept of
a {\em wheeled differential operator}.

The general point of view of this paper is that an adequate framework
for doing noncommutative differential geometry is provided by the
notion of wheelspace. Wheelspaces form a symmetric monoidal category,
and there are natural notions of commutative, Poisson, and Lie algebra
objects in that category, to be called commutative, Poisson, and Lie
{\em wheelgebras}, respectively.  Commutative wheelgebras turn out to
be special cases of \textit{wheeled PROP}s, introduced recently in
\cite{MMS}, which motivates the name.  Given any noncommutative
associative algebra $A$, we define a commutative wheelgebra, $\Ff(A)$,
that resembles the Fock space construction.  Having $m$ inputs and $n$
outputs for a wheeled PROP translates in our setting as having an
$A^{\o m}$-$A^{\o n}$-bimodule structure.

We then apply our abstract categorical definition of differential operator
to the special case of the  monoidal category of wheelspaces.
This way, any
commutative wheelgebra $\cW$ gives rise to a filtered associative
wheelgebra, $\Df(\cW)$, of differential operators on $\cW$. In particular,
differential operators from $\Df(\Ff(A))$ act on $\Ff(A)$ very much
like Heisenberg algebras of creation and annihilation operators act on
Fock spaces.

In the case that $A$ is smooth, it turns out that the associated
graded $\gr \Df(\Ff(A)) \cong \Ff(T_A \DDA)$ (Theorem \ref{dathm}),
and the induced Poisson wheelgebra structure is closely related to the
double Poisson bracket introduced by Van den Bergh as a noncommutative
counterpart of the canonical Poisson structure on $T^*X$, the total
space of the cotangent bundle on a manifold $X$. This further explains
the connection from \cite{Stw} between double Poisson structures and
twisted Poisson structures on tensor algebras.

\subsection{}
In the second half of the paper we transplant various differential
geometric structures related to the notion of
Batalin-Vilkovisky (BV) algebra and BV operator to noncommutative geometry.
To explain this, first let $X$ be a smooth (complex) manifold.
It is well known  that the Schouten-Nijenhuis bracket makes
the space $\La^\bullet T_X$ of polyvector fields on 
 $X$  a Gerstenhaber algebra. 

Assume next that $X$ is a Calabi-Yau manifold of dimension $d$.
Any flat connection on the line bundle $\La^d T_X$ 
 gives rise to
 an odd second-order differential
operator $D: \La^\bullet T_X\to\La^{\bullet-1} T_X$ such that
$D^2=0$;  see \cite{Scn}. The operator $D$,
called the BV operator, satisfies the BV identity
\begin{equation} \label{bvid}
  (-1)^{|\xi| + 1} \{\xi, \eta \} = D(\xi \wedge \eta) - 
D(\xi) \wedge \eta - (-1)^{|\xi|} \xi \wedge D(\eta),
\end{equation}
for $\xi, \eta$ local sections of $\La^\bullet T_X$.
Thus, $D$ gives the  Gerstenhaber algebra $\La^\bullet T_X$
the structure of a BV algebra.
Following Koszul \cite{Kz}, we observe further that any Ricci-flat connection on
 the tangent bundle gives a flat connection on the line bundle $\La^d
 T_X$,
hence gives rise to a BV operator $D$ as above.

In  noncommutative geometry,  vector fields are
replaced by  $\DDA := \Der(A, A \o A)$, the bimodule of
\textit{double derivations} of a  noncommutative algebra $A$;
see \cite{CB}, \cite{VdB}, \cite{CBEG}.
The supercommutative algebra
 $\La^\bullet T_X$ is therefore replaced by  $T_A \DDA$,
the tensor algebra of  the bimodule $\DDA$.
Note, however, that the component  $\La^d T_X,$ of polyvector fields of
top degree, has no  noncommutative counterpart.
Thus, it is not clear {\em a priori} how to extend the Calabi-Yau
geometry outlined above to the  noncommutative setting.

With this in mind, we adapt Koszul's approach and start with a
\textit{bimodule} connection $\nabla$ on $\DDA$, as
defined in \cite{CQ}. We study torsion and  curvature  for such
connections, and also define the
`trace of curvature'.
Thus, we get the notion of a Ricci-flat  bimodule connection  on
$\DDA$. 

At this point we invoke our theory of differential operators on
commutative wheelgebras.  Our main construction associates with any
Ricci-flat, torsion-free bimodule connection $\nabla$ on $\DDA$ a
second-order differential operator, $D_\nabla$, on the commutative
wheelgebra $\Ff(T_A \DDA)$.  Furthermore, we show that such an
operator gives a BV wheelgebra structure that extends the double
Schouten-Nijenhuis structure considered by Van den Bergh \cite{VdB}.

\subsection{The representation functor}
It should be pointed out that the classical theory of differential
operators on a commutative algebra is not a special case of our
theory: differential operators on $\Ff(A)$ for a commutative algebra
$A$ are in general different from usual differential operators on $A$
(although, if we view $A$ itself as a commutative wheelgebra
concentrated in degree zero, i.e., having no inputs or outputs, then
differential operators on it in our sense are the same as in
Grothendieck's sense).

The relationship between  commutative and noncommutative theories 
is provided by the notion of
\textit{representation functor}.
Specifically, associated with a noncommutative algebra $A$, there is
a sequence  $\Rep_{\mathbf{d}}(A),  \  {\mathbf{d}}=1,2,\ldots,$ of
affine schemes (in the conventional sense) parametrizing $
{\mathbf{d}}$-dimensional representations
of $A$. According to  a philosophy advocated by Kontsevich,
any reasonable associative (noncommutative) notion for $A$ should go,
via the representation functor,
to the corresponding usual commutative notion.
Furthermore, the two stories `merge' asymptotically as
${\mathbf{d}}\to\infty$.   

In \S \ref{repsec} we extend the representation functor to act on
wheelgebras and wheeled differential operators.  This way, Theorem
\ref{quivthm} ensures in particular that the action of wheeled
differential operators corresponds, asymptotically, to the natural
action of the quantized necklace Lie algebra from \cite{S,GS}.
Moreover, a noncommutative BV structure $D_\nabla$ on a
smooth algebra $A$ gives rise to an ordinary BV structure on each of
the schemes $\Rep_{\mathbf{d}}(A), \ {\mathbf{d}}\geq1$.


We summarize our main results:
\begin{itemize}
\item Theorems \ref{dot} and \ref{dathm}, which extend 
  the theory of commutative differential operators to the twisted-commutative and wheeled context, which apply to arbitrary associative algebras.
\item Theorem \ref{bvthm}, which proves the formula
  $D_\nabla^2 = i_{\tr(\nabla^2)}$ for smooth associative algebras, thus
  allowing us to deduce the equivalence of Ricci-flat torsion-free
  bimodule connections with wheeled BV structures.
\item Theorem \ref{quivthm}, which shows that path algebras of quivers are
  wheeled Calabi-Yau.
\item Theorem \ref{repthm}, which explains that wheeled differential
  operators $\Df(\Ff(A))$ map to actual differential operators under
  the representation functor.
\end{itemize}

\subsection{Relation to  necklace Lie bialgebras}\label{nlbasec}
Path algebras associated to quivers are free.  Hence, they may be
viewed as noncommutative analogues of flat affine spaces. The affine
space comes equipped with the trivial connection on the tangent
bundle, which is flat and torsion-free.  Correspondingly, for a path
algebra $A$, one has a trivial bimodule connection on $\DDA$, which is
flat and torsion-free.  The latter gives, by our construction, a
wheeled BV operator $D$.

There is an alternative interpretation of the BV operator $D$ for path
algebras. To explain this, we remark first that an important class of
(ordinary supercommutative) BV algebras comes from Lie bialgebras,
as in, e.g., \cite{Gcy}, \S 2.10. We briefly explain this as follows.

Let ${\mathfrak g}$ be a finite-dimensional involutive
Lie bialgebra with cobracket
$\delta\colon {\mathfrak g}\to\La ^2{\mathfrak g}$.
Given a basis $x_1, \ldots, x_n$ of   ${\mathfrak g}$, one has structure constants
for the Lie bracket  defined by the equations
$[x_i,x_j]=\sum_k c_{ij}^k\cdot x_k$  and structure
constants for the cobracket defined by the equations
$\delta(x_k)=\sum_{i,j}
\  f_{ij}^k\
x_i\wedge x_j.
$

Associated with such a  Lie bialgebra, there is a BV operator
on $\La {\mathfrak g}$
defined by the formula
\begin{multline}
D:\ 
a_1\wedge\cdots\wedge
a_n\mto \sum_{1\leq k\leq n}^n(-1)^{k-1}\cdot\delta(a_k)\wedge a_1\wedge\cdots\wedge\widehat
a_k\wedge\cdots\wedge a_n\\
+\sum_{1\leq i<j\leq n}(-1)^{i+j-1}\cdot [a_i,a_j]\wedge
a_1\wedge\cdots \wedge \widehat a_i \wedge \cdots
\wedge \widehat a_j \wedge  \cdots \wedge a_n.
\end{multline}
In coordinates, the above formula reads
\begin{equation}\label{delta}
D=\sum_{i,j,k} f_{ij}^k\cdot x_i\wedge x_j\
\frac{\partial}{\partial x_k}\ +\ \sum_{i,j,k}
c_{ij}^k\cdot x_k\
\frac{\partial}{\partial x_i}\wedge
\frac{\partial}{\partial x_j}.
\end{equation}

It turns out that the wheeled BV operator for a path algebra $A$ may
be viewed as a noncommutative analogue of the operator \eqref{delta}.
Furthermore, this is not merely an analogy.  The
wheelspace-degree-zero part of the wheeled BV operator on ${\mathcal
  F}(A)$ is, in effect, given by \eqref{delta} in the case of
$\mathfrak{g}=A_\cyc,$ where $A_\cyc:=A/[A,A]$ is a super version of
the necklace Lie bialgebra studied by one of us in \cite{S}. In the
case of quivers with one vertex, the resulting BV algebra was studied
by Barannikov (\cite{Barmobvg}, \cite[\S 1]{Barncbvg}, and subsequent
papers) in the context of quantum master equations and $A_\infty$
algebras. There, the relevant BV algebra was denoted as $F$.

In more detail, let $\pr: A \rightarrow A_\cyc, a \mto [a]$ be the
tautological projection. Further, following the strategy of
\cite{Shoch}, \S 5.2 in the ordinary (not super) case, one may lift the necklace
bracket and cobracket on $A_\cyc$ to maps
\begin{equation*}
\{-,-\}:\ A \o A_\cyc \rightarrow A, \text{ and } 
\delta: A \rightarrow A \o A_\cyc,
\end{equation*} respectively.  Then, the action
of $D$ on an element
$u=a_1 \o \cdots \o a_n \o ([b_1] \cdots [b_m])\in A^{\o
  n}\o\SuperSym_\k^m A_\cyc$ is given by the formula (if the $a_i, b_j$ are
homogeneous)
\begin{align} \label{neckbv}
  D(u) &=
  \sum_i \pm a_1 \o \cdots \o a_{i-1} \o \delta(a_i) \o a_{i+1} \o
  \cdots 
\o a_n \o \big([b_1]  \cdots  [b_m]\big)\nonumber \\ 
&+ \sum_{i,j} \pm \sigma'_{i,j} \ldb
a_i, a_j \rdb a_1 \o \cdots \o \hat a_i \o \cdots \o \hat a_j \o \cdots \o a_n \o\,
\big([b_1]  \cdots  [b_m]\big)  \nonumber \\ &+ \sum_{i,j} \pm a_1 \o \cdots \o
a_{i-1} \o \{a_i, [b_j]\} \o a_{i+1} \o \cdots \o a_n \o\, \big([b_1]  \cdots
\widehat{[b_j]} \cdots  [b_m]\big) \\
   &+ \sum_i \pm a_1 \o \cdots \o a_n \o \, \delta([b_{i}]) \cdot \big([b_1]  \cdots 
  [b_{i-1}]\cdot    [b_{i+1}]  \cdots
   [b_m]\big) \nonumber \\ &+ \sum_{i,j} \pm a_1 \o \cdots \o a_n \o\,
  \{[b_i], [b_j]\}\cdot  \big([b_1]  \cdots
  \widehat{[b_i]} \cdots \widehat{[b_j]} \cdots 
  [b_m]\big),\nonumber 
\end{align}
with signs depending on degree and indices.  The twisted-degree-zero
part of the wheeled BV identity (and the equation $D^2=0$) for $D$ is
equivalent to the fact that $(A_\cyc, \{-,-\}, \delta)$ is an
involutive Lie bialgebra.  The twisted-degree-one part of the BV
identity was first discovered in \cite{Shoch}, \S 5.3; it reads:
\begin{equation} \label{ncbvid2}
  \delta(a b) = \delta(a)(b \o 1) \pm (a \o 1) \delta(b) + (1 \o \pr) \ldb a,b \rdb.
\end{equation}

We remark that the deformations of $\delta$ (also satisfying the above
mentioned identities) described in \cite[\S 5.3]{Shoch} correspond, in
our language, to connections other than the trivial one.

\subsection{Future directions and related work}\label{futdirsec}
In the future, we hope to use wheeled differential operators to
construct a \textit{wheeled Heisenberg representation}, as well as a
wheeled version of the Weil representation.  We hope that the latter
will have also a topological analogue, consisting of a projective
action of the mapping class group of a Riemann surface on an
appropriate wheeled Fock space, along the lines of \cite[\S 6.2]{Gcy}.
Constructing such an action is a key ingredient in the approach to
Calabi-Yau algebras arising from fundamental groups of aspherical
3-manifolds, as indicated in \cite{Gcy}, Conjecture 6.2.1.
  
Also, the formalism of noncommutative BV structures seems to
play a role in trying to generalize the notion of Calabi-Yau algebra
along the lines outlined in \cite{EG2}, Remark 1.3.4.

Finally, in \cite{Barmobvg}, Barannikov constructs from any modular
operad a BV-style master equation (\cite[(5.5)]{Barmobvg}) whose
solutions are equivalent to algebras over the Feynman transform of
that operad.  When one sets the modular operad to be the operad
denoted by $\SS[t]$ in \emph{op.~cit.} (\S 9), one obtains a slight
modification of the BV algebra defined above for the case of a quiver
with one vertex, obtained by adding a parameter $t$ and keeping track
of a genus grading.

One may also form a directed analogue of the construction of
\cite{Barmobvg}, replacing modular operads with wheeled PROPs (which
include commutative wheelgebras), by replacing undirected graphs with
directed graphs.  Here, one can additionally keep track of a genus
grading at vertices.

Given a wheeled PROP $P$ with differential $d$, and a
finite-dimensional vector space $V$, one would obtain the same master
equation as in \emph{op.~cit.},
\begin{equation} \label{barbveq}
dm + z D m + \frac{1}{2} \{m, m\} = 0,
\end{equation}
where $m$ is now an element of
\begin{equation}\label{barbvalg}
\sum_{i,j,k}  P(i,j,k) \otimes_{\k[S_i \times S_j]} z^k \cdot \Hom(V^{\otimes i}, V^{\otimes j}).
\end{equation}
Here, $P(i,j,k)$ denotes the part of $P$ with $i$ inputs, $j$ outputs,
and with genus $k$, and $z$ is a formal parameter which keeps track of
the genus grading (which is \emph{not} the same as the parameter $t$
in the operad $\SS[t]$ above). The space $\Hom(V^{\otimes i},
V^{\otimes j})$ is associated to the part of the wheeled PROP with $i$
inputs and $j$ outputs, so that we can compose such operations by
plugging an output of this into the input of another such operation,
or vice-versa.  The bracket
\begin{equation}
\{\,, \}: P(i,j,k) \otimes P(i',j',k') \rightarrow
P(i+i'-1,j+j'-1,k+k')
\end{equation}
now takes a sum over all ways to contract an input of an element in
degrees $(i,j,k)$ with an output of an element in degrees
$(i',j',k')$, to obtain an element in degree $(i+i'-1,j+j'-1,k+k')$.

This is related to our bracket construction of \S \ref{nlbasec} as follows.
Consider a quiver $Q$ with one vertex, and let $V$ be the linear space
generated by the arrows, $Q$.  Consider $V^*$ to be the linear space
generated by the reverse of the arrows, $Q^*$. Then, we can then think
of $\Hom(V^{\otimes i}, V^{\otimes j}) \cong (V^*)^{\otimes i} \otimes
V^{\otimes j}$ as spanned by a pair of an $i$-tuple of reverse edges
and a $j$-tuple of edges. The bracket of two elements $F, G$ in
\eqref{barbvalg} sums over all ways of contracting an arrow of $F$
with its reverse in $G$.  This is analogous to the bracket on
$\SuperSym (\k \dq)_{\cyc}$ considered above (which, as we pointed
out, is recovered when $P = \SS[t]$ and $V = \k \dq$, except with an
extra parameter $t$).

This bracket extends linearly to a map $P \otimes P \rightarrow
P$. Similarly, the operator $D: P(i,j,k) \rightarrow P(i-1,j-1,k+1)$
of \eqref{barbveq} above takes a sum over all ways to contract an
input with an output, and this also extends linearly to all of
$P$. This operator $D$ is a PROP analogue of the wheeled BV operator
$D$ discussed in this paper.  In the one-vertex quiver setting of the
previous paragraph, the operator $D$ acts on an element $F$ of \eqref{barbvalg} by summing over ways of contracting an arrow and a reverse arrow (which becomes the BV bracket of this section in the undirected setting with $P=\SS[t]$ as
discussed above).

\subsection{Acknowledgements} We are grateful to Michel Van den Bergh
for discussions and Peter May for a careful reading and helpful
suggestions.  The first author was partially supported by NSF grant
DMS-0601050, and the second author was partially supported by an NSF
GRF, the University of Chicago's VIGRE grant, a Clay Liftoff
fellowship, and an AIM fellowship.

\subsection{Notation}
In the interests of self-containment, we repeat some standard
definitions in noncommutative geometry. For a reference, see, e.g.,
\cite{G2}.

We will work over a fixed ground field $\k$ of characteristic zero.
The term ``space'' will mean ``$\k$-vector space.''  The terms ``map''
or ``operator,'' unless otherwise specified, refer to any $\k$-linear
maps.  Unadorned tensor products will be taken to be over $\k$
throughout.

We will use cycle notation for permutations: $(a_1,a_2, \ldots, a_n)$
means $a_1 \mapsto a_2 \mapsto a_3 \mapsto \cdots \mapsto a_n$.  
\begin{ntn} \label{tauntn} For any permutation $\sigma \in
  \Sigma_{n}$,
  we define
  $\tau_\sigma: V_1 \o V_2 \o \cdots \o V_n \rightarrow
  V_{\sigma^{-1}(1)} \o V_{\sigma^{-1}(2)} \o \cdots \o
  V_{\sigma^{-1}(n)}$
  as the permutation of components corresponding to $\sigma$.
\end{ntn}
For any algebra $A$, define $A_{\cyc} := A/[A,A]$, where $[A,A]$ is
merely a vector space. When $A$ is equipped with a grading and viewed
as a superalgebra, $[\,, ]$ will be the supercommutator.

For a finitely-generated projective $A$-module $M$, 
denote the dual module $\Hom_A(M,A)$ by $M^\vee$.

For any algebra $A$ over $\k$, we write $A^{\op}$ for the opposite
algebra, and put  $A^e := A \o A^{\op}$.  Note that an
$A^e$-module is the same as an $A$-bimodule whose left and right
actions, when restricted to $\k$, become the same.  We will use
the term ``$A$-bimodule'' to refer to an $A^e$-module.  The notation
$M^\vee$ will primarily be applied to $A^e$-modules.

For any $A$-bimodules $M_1, M_2, \ldots, M_n$, one has a well defined
outer $A$-bimodule action on $M_1 \o M_2 \o \cdots \o M_n$.  If $M$ is
any $A$-bimodule, then the space
$\Hom_{A^e}(M, A^e) = \Hom_{A^e}(M, A \o A)$ is equipped with the
structure of $A$-bimodule coming from the \emph{inner} action on
$A \o A$: namely,
$(a \cdot \phi \cdot b) (m) = (1 \o a) \cdot \phi(m) \cdot (b \o 1)$.
  When we say that $M$ is a projective $A$-bimodule, we mean
  projective as an $A^e$-module. In this case,
  $M^\vee := \Hom_{A^e}(M, A^e)$.

\begin{defn}\label{ddadfn} Let
  $\DDA := \Der(A, (A \o A))$
  for any associative algebra $A$ over $\k$. Elements of $\DDA$ are
  called ``double derivations,'' hence the notation.  
\end{defn}
Similarly to $\Hom_{A^e}(M, A \o A)$, 
the space $\DDA$ is naturally an $A$-bimodule, with
$(a \cdot\phi \cdot b) (x) = (1 \o a) \cdot \phi(x) \cdot (b \o 1)$.

Following \cite{CQ}, let $\Omega^1 A := \ker(\mu_A)$ be the kernel of the
  multiplication map $\mu_A: A \o A \rightarrow A$, 
  considered as an $A$-bimodule using outer multiplication.  
By \cite{CQ}, there is a natural map $d: A \rightarrow \Omega^1 A$, and
for any $A$-bimodule $M$, there is a natural isomorphism
\begin{equation} \label{dstardfn}
d^*: \Hom_{A^e}(\Omega^1 A, M) \iso \Der(A, M), \quad 
d^* \phi = \phi \circ d.
\end{equation}
In particular, setting $M=A \o A$, one
  obtains $\Hom_{A^e}(\Omega^1 A, A \o A) \cong \DDA$.

  So, in the case that $\Omega^1 A$ is finitely-generated and
  projective, so is $\DDA$, and $\Omega^1 A$ and $\DDA$ are dual to
  each other.  The finite-generation follows if we assume that $A$ is
  finitely generated, and the fact that $\Omega^1 A$ is projective
  means precisely that $A$ is smooth.  

\section{Differential operators in noncommutative settings}

\subsection{Grothendieck's definition of differential operators}
Let $A$ be a commutative algebra. We recall:
\begin{defn} \cite{EGA4} \label{grodfn1}
The only differential operator of order $\leq -1$ is the zero
operator.  

Inductively, for any $n \geq 0$,  a linear map $D: A \rightarrow A$
is said to be  a differential operator
of order $\leq n$ if
for all $a \in A$, the operator $[D,a]: b \mapsto D(ab) - a D(b)$
is a differential operator of order $\leq n-1$.
\end{defn}
Equivalently, the above definition says that operators of order $0$ are $\Hom_A(A, A)$, the $A$-module homomorphisms from $A$ to $A$. 

There is an alternative definition to the above which coincides in the
case of commutative algebras, obtained by the same inductive procedure
except defining the order zero operators differently:
\begin{defn}\label{grodfn2}
Differential operators of order zero are left-multiplication by elements of $a$, i.e., the image of $A$ in $\Hom_\k(A, A)$.  Then, inductively, a linear map $D: A \rightarrow A$
is  a differential operator
of order $\leq n$ if
for all $a \in A$, the operator $[D,a]: b \mapsto D(ab) - a D(b)$
is a differential operator of order $\leq n-1$.
\end{defn}

In this paper, we extend the above definitions by using commutative
algebras in more general tensor categories, and so both will coincide
in these cases, and we don't need to prefer one over the other.  
\begin{rem}\label{groth} Definition \ref{grodfn1} does not
  yield a good answer in the case $A = T_\k V$, a free algebra with
  $\dim V \geq 2$: we claim that \textit{all nonzero differential
    operators have order zero}, and are \textit{right-multiplication
    by elements of $A$}. We explain this as follows.

First of all, for any $A$, $D$ has order zero if and only if
$[D,a] = 0$ for all $a$, which implies that
$ D(a) = D(a \cdot 1) = a \cdot D(1), $ for all $a$. Thus, $D$ is
right-multiplication by $D(1)$.

Next, differential operators of order $1$ are those such that
\begin{equation}
[D,a](x) = x \cdot f(a), \quad \forall a,x, \label{doo1}
\end{equation}
for some fixed linear map $f: A \rightarrow A$.  Assuming that $D(1) = 0$
by subtracting right-multiplication by $f(1)$, we get
$
D(a) = D(a \cdot 1) = f(a).
$
Hence, \eqref{doo1} becomes
\begin{equation}
D(ab) = a \cdot D(b) + b \cdot D(a), \quad  \forall a,b \in A.
\end{equation}
Then, checking associativity, we have
\begin{align}
D(abc) &= D((ab)c) = ab \cdot D(c) + c \cdot D(ab) = ab\cdot D(c) + ca\cdot D(b) + cb \cdot D(a) \\
&= D(a(bc)) = a \cdot D(bc) + bc \cdot D(a) = ab \cdot D(c) + ac \cdot D(b) + bc \cdot D(a),
\end{align}
so we deduce that
\begin{equation}
[a,c] \cdot D(b) + [b,c] \cdot D(a) = 0, \quad  \forall a,b,c \in A.
\end{equation}

Setting $b=c$, we get
$
[a,b] \cdot D(b) = 0,$ for any $a,b \in A.$
Thus, if $([A,A] \cdot x = 0) \Rightarrow (x = 0)$ (equivalently,
$(\ldp [A,A] \rdp \cdot x = 0) \Rightarrow (x = 0)$),\footnote{The same
  condition was found in \cite{Stw}, which is the condition under
  which a twisted Poisson algebra structure on $T_\k A$ (in the sense
  of this paper) must yield a double Poisson bracket on $A \o A$ (in
  the sense of \cite{VdB}).}
then we deduce that $D=0$, i.e., all differential operators of order
$\leq 1$ are of order $0$, and inductively, \textit{all} differential
operators are of order $0$.  This condition holds in particular if $A$
is not commutative and has no zero-divisors, such as $A = T_\k V$ for
$\dim V \geq 2$. It also holds if $A$ is a (deformed) preprojective
algebra \cite[Example 2.11]{Stw}.
\end{rem}
\begin{rem}\label{groth2} On the other hand, Definition \ref{grodfn2}
  yields a more interesting space of differential operators in the
  case $A = T_\k V$.  We claim that, in this case, the space of
  differential operators of order $n$ is isomorphic to $\bigoplus_{m
    \leq n} \Der(A)^{\otimes m}$, with action
\begin{multline}\label{gdodfn}
  (D_1 \otimes D_2 \otimes \cdots \otimes D_m) (v_1 v_2 \cdots v_p) \\
  = \sum_{1 \leq i_1 < i_2 < \cdots < i_m \leq p} v_1 v_2 \cdots
  v_{i_1-1} D_1(v_{i_1}) v_{i_1+1} \cdots v_{i_2-1} D_2(v_{i_2}) \cdots \\
  \cdots v_{i_m-1} D_m(v_{i_m}) v_{i_m+1} \cdots v_p.
\end{multline}
In words, the above formula roughly is the sum over all ways of
applying $D_1, D_2, \ldots, D_m$ to distinct $v_{i_1}, v_{i_2},
\ldots, v_{i_m}$, in left-to-right order (i.e., $i_1 < i_2 < \cdots <
i_m$).  Note that it is necessary here that $v_1, \ldots, v_p$ be
elements of $V$ rather than arbitrary elements of $A$.

To prove the above claim, first note that the above operations all
have order $\leq n$: this follows from the identity
\begin{equation} [(D_1 \otimes \cdots \otimes D_m),a](f) = D_1(a) (D_2
  \otimes \cdots \otimes D_m)(f).
\end{equation}
Next, any differential operator of order $\leq n$ is determined by its
restriction to products of $\leq n$ generators. In the case of $TV$,
this means the restriction to $V^{\otimes \leq n}$.  The above
operators realize all linear maps $V^{\otimes \leq n} \rightarrow TV$
uniquely, and hence \eqref{gdodfn} produces an isomorphism from
$\bigoplus_{m \leq n} \Der(TV)^{\otimes m}$ to the space of
differential operators of order $\leq n$.

As a result of the observations of the previous paragraph, in this
case of Definition \ref{grodfn2} one deduces a description of the
differential operators of order $\leq n$ on an arbitrary associative
algebra $A$: if $A$ is presented as a quotient of $TV$, then these are
the subquotient of differential operators on $TV$ of operators
preserving the kernel of $TV \onto A$, modulo operators whose image is
in this kernel.  
\end{rem}
\begin{rem}
  Even in the case of commutative algebras, the distinction between
  Definitions \ref{grodfn1} and \ref{grodfn2} can become important when
  one generalizes to the case of differential operators from an
  $A$-module $M$ to itself; then, in the first definition, the
  order-zero differential operators are $\Hom_A(M,M)$,
  while in the second they are the image of $A$ in
  $\Hom_\k(M,M)$.
\end{rem}

\subsection{Differential operators in symmetric monoidal categories}
Below, we  introduce a notion of differential operator
on a commutative algebra object in an arbitrary symmetric monoidal
category. In the case of the monoidal category of vector spaces
our definition of differential operator reduces to
the standard (Grothendieck's) definition. In this
paper, we will be mostly interested in other monoidal
categories such as the category of $\SS$-modules, cf. \S\ref{twist_sec},
or the category of wheelspaces, cf. \S\ref{wheel_sec}.

Let $(\mathcal{C}, \o, \beta)$ be a $\k$-linear symmetric monoidal
category with product $\o$ and braiding $\beta$, i.e., for any
$X, Y \in \mathcal{C}$, we have $\beta_{X,Y}: X \o Y \iso Y \o X$,
which is functorial in $X$ and $Y$.  We will have in mind the case
$\mathcal{C} = Vect$ throughout, in which case one recovers
Grothendieck's definition, as we will recall after each definition.
For simplicity, we will omit all associativity isomorphisms, and
pretend that $(X \o Y) \o Z = X \o (Y \o Z)$.\footnote{This will
  result in no loss of generality. Moreover, it is a well known result
  of MacLane that any (symmetric) monoidal category is equivalent to a
  strictly associative one.}
Let $C \in \mathcal{C}$ be a commutative algebra object in
$\mathcal{C}$ with multiplication map $\mu: C \o C \rightarrow C$.
\begin{defn} \label{admudefn}
  Let $X \in \mathcal{C}$ be an object.  A morphism
  $\phi: X \o C \rightarrow C$ is called \emph{an action of $X$ on $C$
    by differential operators of order $\leq n$}
  if the following is true: 

For $n=-1$, the morphism must be the zero
  morphism; 

and inductively on $n \geq 0$,
  $[\phi,\mu]: (X \o C) \o C \rightarrow C$ must be an action of $X \o C$
  on $C$ by differential operators of order $\leq n-1$, where
  $[\phi,\mu] := \phi \circ (\Id_X \o \mu) - \mu \circ (\Id_C \o \phi)
  \circ (\beta_{X,C} \o \Id_C)$.

  Graphically, the two terms on the RHS are given by (the diagram is NOT commutative: rather the difference of the two compositions should be a differential operator of order one lower):
\begin{equation} \label{xacte}
\xymatrix{ X \o C \o C \ar[d]_{\beta_{X,C} \o \Id_C} \ar[r]^-{\Id_X \o \mu} & X \o C \ar[r]^\phi & C \\
C \o X \o C \ar[r]^-{\Id_C \o \phi} & C \o C \ar[ur]^{\mu}           
}
\end{equation}
\end{defn}
\begin{rem}\label{catgrodfn2rem}  The above definition also makes sense if $C$ is not necessarily a commutative algebra in $\mathcal{C}$, but only an associative algebra.  However, in this case, it might be better to
  use the following alternative definition which mimics Definition
  \ref{grodfn2} rather than Definition \ref{grodfn1}. Define an action
  of $X$ on $C$ by differential operators of order zero to be a
  composition $\phi = \mu \circ (\iota \times \Id): X \otimes C
  \rightarrow C \otimes C \rightarrow C$. For $n \geq 1$, we then
  define actions by differential operators of order $\leq n$ by the
  same inductive definition as above.  Note that, when $C$ is a
  commutative algebra, this definition is equivalent to the above.
\end{rem}
In the case $\mathcal{C} = Vect$, the category
of $\k$-vector spaces, $C$ must be a commutative
$\k$-algebra. Then, a vector space $X$ together with an action by
differential operators on $C$ means first of all that
$\phi: X \o C \rightarrow C$ is linear, i.e., we have a linear map
$X \mapsto \End_\k(C)$, $x \mapsto \phi_x$, with
$\phi_x(c) := \phi(x \otimes c)$. In terms of $\phi_x$, the two
compositions in \eqref{xacte} become $\phi_x(ab)$ and $a \phi_x(b)$,
so $\phi$ is a differential operator of order $n$ if and only if
$\phi_x(ab) - a \phi_x(b)$ is a differential operator of order $n-1$,
for all $x \in X$ and all $a, b \in C$.

For any $C$-modules $M, N$ in $\mathcal{C}$, one can similarly define a differential
operator action $\phi: X \o M \rightarrow N$ of order $\leq n$ to be
such that $[\phi,\mu]: (X \o C) \o M \rightarrow N$ is a differential
operator action of order $\leq n-1$, without changing anything else.

For any object $X$ and an algebra object $C$,
 let $\Df_{\mathcal C, \leq n}(X, C)\subset
\Hom_{\mathcal{C}}(X\o C,C)$ denote the vector space
of actions of $X$ on $C$ by differential operators of order $\leq n$.

\begin{defn}
Define the \emph{space of differential operators of order $\leq n$}
to be an object $\Df_{\leq n}(C) \in \mathcal{C}$  (if it exists)
 which represents the functor
${\mathcal{C}}\to Vect,\
X \mapsto \Df_{\mathcal{C}, \leq n}(X, C)$.
\end{defn}

By definition, whenever the object $\Df_{\leq n}(C)$
exists, it comes equipped
with an action
$\psi_{\leq n}:$ $ \Df_{\leq n} \o C \rightarrow C$,
on $C$, by differential
operators of order $\leq n$, such that the natural map $\Hom_{\mathcal{C}}(X, \Df_{\leq n})
 \rightarrow \Df_{\mathcal{C}, \leq n}(X, C)$ is an isomorphism for all $X \in \mathcal{C}$.

In the case of $\mathcal{C} = Vect$,
 we see that $\Df_{\leq n}(C)$ is  $\Df_{\leq n}(C)$ is the usual
space of differential operators of order $\leq n$ on $C$.

\begin{rem}
Similarly, one may define $\Df_{\leq n}(M, N)$, for any $C$-modules $M,
N\in \mathcal{C}$.
\end{rem}

\begin{prop}\vi If $\mathcal{C}$ is abelian and has arbitrary limits,
  then there exists $(\Df_{\leq n}(C), \psi_{\leq n})$ and it is unique
  up to unique isomorphism.

\vii There are canonical inclusions
$(\Df_{\leq n}(C), \psi_{\leq n}) \into (\Df_{\leq (n+1)}(C),
\psi_{\leq (n+1)})$.
\end{prop}

\begin{proof}[Sketch of a proof] To prove (i), take the direct sum of all actions $X \o C \rightarrow C$ by differential operators of order $\leq n$, and mod by the kernel.

To prove (ii), 
note that any action by differential operators of order $\leq n$ is also an action of order
$\leq  n+1$. This gives a morphism $j:\
\Df_{\leq n}(C)\to \Df_{\leq (n+1)}(C)$. 

Assume next that 
$X:=\ker j$ is a nonzero object. Then,
the embedding $X\into \Df_{\leq n}(C)$ gives a nonzero
action of $X$ on $C$ by differential operators of order 
$\leq  n$ which is, at the same time,
the zero
action of $X$ on $C$ by differential operators of order 
$\leq  n+1$. But  the condition for an action by differential
operators to be trivial (i.e., the action map is zero) does not depend
on what order it is considered.
We conclude that the morphism $j$ must be injective.
\end{proof}

\begin{defn}\label{Dops} Let $\mathcal{C}$ be abelian with arbitrary limits, and $C \in \mathcal{C}$ a commutative algebra object.
The differential operators $(\Df(C), \phi)$ on $C$ is
defined to be the
direct limit of all $(\Df_{\leq n}(C), \psi_{\leq n})$ (if it exists).
\end{defn}

\medskip

As an example,
let $\mathcal{C} = Rep_G$, the category of representations of a finite
group $G$, and $C$  a commutative algebra with a $G$-action, then
we find that $\Df(C)$ is the usual algebra of differential operators
on $C$, equipped with the canonical $G$-action.  In other words, as
an ordinary filtered associative algebra, $\Df(C)$ is identical with
the usual  Grothendieck construction, but it comes equipped, by definition,
with its $G$-structure.  Similarly, one may say the same thing if
$G$ is an algebraic group and $Rep_G$ is now understood as the
category of $\k[G]$-comodules.  Note that, unlike in the case of $Vect$,
in these examples it is essential that $X$ be allowed to be an arbitrary
object of $\mathcal{C}$ and not merely the unit object (since not all
differential operators are $G$-invariant).

  If $\mathcal{C}$ is the category of super vector spaces, and $C$ is
  a supercommutative algebra, then we get the usual definition of the
  superalgebra of differential operators on $C$.  Note that, in this
  case, we cannot use the original Grothendieck definition
  (taking $C$ to live in $Vect$): for example, if
  $C = \k[\theta]/(\theta^2)$ is the supercommutative algebra with
  $\theta$ odd, then differential operators according to the original
  definition include only the even operators (i.e., operators of the
  form $\lambda + \mu \theta \partial_{\theta}$ for $\lambda, \mu \in \k$), and do not include
  $\partial_\theta$ itself.
Similarly, our definitions apply in the case where 
 $\mathcal{C}$ is the category of dg-vector spaces.

Finally, given a commutative algebra $A$, let
$\mathcal{C}:= A\mmod$ be the category of left
$A$-modules, with monoidal structure given by
the tensor product over $A$. Let
$C$ be a commutative algebra object in $A\mmod$, that is
a commutative $A$-algebra. In this case,
our definition of ${\mathcal{D}}(C)$ gives
the algebra of {\em relative} differential operators
with respect to the projection $\Spec C\to\Spec A$.

\subsection{Twisted algebras}\label{twist_sec}
We begin with the definitions of twisted commutative algebras (which
date back to at least the 1950's in algebraic topology; see, e.g.,
\cite{Bar,J}).  The purpose is to provide a reasonable extension of
differential operators to \textit{tensor algebras} (which are clearly
not commutative, but are twisted-commutative).  

For the reader familiar with operads (by which we will always mean
\emph{linear} operads, i.e., operads $\mathcal{O}$ such that
$\mathcal{O}(n)$ is a vector space for all $n \geq 0$), a twisted
algebra over an operad $\mathcal{O}$ is the same as an
$\mathcal{O}$-algebra in the symmetric monoidal category of
$\SS$-modules (also known as symmetric sequences of vector spaces).
Recall that an $\SS$-module is a $\Z_{\geq 0}$-graded vector space
$\cV = \bigoplus_{m \geq 0} \cV(m)$, equipped with an action of $S_m$ on
$\cV(m)$ for all $m$.  The category of $\SS$-modules is a symmetric
monoidal category, with
$\cV \o_{\SS} \cW = \bigoplus_{p \geq 0} (\cV \o_{\SS} \cW)(p)$
where
$(\cV \o_{\SS} \cW)(p):=
\bigoplus_{m=0}^{p} \Ind_{S_m
  \times S_{p-m}}^{S_p} (\cV(m) \otimes \cW(p-m))$.
The braiding is given by
$\cQ \o_\SS \cP \iso \cP \o_\SS \cQ, \quad q \o p \mapsto (12)^{p,q} (p \o
q)$.
Moreover, any ordinary vector space may be viewed as an $\SS$-module
concentrated in degree zero, and hence we can tensor $\SS$-modules by
vector spaces, so viewed.
Thus, any
operad $\mathcal{O}$ defines both a category of twisted algebras
($\mathcal{O}$-algebras in $\SS$-modules) and a category of usual
algebras ($\mathcal{O}$-algebras in vector spaces).

We give an explicit definition in terms of permutations, which, although a
bit complicated, is important for computations.


\begin{ntn} \label{idefn} For any $i_1, \ldots, i_\ell \geq 1$,
  let $i_{i_1, \ldots, i_\ell}: S_\ell \into S_{i_1+\ldots+i_\ell}$ be
  the monomorphism which sends a permutation $\sigma \in S_\ell$ to
  the permutation of $S_{i_1+\ldots+i_\ell}$ that permutes the cells of
  the partition
  $(\{1, \ldots,i_1\},\{i_1+1, \ldots,i_1+i_2\}, \ldots, \{i_1+\ldots+i_{\ell-1}+1,
  \ldots, i_1+\ldots+i_\ell\})$,
  preserving the order in each cell of the partition, by 
  rearranging all the cells according to $\sigma$.

We will frequently use the shorthand
\begin{equation}
\sigma^{i_1, \ldots,i_\ell} := i_{i_1, \ldots,i_\ell}(\sigma).
\end{equation}
In particular, $(12)^{m,n}$ is the permutation
$(1,2, \ldots,m,m+1, \ldots,m+n) \mapsto (m+1, \ldots,m+n,1, \ldots,n)$.
\end{ntn}

\begin{defn} \label{tcad} A twisted associative algebra $\cA := \bigoplus_{m \geq 0} \cA(m)$ 
is:
\begin{enumerate}
\item a graded associative algebra over $\k$ with multiplication $\mu: \cA \o \cA \rightarrow \cA$, 
      together with
\item an action of $S_m$ on $\cA(m)$ for all $m \geq 1$, such that
\item $\mu: \cA(m) \o \cA(n) \rightarrow \cA(m+n)$ is a map of $S_m \times S_n \subset S_{m+n}$-modules.
\end{enumerate}

A twisted \textbf{commutative} algebra is a twisted associative algebra
$A$ such that
\begin{equation} \label{twcomm}
(12)^{m,n} \mu(a \o b) = \mu(b \o a), \quad \forall a \in \cA(m), b \in \cA(n).
\end{equation}\end{defn}

In other words, for a twisted commutative algebra we have the commutative diagram:
$$
\xymatrix{
\cA(m) \o \cA(n) \ar[rd]^{\mu^{\op}} \ar[r]^\mu & \cA(m+n) \ar@{->}[d]^{(12)^{m,n}}_\sim
\\
& \cA(m+n).
}
$$

\begin{rem}
In general, a twisted algebra $\cA$ over
an operad $\mathcal{O}$  may be explicitly defined as follows:
$
\cA = \bigoplus_{m \geq 0}\ \cA(m),
$
with an $S_m$ action on each $\cA(m)$ (i.e., $A$ is an $\SS$-module), 
equipped with a map 
\begin{equation} \label{Oaact}
\bigoplus\nolimits_{m \geq 0}\ \mathcal{O}(m) \o_{S_m} (\cA^{\o m}) \rightarrow \cA,
\end{equation}
which descends from maps
\begin{equation}
\mathcal{O}(m) \o \cA(i_1) \o \cdots \o \cA(i_m) \rightarrow \cA(i_1+\cdots+i_m)
\end{equation}
of
$S_{i_1} \times \cdots \times S_{i_m} \subset S_{i_1+\cdots +i_m}$-modules
(under the standard embedding).
Also, the map \eqref{Oaact} is required to satisfy an associativity
condition.  
\end{rem}

  Any ordinary commutative algebra is
an example of twisted commutative algebra
  concentrated in degree zero.

Another example of a  twisted commutative algebra
is the graded space $\bigoplus_{m \geq 0} \k[S_m]$ where 
multiplication is given by tensor product
(using the standard inclusion
$(S_i \times S_j) \subset S_{i+j}$) and where the $S_n$-action 
is given by \emph{conjugation}.


For any $\SS$-module $\cM$, let $\Syns_{\SS} \cM$ denote the symmetric
algebra in the category of $\SS$-modules generated by $\cM$.
Explicitly, we have
\begin{equation}
\Syns_{\SS} \cM = \bigoplus\nolimits_{m \geq 0}\  (\cM^{\otimes_\SS m})_{S_m}.
\end{equation}
This is a free  twisted commutative algebra.

\begin{ntn}
Given  a vector space $V$ and $i \in\{0,1\}$, let $V_{\SS,i}$ denote the $\SS$-module concentrated in degree $i$ with $(V_{\SS,i})_i = V$.
\end{ntn}

It is clear that 
$\Syns_{\SS} V_{\SS,0}=\Syns_\k V$ is the ordinary symmetric
algebra of the vector space $V$.
Similarly, one has
$\Syns_{\SS} V_{\SS,1}=T_\k V$
is the ordinary tensor algebra of the vector space $V$.
Here,  $T_\k V$ is equipped with the standard grading
and the  $S_n$-action on $V^{\o n}$
is given by permutation of components.

\subsection{Differential operators on twisted-commutative algebras}
 Let $\mathcal{C} = \SS\mmod$ be the monoidal
category of $\SS$-modules. Then we
  may rewrite our general definition
of differential operators in terms of individual
  operators as follows.  

A morphism $X \o C \rightarrow C$ is an
  action of $X$ on $C$ by differential operators of order $\leq n$
  iff, for all homogeneous $x \in X$, the map
  $\ds C \mathop{\cong}^{vs.} \{x\} \o C \rightarrow C$ is a
  \emph{differential operator $C \rightarrow C$ of degree $|x|$ and
    order $\leq n$},
  where the latter is defined as follows.  First, an \emph{operator
    $C \rightarrow C$ of degree $m \geq 0$}
  is a degree-$m$ morphism of graded vector spaces such that the
  restriction to $C_k \rightarrow C_{m+k}$ is a morphism of
  $S_k$-modules where the action on $C_{m+k}$ is via
  $S_k \cong (\{\Id_{S_m}\} \times S_k) \into S_{m+k}$, i.e., the
  composition has the property that
  $\sigma \mapsto \sigma', \sigma'(i) = \begin{cases} i, & \text{if }
    i \leq m, \\ m+\sigma(i-m), & \text{otherwise}.\end{cases}$

Next, a differential operator of degree $m$ and order $\leq -1$ is the
zero morphism. Inductively, $D: C \rightarrow C$ is a differential
operator of degree $m$ and order $\leq n$ iff, for all homogeneous
$c \in C$, $[D,c]: C \rightarrow C$ is a differential operator of
degree $m+|c|$ and order $\leq n-1$.  Here, $[D,c]$ sends any
homogeneous $a \in C$ to $D(ca) - (21)^{|c|,m,|a|} cD(a)$, and extends
linearly.

The above theory works better than Grothendieck's original definition on
tensor algebras.  As opposed to Remark \ref{groth}, we have the
following classification of differential operators in the
\textit{twisted-commutative} sense on $T_\k V$:
\begin{prop} \label{diffopstkv}
The twisted associative 
algebra $\Df(A)$ of differential operators on $A = T_\k V$ is the
quotient of the free twisted associative algebra $T_{\SS} (\End(V)_{\SS,0} \oplus V_{\SS,1})$
by the commutation relations (where $[a,b] := a b - \sigma^{|a|,|b|} b a$ is the
\textbf{twisted} commutator): for all $v, w \in V$ and $\phi, \psi \in \End(V)$,
\begin{equation} \label{dfarelns}
[\lambda_v, \lambda_w] = 0, \quad [D_\phi, \lambda_v] = \lambda_{\phi(v)}, \quad [D_\phi, D_\psi] = D_{[\phi, \psi]},
\end{equation}
where, for each $\phi \in \End(V)$, we denote by $D_\phi$ the
corresponding image in $T_{\SS} (\End(V)_{\SS,0} \oplus V_{\SS,1})$, which
acts as a derivation on $V$ by applying $\phi$, and we denote by $\lambda_v$, 
for each $v \in V$, the operator of left-multiplication by $v$, which is the
summand of $V_{\SS,1}$ in $T_{\SS} (\End(V)_{\SS,0} \oplus V_{\SS,1})$.
\end{prop}

In other words, the proposition may be interpreted as saying that
$\Df(A)$ is the universal enveloping twisted associative
algebra
of the twisted Lie algebra $\End(V)_{\SS,0} \oplus V_{\SS,1}$ with bracket
\eqref{dfarelns}.

\begin{proof}
  Inductively, a differential operator of order $\leq n$ on $T_\k V$
  is determined by its restriction to $V^{\o \leq n}$. By compatibility
  with permutations, this restriction must be describable using
  polynomials of degree $\leq n$ in operators of the form $D_\phi$,
  together with tensoring everything on the left by linear
  combinations of elements $a_1 \o \cdots \o a_\ell$.  Such an element
  is easily verified to be a differential operator.  From this, the
  above description easily follows.
\end{proof}



\subsection{Differential operators are almost-commutative}
By \emph{almost-commutative}, we will always mean filtered associative with
a commutative associated graded.

Let $\mathcal{C}$ be any abelian symmetric monoidal category with arbitrary
limits. Parallel
to the case of vector spaces, we have
\begin{thm} \label{dot} The differential operators $\Df(C)$ have the
  canonical structure of filtered algebra in the category
  $\mathcal{C}$. Moreover, the associated graded $\gr \Df(C)$ is a
  commutative algebra in $\mathcal{C}$. The commutator on
  $\Df(C)$ induces a Poisson algebra structure in $\mathcal{C}$ on
  $\gr \Df(C)$.
\end{thm}
Here, a \emph{filtered algebra object in $\mathcal{C}$} is an object $A \in \mathcal{C}$ 
equipped with a filtration in $\mathcal{C}$, $0 \subset A_{\leq 0} \subset A_{\leq 1} \subset \cdots$, such that $A_{\leq i} \cdot A_{\leq j} \subset A_{\leq i+j}$ for all $i,j$.  The associated graded of an almost-commutative algebra in
$\mathcal{C}$ is \emph{always} Poisson (via the commutator), so the final assertion of the theorem
is tautological.
\begin{proof}[Sketch of a proof]
By definition, $\Df(C)$ is a filtered object of $\mathcal{C}$. We have
to exhibit a multiplication $\circ: \Df(C) \o \Df(C) \rightarrow \Df(C)$ in $\mathcal{C}$, show that it respects the filtration, and that the commutator $D \circ D' - D' \circ D$ of operators of order $\leq n, n'$ has order
$\leq n+n'-1$.

These can be proved using the inductive definition of differential
operators, similarly to the case where $\mathcal{C} = Vect$.  Let
$D: X \o C \rightarrow C$ and $D': X' \o C \rightarrow C$ be
differential operators.  Then, one has the following
\begin{lemma}
Let $(\ad\ \mu)(f) := [f, \mu]$ for any operator (see Definition \ref{admudefn}). Then, as operators $(X \o X' \o C^{\o p}) \o C \rightarrow C$, we have
\begin{equation}
(\ad\ \mu)^p(D \circ D') = \sum_{m+n=p} (\ad\ \mu)^m(D)  \circ (\Id \o (\ad\ \mu)^n(D')) \circ \beta^{(m)},
\end{equation}
where
\begin{equation}
\beta^{(m)}(X \o X' \o C^{\o p} \o C) = (X \o C^{\o m}) \o (X' \o
C^{\o n}) \o C
\end{equation}
is the composition of braidings which preserves the left-to-right
order of the factors of $C$.
\end{lemma}
The proof of the lemma is the same as in the case $\mathcal{C} = Vect$, where
it becomes  (again using the right adjoint action $(\ad\ a) (D) = [D, a]$):
\begin{equation}
[ \cdots [D \circ D', a_1],a_2],\cdots,a_{p}] = 
\sum_{I \subset \{1, \ldots,p\}} \bigl( \prod_{i \in I} \ad\  a_i\bigr)(D) \circ \bigr(\prod_{j \in \{1, \ldots,p\}\setminus I} \ad\  a_j \bigr)(D').
\end{equation}
We omit the proof.

The lemma shows that the filtration is multiplicative (we get a
filtered algebra).  To show almost-commutativity, we only need to
notice that two operators of order zero must commute with each other,
since $C$ is a commutative algebra in $\mathcal{C}$.
\end{proof}

\subsection{Twisted Poisson algebras}
The notion of twisted Poisson algebra is closely related to Van den
Bergh's double Poisson algebras, as will be explained in Section
\ref{dwps} (cf.~\cite{Stw}).  As before, a twisted Poisson algebra is
the same as an algebra, in the category of $\SS$-modules, over the
Poisson operad. We give, however, a more explicit definition:
\begin{defn} \label{twpd}
  A twisted (commutative)\footnote{For us, Poisson algebras are commutative
  by definition.}
  Poisson algebra is a twisted commutative algebra
  $A = \bigoplus_{m \geq 0} A[m]$ with a graded bracket
  $\{\, ,\}: A \o A \rightarrow A$ satisfying
\begin{enumerate}
\item $\{\, ,\}: A[m] \o A[n] \rightarrow A[m+n]$ is a map of $S_m \times S_n \subset S_{m+n}$-modules,
\item $\{\, ,\}$ is a twisted Lie bracket: For $a \in A[m], b \in A[n], c \in A[p]$,
\begin{gather} \label{twskewsym}
\{a,b\} = -(12)^{n,m}\{b,a\}, \\ \label{twjac}
\{\{a,b\},c\} + (132)^{p,m,n} \{\{c,a\},b\} + (123)^{n,p,m} \{\{b,c\},a\} = 0;
\end{gather}
\item $\{\, ,\}$ is twisted Poisson: this means
\begin{equation}
\{a,bc\} = \{a,b\}c + (12)^{n,m,p} b\{a,c\}. \label{twleib}
\end{equation}
\end{enumerate}
\end{defn}
As a consequence of Proposition \ref{dot}, we have
\begin{cor}
  Let $A$ be a filtered twisted associative algebra with twisted
  commutative $\gr A$. Then $\gr(A)$ has a natural twisted Poisson algebra
  structure.
\end{cor}



\begin{prop} For $A = T_\k V$, $\gr \Df(A) \cong \Syns_{\SS-mod}(\End(V)_{\SS,0} \oplus V_{\SS,1})$, with the twisted Poisson bracket given by
\begin{equation}
\{\lambda_v, \lambda_w\} = \{D_\phi, D_\psi\} = 0, \quad \{D_\phi, \lambda_v\} = \phi(v).
\end{equation} 
\end{prop}
The proof is straightforward.

\section{Differential operators on associative 
algebras and wheelgebras}
\label{wps} The main goal of this section is to extend the formalism
of differential operators from $T_\k M$ (for which the previous
section is applicable) to $T_A M$, where
$A$ is an associative algebra and $M$ is an $A$-bimodule.
  This requires introducing  a new algebraic
structure called a \textit{wheelgebra}, which will be fundamental to
the differential geometry of associative algebras. Precisely, we will
define, for any associative algebra $A$, a \emph{commutative wheelgebra}
$\Ff(A)$, and define differential operators on $A$ in terms of categorical
differential operators on $\Ff(A)$.

\subsection{Wheelspaces}\label{wheel_sec} Given an associative
algebra $A$, one can use multiplication in $A$ to
construct, for each $m\geq 1$, various \emph{contraction}
maps $A^{\o m} \rightarrow A^{\o (m-1)}$.
Abstract  properties of such  contraction maps
are conveniently formalized in the following

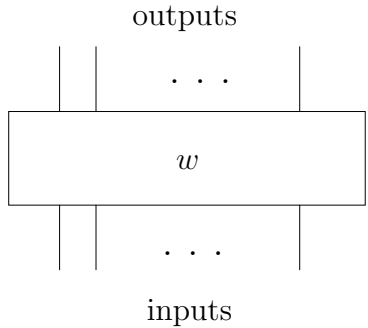
\begin{figure}[hbt]
\begin{center}
\input{chip.eepic}
\caption{Diagrammatic representation of $w \in \cW(m)$}
\label{chipfig}
\end{center}
\end{figure}

\begin{figure}[hbt]
\begin{center}
\input{chipmuij.eepic}
\caption{Diagrammatic representation of $\mu_{i,j} w \in \cW(m-1)$}
\label{chipmuijfig}
\end{center}
\end{figure}
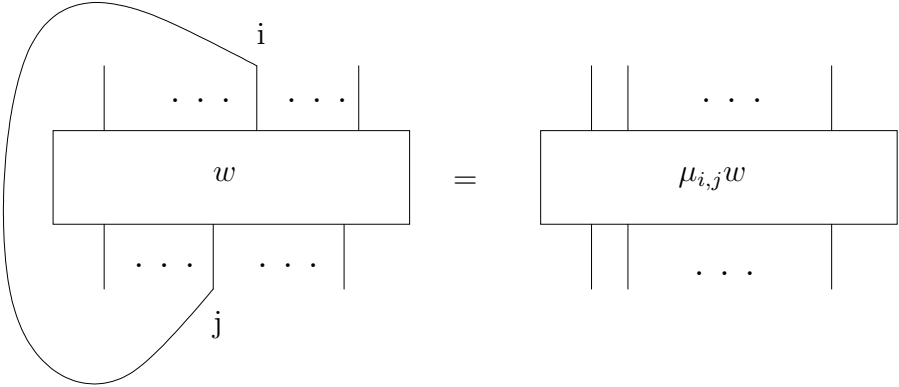

\begin{figure}[hbt]
\begin{center}
\input{chipperm.eepic}
\caption{Diagrammatic representation of $((1,2) \times (m-1,m)) w \in \cW(m)$}
\label{chippermfig}
\end{center}
\end{figure}
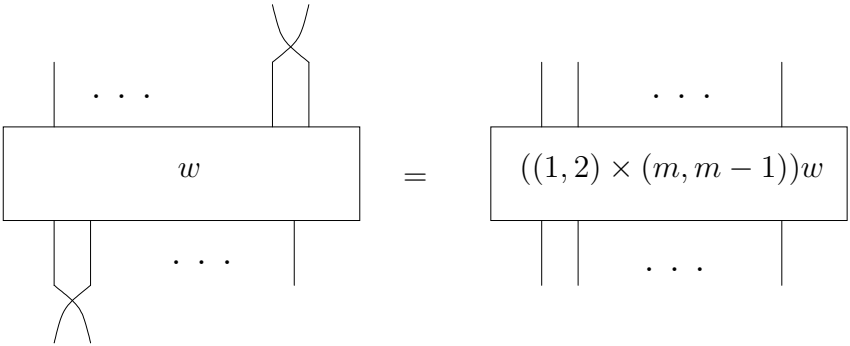

\begin{defn}
A \textit{wheelspace} is 
a $\Z_{\geq 0}$-graded vector space $\cW = \bigoplus \cW(m)$, with $S_m \times S_m$-actions on $\cW(m)$, and contraction maps
\begin{equation}
\mu_{i,j}: \cW(m) \rightarrow \cW(m-1)
\end{equation}
which are morphisms of $S_{m-1} \times S_{m-1}$-modules, where the action
on the left is by permuting the sets $(\{1, \ldots,m\} \setminus \{j\}) \times (\{1, \ldots,m\} \setminus \{i\})$, satisfying the following condition:
\begin{equation} \label{gperm}
\mu_{i,j} = \mu_{1,1} \circ ((1,2, \ldots,j) \times (1,2, \ldots,i)).
\end{equation}
Also, the contraction
maps must satisfy the associativity condition
\begin{equation}
\mu_{i,j} \circ \mu_{k,\ell} = \mu_{k', \ell'} \circ \mu_{i', j'},
\end{equation}
where, if $i < k$, then $(i', k') = (i, k-1)$, and otherwise $(i', k')
= (i+1, k)$; the same relationship holds between $(j, \ell)$ and
$(j', \ell')$.
\end{defn}
Informally, we think of elements $w \in \cW(m)$ as represented
by a picture as in Figure \ref{chipfig}, with a total of $m$ inputs
and $m$ outputs.
The contraction operation $\mu_{i,j}$ is diagrammatically represented
by joining the $i$-th output to the $j$-th input, as in Figure \ref{chipmuijfig}.
The permutation action is diagrammatically represented by permuting
inputs and outputs.  For instance, $((1,2) \times (m-1,m)) w$ is
represented in  Figure \ref{chippermfig}.
Then, the compatibility and associativity conditions for the permutations
and contractions say that any combination of operations applied to $w$ depends
only on the resulting diagram, up to equivalence (where equivalence means up to homotopy of the legs preserving the left-to-right order of their top endpoints, as well as of their bottom endpoints).

Any wheelspace is, in particular, an $\SS$-module, with diagonal
permutation action ($\sigma$ acts on $w$ by
$(\sigma \times \sigma)(w)$).  Then, the category of wheelspaces forms
a symmetric monoidal category in a canonical way such that the
forgetful functor $\text{wheelspaces} \rightarrow \SS\mmod$ is
symmetric monoidal, and
such that, for $w \in \cW(m), v \in \cW(n)$, and any
$i,j \in \{1,\ldots, m\}$ and $k,l \in \{1,\ldots, n\}$,
\begin{equation}
\mu_{i,j} (w \o v) = (\mu_{i,j} w) \o v, \quad \mu_{k+m,\ell+m}(w \o v) = w \o (\mu_{k,\ell} v).
\end{equation}
Slightly informally, the definition of wheelgebras is then the following:
\begin{defn} \label{whgdfn} A (commutative, associative, Lie, etc.)
  wheelgebra is a wheelspace which is also a twisted (commutative,
  associative, Lie, etc.)  algebra, such that the algebraic operations
  (multiplication, bracket, etc.) are morphisms of wheelspaces, i.e.,
  the operations $\cW(m_1) \otimes \cdots \otimes \cW(m_n) \rightarrow \cW(m_1 + \cdots + m_n)$ are morphisms of $S_{m_1 + \cdots + m_n} \times S_{m_1 + \cdots + m_n}$-modules which commute with contractions $\mu_{i,j}, i,j \in \{1,2,\ldots, (m_1+\cdots+m_n)\}$.
\end{defn}
For example, an associative wheelgebra $\cW$ is a wheelspace equipped
with an associative multiplication
$\cW(m) \otimes \cW(n) \rightarrow \cW(m+n)$ which is a morphism of
$S_{m+n} \times S_{m+n}$-modules compatible with the contraction.

The word \emph{twisted} above is necessary above since a commutative
wheelgebra is not an ordinary commutative algebra, only a twisted
commutative algebra (for the associative case, one in fact has an
ordinary associative algebra because the braiding does not enter into
the associativity axiom).  Recall \eqref{twcomm} for the explicit
definition of twisted commutative algebra. 

A slightly different way to state the above definition is as follows:
A commutative wheelgebra is a twisted-commutative algebra together
with an extension of the $S_m$ action in degree $m$ to an
$S_m \times S_m$-action such that the original action is obtained by
the diagonal embedding $S_m \into S_m \times S_m$, and together with
contractions $\cW(m) \rightarrow \cW(m-1)$ which are compatible with
multiplication operations. For any other type of algebra, replace all
instances of ``commutative'' here by some other type of algebra.

For the reader familiar with operads, we can state the above
definition more precisely as follows:
\begin{defn}\label{whopgdfn}
Given any operad $\mathcal{O}$, an $\mathcal{O}$-wheelgebra is an algebra over $\mathcal{O}$ in the symmetric monoidal category of wheelspaces.  
\end{defn}
As always, by an operad we mean a usual $\k$-linear operad (not an operad
in the category of wheelspaces).  The above definition then makes
sense because the category of wheelspaces is ``tensored over vector
spaces'': this means that the tensor product $V \o W$ of a vector space
$V$ and a wheelspace $\cW$ is well defined.  Precisely, this is given by
the symmetric monoidal functor
$\text{Vect} \rightarrow \text{Wheelspaces}$ sending a vector space to
the corresponding wheelspace concentrated in degree zero.  Then, a
$\mathcal{O}$-wheelgebra structure on a wheelspace $\cW$ is a map of
wheelspaces
\begin{equation}\label{owheqn}
  \bigoplus_{m \geq 0} \mathcal{O}(m) \o_{S_m} (W^{\o_{\text{wh}}m} ) \rightarrow W,
\end{equation}
which satisfies the same associativity and unit constraints as in the
usual setting. The notation $\otimes_{\text{wh}}$ is to remind the
reader that this is the tensor product in the category of wheelspaces.

\begin{rem} \label{wheelprop} A commutative wheelgebra is a special
  case of a (nonunital) wheeled PROP \cite{MMS}.  Namely, a
  commutative wheelgebra is a nonunital wheeled PROP satisfying the
  condition that $\Hom(m,n) = 0$ for any $m \neq n$.  In the language
  of commutative wheelgebras, a wheeled PROP unit is an element $t$ in
  degree one such that $\mu_{1,i}(t x) = x = \mu_{i,1}(t x)$ for all
  $i \neq 1$.  Furthermore, one may define a bigraded version of
  wheelspaces where, in bidegree $(m,n)$, one has an action of
  $S_m \times S_n$. This also forms a symmetric monoidal category, and
  in this category, commutative algebras are precisely nonunital
  wheeled PROPs.  Then, a unit is as before, but now placed in
  bidegree $(1,1)$.
\end{rem}

\subsection{The commutative wheelgebra $\Ff(A)$}
We now describe our most important example of a wheelgebra, $\Ff(A)$.
We emphasize that the definition is partly motivated by the
application given in, e.g., \S \ref{ovwsec}, and the reader is
encouraged to glance there to help understand the reason for it.  In
one sentence, $\Ff(A)$ is the commutative wheelgebra freely generated
by $A$ in degree $1$, subject to the relation that the contraction
operation $A \o A \rightarrow A$ is the multiplication on
$A$. Explicitly, as a twisted-commutative algebra, it is given by


\begin{equation}  \label{ffadfn}
\Ff(A) := \bigoplus_{m \geq 0} \Ff_m(A), \quad \Ff_m(A) := (\Ind_{S_m}^{S_m \times S_m} A^{\o m}) \o \Syns_\k A_\cyc,
\end{equation}
where the $\Ind$ uses the \textbf{diagonal embedding} $S_m \into (S_m \times S_m)$, i.e.,
$\sigma \mapsto (\sigma \times \sigma)$. \emph{Note that we will always use the diagonal embedding in expressions $\Ind_{S_m}^{S_m \times S_m}$ from now on, without mentioning it.}  For convenience, let us write
the $S_m \times S_m$-action on $\Ff_m(A)$ by having the first factor of
$S_m$ act on the left and the second act on the right. 

In the pictorial representation of Figure \ref{chipfig}, the $i$-th
input corresponds to the tensor factor which is changed by the $i$-th
left $A$-module structure, as follows: Given the element
$\sigma_L (f_1 \o \cdots \o f_m) \sigma_R$, plugging $a \in A$ into
the $i$-th input of the diagram should be thought of as replacing
$f_{\sigma_L^{-1}(i)}$ with $a f_{\sigma_L^{-1}(i)}$, and similarly
plugging $a$ into the $i$-th output should be thought of as replacing
$f_{\sigma_R^{-1}(i)}$ with $a f_{\sigma_R^{-1}(i)}$. This
interpretation translates the pictorial identities (Figures
\ref{chipmuijfig} and \ref{chippermfig}) into the operations on
$\Ff(A)$ as defined in \eqref{ffadfn}.  We explain these operations
precisely as follows.

The $S_m \times S_m$-module structure in degree $m$ (i.e., on
$\Ff_m(A)$) is immediate from \eqref{ffadfn}.

The commutative multiplication (let us call it $\cdot$ to distinguish it from
$\o_\k$) is part of the aforementioned twisted-commutative algebra structure,
and is given by
\begin{equation}
((\sigma_L X \sigma_R)\o Y) \cdot ((\sigma_L' X' \sigma_R') \o Y') =
[(\sigma_L \times \sigma_L') (X \o X') (\sigma_R \times \sigma_R')] \o (Y \& Y'),
\end{equation}
for $X \in A^{\o m}, X' \in A^{\o n}$ and any
$Y, Y' \in \Syns(A_{\cyc})$.  

The contraction operations $\mu_{i,j}$ are obtained by multiplying the $i$-th
output to the $j$-th input.  This means that
\begin{multline}
\mu_{i,j} (a_1 \o \cdots \o a_n) \\= \begin{cases} \sigma_L (a_{i} a_j \o a_1 \o \cdots \o \hat a_{i} \o \cdots \o \hat a_{j} \o \cdots \o a_n) \sigma_R, & \text{if $i \neq j$}, \\
(a_1 \o \cdots \o \hat a_i \o \cdots \o a_n) \o [a_i]_\cyc, & \text{if $i = j$},
\end{cases}
\end{multline}
where, in the first equation, $\sigma_L$ and $\sigma_R$ are the
permutations such that the order of the inputs and outputs have the
$a_k$ in increasing-$k$ order, except that the term $a_i a_j$ is
considered as having $i$ as an input, and as $j$ as an output.  The
$[a_i]_\cyc$ in the second line means the image of $a_i$ in $A_\cyc$.
Finally, the hats ($\hat a_i$ and $\hat a_j$) mean that the given terms
($a_i$ and $a_j$) are \textbf{omitted} from the tensor product.

The above formulas determine all of the contraction action by
compatibility with permutations \eqref{gperm}, linearity, and the
condition $\mu_{i,j} (X \o Y) = \mu_{i,j}(X) \o Y$ if $Y \in \Sym
A_\cyc$.

\begin{rem} The commutative wheelgebra $\Ff(A)$ is not merely
a nonunital wheeled PROP (cf.~Remark \ref{wheelprop}), but in fact is
unital, with unit $1 \in A \subset \Ff_1(A)$.
\end{rem}
\begin{rem} In \cite{MMS}, the notion of algebra over a wheeled PROP
  was defined (which generalizes the notion of algebra over a PROP,
  and of algebra over an operad). Since, as in Remark \ref{wheelprop},
  a commutative wheelgebra is a special kind of wheeled PROP, one
  obtains a notion of algebra over a commutative wheelgebra. In the
  case
  $\cW = \Ff(A)$, finite-dimensional algebras $S$ over $\cW$ are the same as
  \textit{representations of $A$} with underlying 
  vector space $S$ (since the compatibility with contractions says
  that $\rho([a]_\cyc) = \tr \rho(a)$).  If we generalize $S$ from a
  finite-dimensional vector space to a finitely-generated free module
  over a commutative $\k$-algebra $B$, we get families of representations
  parametrized by $B$.
  This theory (and its application to polyvector fields, differential
  operators, BV structures, and so on) will be discussed in
  \S \ref{repsec}.
\end{rem}

\subsection{Differential operators on commutative wheelgebras}
We may now invoke the categorical definition of differential operators:
\begin{defn}
  For any commutative wheelgebra $\cW$, let $\Df(W)$ denote the almost-commutative
  wheelgebra of differential operators on $\cW$.
\end{defn}
As before, the associated graded of an almost-commutative wheelgebra
is a Poisson wheelgebra.

Explicitly, a differential operator of twisted-degree $i$ on a
commutative wheelgebra $\cW$ is a differential operator $D$ of
twisted-degree $i$ on the underlying twisted-commutative algebra of
$\cW$, which satisfies the following condition:
\begin{equation}
D \circ \mu_{j,k} = \mu_{i+j,i+k} \circ D.
\end{equation}
The contractions of $\Df(W)$ are given by the equation
\begin{equation}
\mu_{i,j}(D) = \mu_{i,j} \circ D.
\end{equation}

\begin{exam} For any associative algebra $A$, one may consider
the wheelgebra $\Df(\Ff(A))$, which we will sometimes call the \emph{wheeled differential operators} on $A$.
\end{exam}

\begin{exam} The differential operators $\Df(\Ff(T_\k V))$ include, as
  a subalgebra, the operators $\Df(T_\k V)$: the latter are exactly
  the differential operators which preserve the \textbf{non-wheeled}
  subalgebra of $\Ff(T_\k V)$ generated by $V$ (this subalgebra is
  isomorphic as a twisted algebra to $T_\k V$).  Note also that all
  elements of $\Df(\Ff(T_\k V))$ are determined by their restriction
  to the aforementioned subalgebra.
\end{exam}

\subsection{Lie and Poisson wheelgebras}\label{lpwgbas}
A Lie or Poisson wheelgebra is a Lie or Poisson algebra in the category
of wheelspaces, respectively.  

Explicitly, a Lie wheelgebra bracket on
a wheelspace $\cW$ is a twisted Lie bracket
  $\{-,-\}: W \o W \rightarrow W$  satisfying
\begin{equation} \label{poiscontr}
\{\mu_{i,j} a, b\} = \mu_{i,j} \{a, b\}, \quad 1 \leq i,j \leq |a|.
\end{equation}
A Poisson wheelgebra is a commutative wheelgebra equipped with
a Lie wheelgebra bracket, which satisfies the Leibniz rule.
The Leibniz rule may equivalently be stated as saying that the
underlying twisted structures form a twisted Poisson algebra.

We remark that there is a construction of universal enveloping
wheelgebra of Lie wheelgebras and a PBW theorem for Lie wheelgebras,
which is completely analogous to the usual case (for the twisted,
rather than wheeled, setting, this has been known for a long time in
topology; see, e.g., \cite{Bar}).  

Similarly, one may construct a Koszul complex for a Lie wheelgebra
(giving a resolution of the ``augmentation'' module of the universal
enveloping wheelgebra).  Specifically, the analogue of
$\Sym \mathfrak{g}$ is $\Syns_{\text{wh}}(L)$, where for an arbitrary
wheelspace $L$, $\Syns_{\text{wh}}(L)$ is defined to be the free
commutative wheelgebra generated by $L$ (this means, the commutative
wheelgebra satisfying the universal property that
$L \into \Syns_{\text{wh}}(L)$ is an embedding of wheelspaces; this
can also be presented explicitly with some work).  Similarly, the
analogue of $\La \mathfrak{g}$ is the \textit{super} version
$\SuperSym_{\text{wh}}(L(1))$, where the braiding is by superbraidings
$\tau_\sigma$ with respect to the grading $|L(1)|=1$, and uses
and one uses supercommutators rather than commutators. The result is a
supercommutative wheelgebra.  Then, if $L$ is a Lie wheelgebra,
$\SuperSym_{\text{wh}}(L \oplus L(1))$ (with $|L|=0, |L(1)|=1$) may be
equipped with a Koszul differential, summing with sign over applying
the bracket, for $f \in \Syns_{\text{wh}}(L)$ and
$x_1, \ldots, x_m \in L(1)$:
\begin{multline}
d\bigl(f \otimes \bigl(x_1 \o \cdots \o x_m\bigr) \bigr) \\ = \sum_{i=1}^m (-1)^{i+1} (\Id_{S_{|f|}} \times \sigma_i) \bigl(f \cdot x_i \otimes \bigl(x_1 \o \cdots \hat x_i \cdots \o x_m\bigr)\bigr)(\Id_{S_{|f|}} \times \sigma_i) \\ - \sum_{i < j} (-1)^{i+j} (\Id_{S_{|f|}} \times \sigma_{ij}) \bigl(f \o \bigl(\{x_i, x_j\} \o x_1 \o \cdots \hat x_i  \cdots \hat x_j \cdots \o x_m\bigr)\bigr)(\Id_{S_{|f|}} \times \sigma_{ij}),
\end{multline}
where $\sigma_i$ is the inverse of the permutation sending
$x_1 \o\cdots \o x_m$ to $x_i \o x_1 \o \cdots \hat x_i \cdots \o x_m$
and $\sigma_{ij}$ is the inverse of the permutation that sends
$x_1 \o \cdots \o x_m$ to
$x_i \o x_j \o x_1 \o \cdots \hat x_i \cdots \hat x_j \cdots \o x_m$.

Details will be discussed elsewhere.

\subsection{Wheeled Poisson structure on the double cotangent bundle}
\label{dwps}
We briefly recall the {\em even} version of Van den Bergh's double
Schouten-Nijenhuis bracket \cite{VdB} (we will explain things slightly
differently from \cite{VdB}).

An initial idea is to try to produce an
analogue for $\DDA$ of the standard Lie bracket of vector fields.
Naively, this involves a formula of the form
$[\xi, \eta] = \xi \circ \eta - \eta \circ \xi$. (One certainly wants
a commutator-type formula if this is to come from the associated
graded of a space of differential operators.)  Of course, as written,
this cannot work, since $\eta(A) \subset A \o A$, and so after
applying $\eta$, we cannot apply $\xi$: we have to choose either
$\Id \o \xi$ or $\xi \o \Id$.  So, there are four ways to apply one of
$\xi, \eta$ and then the other.  Our goal is to obtain another
derivation. It turns out that there is a two-dimensional vector space
of linear combinations of such compositions that always yield derivations
$A \rightarrow (A \o A \o A)$, using the \textit{outer} bimodule
action on $A \o A \o A$. This vector space is spanned by
\begin{gather} \label{vdbelts}
(\xi \o \Id) \circ \eta - (\Id \o \eta) \circ \xi, \quad  (\Id \o \xi) \circ \eta - (\eta \o \Id) \circ \xi.
\end{gather}
These two elements form a \textit{canonical} basis, up to scaling both
by the same nonzero constant, as
will be explained after the following definitions.
\begin{defn}[\cite{VdB}] Define the brackets 
\begin{gather} \label{vdblsub}
\ldb -, -\rdb_L: \DDA \o \DDA \rightarrow \DDA \o A, \\
\ldb -, -\rdb_R: \DDA \o \DDA \rightarrow A \o \DDA \label{vdbrsub}
\end{gather}
to be the unique elements such that the associated derivations $A \rightarrow A \o A \o A$ (viewing $A \o A \o A$ as an $A$-bimodule using the \textbf{outer} action) are
\begin{gather}
\ldb \xi, \eta \rdb_L^\sim = (\xi \o \Id) \circ \eta - (\Id \o \eta) \circ \xi,
\ldb \xi, \eta \rdb_R^\sim = (\Id \o \xi) \circ \eta - (\eta \o \Id) \circ \xi.
\end{gather}
\end{defn}
\begin{defn}[\cite{VdB}] A double Poisson algebra is an associative algebra $A$ together with a $\k$-bilinear bracket $\ldb -,- \rdb: A \o A \rightarrow A \o A$ 
satisfying the conditions
\begin{gather}
\ldb \xi, \eta_1 \eta_2 \rdb = \ldb\xi, \eta_1\rdb \eta_2 +
\eta_1 \ldb\xi, \eta_2\rdb, \label{dleib} \\
\ldb-,-\rdb = - \tau_{(21)} \circ \ldb-,-\rdb \circ  \tau_{(21)}, \label{dskew}\\
\ldb\xi, a \eta\rdb = a \ldb\xi, \eta\rdb + \ldb\xi, a\rdb \eta, \label{dpoiss} \\
\sum_{i=0}^{2} \tau_{(123)^i} \circ \bigl(\ldb-,-\rdb^{23} \circ \ldb-,-\rdb^{12}\bigr) \circ \tau_{(123)^{-i}}= 0. \label{djac}
\end{gather}
\end{defn}
\begin{rem} \label{dbpref} Note that the resulting map
  $A \o A \rightarrow A \o A$ is a derivation from the first component
  to the \textit{inner} action on the target, and a derivation from
  the second component to the \textit{outer} action on the target.  In
  some derivation sense, this is saying that
  \textit{left-multiplication components are reversed}: multiplying on
  the left of the first component on the domain $A \o A$ is sent to
  multiplying the left of the second component of $A \o A$ in the
  image, up to a derivation term. More precisely, for such a derivation $D$,
\begin{equation}
D(a b_1 \o b_2) = (1 \o a) D(b_1 \o b_2) + D(a \o b_2) (b_1 \o 1),
\end{equation}
and saying that ``left-multiplication components are reversed'' is
saying that, if we compare the LHS with the first term on the RHS, the
$a$ changes from the first to the second component.  On the other
hand, \textit{right-multiplication components are preserved}, in the
sense that, comparing the LHS with the second term on the RHS, the
$b_1$ remains in the first component.  We interpret this as saying
that \textit{the double bracket prefers right-multiplication}.  We
will eliminate this preference when we introduce the wheeled Poisson
version of the above.  \erem

Following Van den Bergh, we  introduce a double bracket, which 
is an {\em even} version of the {\em odd} double  bracket,
the   Schouten-Nijenhuis double bracket, considered in \cite{VdB}.
\begin{defn}\cite{VdB} \label{dpbdfn}
Define $\ldb -,- \rdb: (T_A \DDA)^{\o 2} \rightarrow (T_A \DDA)^{\o 2}$
to be the unique double Poisson bracket such that
\begin{gather} \label{adbe}
\ldb a, b \rdb = 0 \quad \text{if $a,b \in A$}; \\ \label{aphidbe}
\ldb \xi, b \rdb = \xi(b) \quad \text{if $\xi \in \DDA, b \in A$}; \\
\ldb \xi, \eta \rdb = \ldb \xi, \eta \rdb_L + \ldb \xi, \eta \rdb_R, \quad \text{if $\phi, \psi \in \DDA$}.
\end{gather}
\end{defn}
The above double bracket $\ldb -, - \rdb$ is the \textit{unique}
double Poisson bracket, up to a constant factor, satisfying
\eqref{adbe}, \eqref{aphidbe} such that $\ldb \xi, \eta \rdb$ is given
by a fixed linear combination of elements \eqref{vdbelts} written to
lie in either $A \o \DDA$ or $\DDA \o A$, and such that this works for
\textit{any} algebra $A$.  This explains why the elements in
\eqref{vdbelts} form a canonical basis up to scaling (we must have the
first map to $\DDA \o A$ and the second map to $A \o \DDA$ so as to
preserve the double Leibniz rule \eqref{dleib}).

In this article, we will use a modification of Van den Bergh's double
bracket above.  There are two motivations. One is to correct the
``preference'' that the double bracket makes for preserving
right-multiplication components (Remark \ref{dbpref}), by allowing one
to keep track of \textit{both} left- and right-multiplication
components. The other motivation is to incorporate this into the
general \textit{wheelgebra} framework that allows one to relate
to differential operators, (almost-)commutative structures, etc.
\begin{defn} \label{wpbd}
Let the \textbf{wheeled Poisson bracket} on $\Ff(T_A \DDA)$ be
the unique one such that, for all $f,g \in T_A \DDA \subset \Ff_1(T_A \DDA)$,
\begin{equation}
\{f, g\} = (12) \ldb f, g \rdb.
\end{equation}
 Here,
$\ldb -,-\rdb : T_A \DDA \o T_A \DDA \rightarrow T_A \DDA \o T_A
\DDA$
is the double Poisson bracket of Definition \ref{dpbdfn}.
\end{defn}
Using the above definition, it is now true that the wheeled Poisson
bracket on $\Ff(T_A \DDA)$ transforms the $i$-th left-multiplication
on $(f \o g)$, i.e., left-composing $(f \o g)$ by
$[\Id \o (1^{\o (i-1)} \o a \o 1^{\o n-i}) \o \Id]$ for $a \in A$, to
the $i$-th left-multiplication on the element $\{f, g\}$ plus an
appropriate derivation term (i.e., a term obtained from (1) applying
$g$ to $a$ as a sum over the derivations in $\DDA$ appearing in $g$,
(2) tensoring by $f$, and (3) applying a sign and a permutation of
components); similarly for right-multiplication.  This corrects the
double Poisson phenomenon of Remark \ref{dbpref}, and preserves
additional information.
\begin{rem} \label{dpwprem} It is \textit{not} true that a wheeled
  Poisson bracket on $\Ff(A)$ is equivalent to a double bracket on
  $A$. Wheeled Poisson brackets are a generalization of double Poisson
  brackets.  Precisely, the former corresponds to wheeled Poisson
  brackets on $\Ff(A)$ of the form
  $\{a, b\} \in (12) (A \o A) \subset \Ff_2(A)$ for all
  $a, b \in A \subset \Ff_1(A)$.  In the smooth case, in the
  description of Theorem \ref{dathm}, wheeled Poisson brackets are
  given by elements of $\Ff_2(T_A \DDA)[2]$ which have zero
  Schouten-Nijenhuis bracket with themselves,  while double
  Poisson brackets are given by the subset of such elements in
  $(T_A \DDA)_{\cyc}[2]$ (this description of double Poisson brackets
  was given in \cite{VdB}).
\end{rem}
\subsection{Lie wheelgebroids and structure of $\Df(\Ff(A))$} We are
interested in the part of $\Ff(T_A \DDA)$ that involves only one copy of $\DDA$.
Thus, we consider
the space
\begin{equation}
\WDDA := \Ff_{\geq 1}(A) \o_{A^e} \DDA \cong \Ff(T_A \DDA)[1],
\end{equation}
where the $[1]$ means degree-one with respect to $|\DDA|=1$, $|A|=0$.
Here, $\Ff_{\geq 1}(A) \o_{A^e} \DDA$ means the span of elements
$a \o_{A^e} \xi$, for $a \in \Ff_{\geq 1}(A)$ and $\xi \in \DDA$; the
multiplication $\o_{A^e}$ can also be expressed as first taking $\o$
and then performing the two appropriate contractions (here, $\mu_{1,|a|} \circ \mu_{|a|+1, 1}$).

It is easy to see that the constructions of the previous section
restrict to the statement that $\WDDA$ is a Lie wheelgebra (this can
also be obtained without using $T_A \DDA$ at all: one needs only
$\ldb \xi, \eta\rdb$ for $\xi, \eta \in \DDA$, and the action of
$\DDA$ on $A$).  Furthermore, $\Ff(T_A \DDA)$ may be interpreted as 
$\Syns_{\Ff(A)} \WDDA$;
one then has the wheeled analogue of the usual result that this construction
produces a Poisson algebra.  

Also, there is a construction of universal enveloping wheelgebra of a
Lie wheelgebra (or wheelgebroid), and a ``wheeled PBW theorem'' which says
that 
$$\gr U_{\Ff(A)} \WDDA \cong \Syns_{\Ff(A)} \WDDA.$$ 
This is
another way to obtain the same wheeled Poisson bracket.  Moreover, in
the case that $A$ is smooth, $\Df(\Ff(A)) \cong U_{\Ff(A)} \WDDA$ (see
Theorem \ref{dathm} below).  This is the analogue of the classical
statement that differential operators on a manifold are the universal
enveloping algebra of the Lie algebroid of vector fields over the
space of functions.  To make these statements precise, we must
interpret $\WDDA$ not merely as a Lie wheelgebra, but as a Lie
wheelgebroid over $\Ff(A)$:
\begin{defn}
A wheelmodule $\cM$ over a wheelgebra $\cW$ is a module over $\cW$ in the
category of wheelspaces.
\end{defn}
Explicitly, a wheelmodule $\cM$ over a wheelgebra $\cW$ is a wheelspace together with
a $\k$-linear map
$
\rho: \cW \o \cM \rightarrow \cW
$
such that 
\begin{enumerate}
\item $\cW(m) \o \cM(n) \rightarrow \cW(m+n)$ is a morphism of $S_{m+n} \times S_{m+n}$-modules;
\item $\rho((w_1 \o w_2) \o x) = \rho(w_1 \o \rho(w_2 \o x))$ (for $w_1, w_2 \in W, x \in \cM$, 
\item $\rho(\mu_{i,j}(w) \o x) = \mu_{i,j} \rho(w \o x), \quad \rho(w \o \mu_{i,j}(x))
  = \mu_{i+|w|,j+|w|} \rho(w \o x)$ (for $w \in W, x \in \cM$).
\end{enumerate}
When there is no confusion, we will write $wx := \rho(w \o x)$.

A \emph{Lie wheelgebroid} is a Lie algebroid in the category of
wheelspaces. We also give this definition explicitly, since the
categorical definition of Lie algebroids is probably less familiar.
\begin{defn}
A \textit{Lie wheelgebroid} $L$ over a commutative wheelgebra $\cW$ is
a wheelmodule over $\cW$ together with 
\begin{enumerate}
\item an action $\theta: L \o W \rightarrow W$, which satisfies 
(denoting $\theta_x(w) := \theta(x,w)$)
\begin{gather}
w_1 \theta_x(w_2) = \theta_{w_1 x}(w_2), \\
\theta_x(w_1 w_2) = \theta_x(w_1) w_2 + (\sigma \times \Id) w_1 \theta_x(w_2),
\end{gather}
where $\sigma \in S_{|w_1|+|x|}$ is the permutation which sends $1, \ldots, |w_1|$ to $|x|+1, \ldots, |x|+|w_1|$, and is increasing on $1, \ldots, |w_1|$ and 
on $|w_1|+1, \ldots,|w_1|+|x|$;
\item a Lie bracket in wheelspaces,
 $\{-,-\}: L \o L \rightarrow L$, satisfying the Leibniz rule
\begin{equation} \label{wgbleib}
\{w x, y\} = (w \{x, y\}) + \omega (\theta_y(w) x),  \quad w \in W, \  x,y \in L,
\end{equation}
where $\omega \in S_{|w|+|x|+|y|}$ is the permutation which sends
$1, \ldots, |y|$ to the end ($|x|+|w|+1, \ldots,|x|+|w|+|y|$), and is increasing on $1, \ldots,|y|$ and on $|y|+1, \ldots,|y|+|x|+|w|$.
\end{enumerate}
\end{defn}
It is easy to see that $\WDDA$ is a Lie wheelgebroid over $\Ff(A)$.
Explicitly, the action of $\Ff(A)$ on $\WDDA$ is the tensor
product; the action of $\WDDA$ on $\Ff(A)$ is by applying the element
of $\DDA$ and multiplying the rest accordingly.  The Lie structure on
$\WDDA$ is the wheeled Poisson bracket.  Namely, for two elements
$\xi, \eta \in \DDA$, we apply
$(\xi \o \Id) \circ \eta - (\Id \o \eta) \circ \xi$ and similarly the
other way, to get a new element of $(\DDA \o A) \oplus (A \o \DDA)$;
to generalize to bracketing two elements of $\WDDA$ involves also
summing over applying $\xi$ to copies of $A$ or $A/[A,A]$, and
similarly for $\eta$.

For any Lie wheelgebroid $L$ over a commutative wheelgebra $\cW$, one
may define the universal enveloping wheelgebra $U_W L$ in a
straightforward way. We remark also that the Koszul complex construction
of Section \ref{lpwgbas} can be done for Lie wheelgebroids.
\begin{thm}\label{dathm} For any  {\em smooth} associative algebra $A$,
we have
\begin{enumerate}
\item[$\mathsf{(i)}$] The wheelgebra 
$\Df(\Ff(A))$ is canonically isomorphic to
  the universal enveloping wheelgebra of the Lie wheelgebroid $\WDDA$.
\item[$\mathsf{(ii)}$]  The 
  Poisson wheelgebra $\gr \Df(\Ff(A))$ is isomorphic to $\Ff(T_A \DDA)$,
  equipped with the wheeled Poisson bracket.
\end{enumerate}
\end{thm}
We will only prove (ii). It is not difficult to deduce (i) from this.

Note that, if the algebra $A$ is not smooth, we still have a natural
map of wheeled Poisson algebras
$\Ff(T_A \DDA) \rightarrow \gr \Df(\Ff(A))$.  The reason is that
$\DDA$ and $A$ always act on $\Ff(A)$, regardless of whether $A$ is
smooth. For example, any derivation $\theta \in \DDA$ acts by
\begin{equation}
\theta(a_1 \o a_2 \o \cdots \o a_n) = \sum_i \sigma_i \bigl(a_1 \o a_2 \o \cdots \o a_{i-1}
\o \theta(a_i) \o a_{i+1} \o \cdots \o a_n \bigr),
\end{equation}
where $\sigma_i \in S_{n+1}$ is an appropriate permutation.  The
action of $A$ is by tensoring on the left. This action extends to
$T_\k (\DDA \oplus A)$ by replacing tensor product of elements by
composition of the corresponding operators.  One may construct from
this a $\k$-linear map $\Ff(T_A \DDA) \rightarrow \Df(\Ff(A))$ (using
twisted symmetrization, which in fact yields a map of $\SS$-modules),
which induces a map $\Ff(T_A \DDA) \rightarrow \gr \Df(\Ff(A))$ of
commutative wheelgebras.  However, the map need not be injective or
surjective.  Similarly, we have a natural map of wheelgebras
$U_{\Ff(A)} \WDDA \rightarrow \Df(\Ff(A))$, which need not be
injective or surjective.


Before proving the theorem, we need to define the notion of principal symbol.
\begin{defn}
Define the \textit{$n$-th principal symbol}
of a differential operator $D \in \Df(B)_{\leq n}$ on a
\textit{twisted-commutative algebra} (or commutative algebra in any
appropriate symmetric monoidal category) 
$B$  to be a map  $\Gamma_n(D): B^{\o n} \rightarrow B,$ given by
$b_1 \o \cdots \o b_n\mto [[[D, b_1], b_2], \ldots, b_n](1)$.
The term ``principal symbol'' of a differential operator $D$ refers to
$\Gamma_n(D)$ where $D$ has order $\leq n$ but not order $\leq n-1$.
\end{defn}

One can prove by a straightforward computation that the principal
symbol is given by the following explicit formula (for the twisted case):
\begin{equation}
\Gamma_n(D) (b_1 \o \cdots \o b_n)
= \sum_{I \subset \{1, \ldots,n\}} (-1)^{|I|}\cdot (\Id_{S_{|D|}} \times
\sigma_I)\circ
 D\bigl(\prod_{i \in I} b_i\bigr) \prod_{j \notin I} b_j,
\end{equation}
where $\sigma_I \in S_n$ is the permutation which, if considered to
act by permuting the indices, would rearrange them in order
$1, \ldots,n$ (i.e., $\sigma_I$ sends the ordered set $I+I^c$ to
$(1,2, \ldots,n)$, where $I^c$ is the complement of $I$, $+$ is
concatenation, and $I, I^c$ are assumed to have the increasing order).
Also, the products above are in increasing order of indices from left
to right.  In the arbitrary categorical setting, we replace the
permutation $\Id_{S_{|D|}} \times \sigma_I$ by the corresponding
composition of braidings.

Note that the $n$-th principal symbol $\Gamma_n(D)$ is zero if and only if $D$ is
a differential operator of order $\leq n-1$.
\begin{proof}[Proof of Theorem \ref{dathm}] We only prove part (ii);
part (i) easily follows.
 Recall that
$\Df(A)_{\leq n}$ denotes the subspace of operators of order $\leq n$. 
We may consider the map
\begin{equation}
  \Df(A)_{\leq n} \stackrel{\Gamma_n}\too \Hom_{\k}(\Ff(A)^{\o n}, \Ff(A)),
\end{equation}
whose kernel is clearly the differential operators of order $\leq
n-1$,
and hence induces an embedding
$\gr_n \Df(A) \into \Hom_{\k}(\Ff(A)^{\o n}, \Ff(A))$. Let us study
the image.  It is straightforward to check that the image is a
derivation in all tensor components, in the sense that (for
homogeneous $a_1, \ldots, a_n, a_i' \in \Ff(A)$)
\begin{multline}
\Gamma_n(D)(a_1 \o \cdots \o a_{i-1} \o a_i a_i' \o a_{i+1} \o \cdots \o a_n)
= \mu_{|D|+i, |D|+(i+1)} (\Id_{S_{|D|}} \times \sigma_1)
 \cdot \\ \bigl(\Gamma_n(D) \cdot (a_1 \o \cdots \o a_{i-1} \o a_i' \o a_{i+1} \o \cdots \o a_n) \otimes  a_i\bigr) (\Id_{S_{|D|}} \times \sigma_1)
\\ + \mu_{|D|+i, |D|+(i+1)} (\Id \times \sigma_2) \bigl(\Gamma_n(D)(a_1 \o \cdots \o a_{i-1} \o 
a_i \o a_{i+1} \o \cdots \o a_n \bigr) \otimes a_i'\bigr) (\Id \times \sigma_2),
\end{multline}
where $\sigma_1, \sigma_2 \in S_{|a_1|+ \cdots + |a_n| + |a_{i}'|}$
are the respective permutations which
would rearrange the $a$ symbols back to the order
$a_1, \ldots, a_i, a_i', a_{i+1}, \ldots, a_n$.

The image is determined by its restriction to
$\Hom_{\k} (A^{\o n}, \Ff_{\geq n}(A))$, in view of the compatibility of contraction
maps on the domain and the image.   Using the rightmost $(A^e)^{\o n}$-bimodule structure
on $\Ff_{\geq n}(A)$, we see that these restrictions land in
\begin{multline}
\Der(A^{\o n}, \Ff_{\geq n}(A))=
\Der(A, \Der(A^{\o (n-1)}, \Ff_{\geq n}(A)))\\=\DDA \o_{A^e}
\Der(A^{\o(n-1)}, \Ff_{\geq n}(A)),
\end{multline}
where the last isomorphism is due to the smoothness of $A$.
This way, by induction,
we get $\Der(A^{\o n}, \Ff_{\geq n}(A)) \cong \DDA^{\o n} \o_{(A^e)^{\o n}} \Ff_{\geq n}(A)$.
Thus, we obtain an embedding
\begin{equation}
\gr_n \Df(A) \into (\DDA^{\o n})  \otimes_{(A^e)^{\o n}} 
(\Ff_{\geq n}(A)).
\end{equation}

It remains to show that this map is surjective.  To see this, we may
define, for any $\theta_1, \ldots, \theta_n \in \DDA$, and any
$X \in \Ff_{\geq n}(A)$, a differential operator in $\Df(\Ff(A))$ mapping to
$(\theta_1 \o \cdots \o \theta_n) \o_{(A^e)^{\o n}} X$, as follows.
For any $i$, we may consider the operator
$\theta_i: \Ff_\bullet(A) \rightarrow \Ff_{\bullet + 1}(A)$, given by
\begin{multline}
\theta_i(a_1 \o \cdots \o a_m) \\ 
= \sum_{j=1}^m (1,j) \bigl( \theta_i(a_j)' \o a_1 \o a_2 \o \cdots \o a_{j-1} \o \theta_i(a_j)'' \o a_{j+1} \o \cdots \o a_m\bigr),
\end{multline}
extended so as to be an $S_m$-module map $\Ff_m(A) \rightarrow \Ff_{m+1}(A)$ (viewing
$\Ff_{m+1}(A)$ as an $S_m$-module by the composition map $S_m \iso S_m \times S_1 \into S_{m+1}$),
and so as to be compatible with contractions.  

Then, for any $X \in \Ff_{p}(A)$ with $p \geq n$,
we define a differential operator $D$ mapping to
$(\theta_1 \o \cdots \o \theta_n) \o_{(A^e)^{\o n}} X$ as follows: for
any $Y \in \Ff_m(A)$, first consider
$Y' = \theta_1 \circ \theta_{2} \circ \cdots \circ \theta_n(Y) \in
\Ff_{m+n}(A)$.
Then, we let $D(Y) := X \o_{(A^e)^{\o n}} Y' \in \Ff_{\geq m}(A)$, considering
$Y'$ as an element of an $(A^e)^{\o n}$-module via the first $n$ components of the
$(A^e)^{\o (n+m)}$-action, and similarly $X$ via the first $n$ components
of the $(A^e)^{\o (n+m)}$-action. The result $D(Y) \in \Ff_{(p-n)+m}(A)$
has $(A^e)^{\o (p-n)+m}$-action given by: first the remaining $A^{\o(p-n)}$-action
on $X$, and then the remaining $A^{\o m}$-action on $Y'$.
Also, here, tensoring over $(A^e)^{\o n}$ is the same as taking a
tensor product over $\k$ and applying the appropriate $2n$ contraction
maps.  

It is not difficult to show that the above $D$ indeed maps to
$(\theta_1 \o \cdots \o \theta_n) \o_{(A^e)^{\o n}} X$ under the
principal symbol map.  Furthermore, it follows from the construction
that the obtained identification is compatible with permutations and
contractions (for example, applying such operations to the last $p-n$
components of $X$ above is the same as applying these operations to
the first $p-n$ components of the result).

It remains to consider the Poisson structure.  It
suffices to show that the induced wheeled Poisson bracket agrees with
the double Poisson bracket from \cite{VdB} when restricted to
$T_A \DDA \o T_A \DDA$ (it is clear that the restriction to $A \o A$
determines the bracket).  In fact, it is enough to check the
restrictions to $A \o A$, $\DDA \o A$, and $\DDA \o \DDA$.  The first
restriction is clearly zero, for degree reasons.  We have to show
that, if $\xi, \eta \in \DDA$, and $\phi, \psi \in \Df(A)_{\leq 1}$
satisfy $\Gamma_1(\phi) = \xi, \Gamma_1(\psi) = \eta$, then
\begin{equation}
[\phi, a] = \xi(a), \quad \Gamma_1[\phi, \psi] = \{\xi, \eta\}, \quad \forall a \in A.
\end{equation}
The first identity is immediate from the definition of $\Gamma_1$.  For
the second, we recall the definition of $\ldb -,- \rdb$ and $\{-,-\}$.
This says that $\{ \xi, \eta \}$, viewed as a map
$A \rightarrow \Ind_{S_3}^{S_3 \times S_3} (A \o A \o A)$ (sending
the input and output for $A$ in the domain to the third input and
third output in the image), sums over all ways to apply $\eta$
first and then $\xi$. Since we may choose $\phi$ and $\psi$ to
effectively sum over applying $\xi$ and $\eta$, respectively, to all
$A$ components that appear, this proves the desired equality.
\end{proof}

\section{Torsion of bimodule connections on $\DDA$ and $\Omega^1 A$} \label{ts}
\subsection{Sign conventions} \label{signntns}
Unless otherwise specified, we will work in the super $\Z \oplus \Z$-graded
context, with bidegrees $|A|=(0,0), |\DDA| = (1,0), |\Omega^1 A|=(0,-1)$.

Then, we use the corresponding superbraiding,
\begin{equation} \label{spbdeq}
\tau_{(21)}(M \o N) := s(|M|,|N|) (N \o M),
\end{equation}
where $s((a,b),(c,d)) = (-1)^{ac + bd}$.

The grading $|A|$ used in the definition of $\tau$ above will be called
the ``$\tau$-grading'' (to distinguish from other gradings that will arise,
e.g., $\SS$-module gradings).
\begin{ntn} \label{bartauntn}
Let $\bar \tau_\sigma$ be the non-super version of
 Notation \ref{tauntn}: $\bar \tau_{(21)}(A \o B) = B \o A$ for all $A, B$.
\end{ntn}

\subsection{Connections on a bimodule}
This subsection is a reminder of basic facts and definitions from \cite{CQ};
accordingly, we omit the citations of \cite{CQ}.
Let $A$ be an associative algebra over $\k$.  
\begin{defn} \label{lcrcd}
  Let $M$ be any $A$-bimodule.  Then a left connection $\nabla_\ell$
  on $M$ is a right $A$-module map
  $\nabla_\ell: M \rightarrow \Omega^1 A \o_A M$ satisfying
  $\nabla_\ell(am) = a \nabla_\ell(m) + da \o m$ for any
  $m \in M, a \in A$.  Similarly, a right
  connection $\nabla_r$ on $M$ is a left $A$-module map
  $\nabla_r: m \rightarrow M \o_A \Omega^1 A$ satisfying
  $\nabla_r(ma) = \nabla_r(m) a + m \o da$.
\end{defn}

Recall that being a right $A$-module map means
  $\nabla_\ell(ma) = \nabla_\ell(m) a$.
\begin{defn} A connection $\nabla=(\nabla_\ell, \nabla_r)$ on an
  $A$-bimodule $M$ is a collection of a left connection $\nabla_\ell$
  and a right connection $\nabla_r$.
\end{defn}
\begin{rem} There is also a similar notion of connection on a left
or right $A$-module.  However, a connection on an $A$-bimodule is \textbf{not}
the same as a connection on the corresponding (left or right) $A^e$-module:
the latter is a finer notion. See \cite[\S 8]{CQ}.
\end{rem}
\begin{prop} A connection on an $A$-bimodule $M$ exists if and only if $M$ is
 a projective $A^e$-module. 
\end{prop}

\begin{defn} \label{ondfn}
For each $n \geq 1$, let $\Omega^n A := (\Omega^1 A)^{\otimes_A n}$, $\Omega^0 A := A$, and
$\Omega A := \bigoplus_{n \geq 0} \Omega^n A$. If $\eta \in \Omega^n A
\subset \Omega A$, 
we write $|\eta|=n$ and say $\eta$ has degree $n$.
\end{defn}

In the proposition below,
 $|\eta| = n$ means $\eta \in M \o_A \Omega^n A$ (or
  $\Omega^n A \o_A M$ for that matter).
  \begin{prop} \label{cqop1} Any left connection on $M$ extends
    uniquely to an operator of degree one on $\Omega A \o_A M$ by the
    condition
    $\nabla(\omega \eta) = d(\omega) \eta + (-1)^{|\omega|}\omega
    \nabla(\eta)$.
    Similarly, any right connection on $M$ extends uniquely to
    $M \o_A \Omega A$ by
    $d(\eta \omega) = \nabla(\eta) \omega + (-1)^{|\eta|} \eta
    d(\omega)$.
  \end{prop}
\subsection{Torsion of a connection on $\DDA$}  \label{tcddas}
\label{tcs}
We first recall the definition of torsion in the classical case:
let $X$ be a manifold and $\nabla$ a connection on the tangent bundle
$T_X$.  The connection $\nabla$ induces an operator of degree one on
$T_X \o \Omega_X$, also denoted by $\nabla$, defined by
$\nabla(\xi \o \omega) = \nabla(\xi) \wedge \omega + \xi \o d\omega$
for $\xi \in \Gamma(U, T_X), \omega \in \Gamma(U, \Omega_X^n)$ for any
$n$ and any open subset $U \subset X$.  Restricting to degree-one forms, one notices that
$T_X \o \Omega_X^1 \cong T_X \o T^*_X \cong \End(T_X)$. So, one can
consider the element $\iota \in \Gamma(X,T_X \o \Omega_X^1)$ corresponding to the
element of $\End(T_X)$ which is the identity on fibers.  Then, the
torsion $\tau(\nabla)$ of the connection $\nabla$ is given by
$\tau(\nabla) := \nabla(\iota) \in \Gamma(X,T_X \o \Omega_X^2)$.

Equivalently, the torsion may be defined by
$\tau(\nabla)(\xi, \eta) = \nabla_\xi \eta - \nabla_\eta \xi - \{\xi,
\eta\}$,
where $\{-,-\}$ is the Lie bracket of vector fields.  A connection is
\textbf{torsion-free} if its torsion is zero.

We wish to imitate this in our setting.  
Consider the compositions
\begin{multline}
  \mu \circ (\nabla_\ell \o_A \Id): M \o_A \Omega^n A \traa^{\nabla_\ell \o_A
    \Id} \Omega^1 A \o_A M \o_A \Omega^n A
  \\ \tra^{\mu} M \o_{A^e} (\Omega^n A \o_A \Omega^1 A) \cong M \o_{A^e} \Omega^{n+1} A,
\end{multline}
\begin{equation}
  \mu \circ \nabla_r: M \o_A \Omega^n A \tra^{\nabla_r} M \o_A
  \Omega^{n+1} A \tra^{\mu} M \o_{A^e} \Omega^{n+1} A.
\end{equation}

\begin{lemma}
 For any connection $\nabla = (\nabla_\ell, \nabla_r)$ on an $A$-bimodule
  $M$, the operator 
$$(-1)^n \mu \circ (\nabla_\ell \o_A \Id) + \mu \circ \nabla_r: M \o_A
 \Omega^n A \rightarrow M \o_{A^e} \Omega^{n+1} A$$
 factors through the multiplication $\ds M \o_A \Omega^n A \tra^\mu M \o_{A^e} \Omega^{n} A$, yielding a well defined map
$\nabla: M \o_{A^e} \Omega^n A \rightarrow M \o_{A^e} \Omega^{n+1} A$.
\end{lemma}
\begin{proof} First, note that since
  $\nabla_\ell: M \rightarrow \Omega^1 A \o_A M$ is a right $A$-module
  map, the map
  $\nabla_\ell \o_A \Id: M \o_A \Omega^n A \rightarrow \Omega^1 \o_A M
  \o_A \Omega^n A$
  is well defined. Then, it remains to check that
  $\nabla' := (-1)^n \mu \circ$ $(\nabla_\ell \o_A \Id) + \mu \circ \nabla_r$
  satisfies $\nabla'(am \o_A \omega) = \nabla'(m \o_A \omega a)$.
  This follows because
  $\nabla_\ell(am) - a \nabla_\ell(m) = da \o_A m$ and
  $\nabla_r(m \o_A \omega a) - \nabla_r(m \o_A \omega) a =
  (-1)^{|\omega|} m \o_A \omega da$.
\end{proof}

Using the  identification
  $(\Omega^1 A \o_A \DDA) \o_{A^e} \Omega^1 A \cong \DDA \o_{A^e}
  \Omega^2 A$, we obtain
\begin{cor} Given any bimodule connection
  $\nabla=(\nabla_\ell, \nabla_r)$ on $\DDA$, there is a well defined
  map 
\begin{multline}
\DDA \o_{A^e} \Omega^1 A \rightarrow \DDA \o_{A^e} \Omega^2 A,
\\
\xi \o_{A^e} \omega \mapsto -\nabla_\ell(\xi) \o_{A^e} \omega +
  \nabla_r(\xi) \o_{A^e} \omega + \xi \o_{A^e} d\omega.
\end{multline}
\end{cor}
For the rest of the paper, we will assume that $\Omega^1 A$ is a
finitely-generated projective $A$-bimodule. 
\begin{defn}\label{iodfn} Let $\iota \in \DDA \o_{A^e} \Omega^1 A \cong \End_{A^e}(\DDA)$ correspond
to the identity element. If $\nabla$ is a connection on $\DDA$,
  then, the torsion $\tau(\nabla) \in \DDA \o_{A^e} \Omega^2 A$ of
  $\nabla$ is defined by $\tau(\nabla) := \nabla(\iota)$. A connection
  is \textbf{torsion-free} if the torsion is zero.
\end{defn}
 In this case, one
can identify the space
$\DDA \o_{A^e} \Omega^2 A$ in which the torsion is defined with two
other spaces:
\begin{enumerate}
\item[(i)] $\Hom_{A^e \o A^e}((\DDA \o \DDA)_{1,2}, (\DDA \o A)_{\inn, \out})$.
\item[(ii)] $\Hom_{A^e}(\Omega^1 A, \Omega^2 A)$;
\end{enumerate}
In space (i), the notation $1,2$ means that the first $A^e$ acts on
the first $\DDA$ component and the second $A^e$ acts on the second
such component, while in $(\DDA \o A)_{\inn, \out}$ the first $A^e$ has
the inner $A^e$-action and the second has the outer $A^e$-action.

In the following subsections, we will provide interpretations of 
torsion using each of these spaces.  The first will provide an
analogue of the classical formula
$\tau(\nabla)(\xi, \eta) = \nabla_\xi \eta - \nabla_\eta \xi -
\{\xi, \eta\}$,
replacing the Lie bracket with Van den Bergh's double
Schouten-Nijenhuis bracket \cite{VdB}. The second will show
equivalence with the definition of torsion for connections on
$\Omega^1$ given in \cite{CQ}, and explain how to pass from
connections on a module $M$ to connections on its dual $M^\vee$.
\subsection{The double and wheeled Schouten-Nijenhuis bracket}\label{dsnbs}
We will need the odd version of the double and wheeled Poisson bracket
from Section \ref{dwps}.  These are straightforward generalizations,
where we now consider $|\DDA| = 1, |A|=0$, and use \textit{superbraidings}
 $\tau_\sigma$ (or equivalently, we use the bigrading
explained in \S \ref{signntns}). The $\SS$-module structure of
$T_\k (T_A \DDA)$ is defined using the same signed permutations, and
hence, $\Ff(T_A \DDA)$ will be defined as in \eqref{ffadfn}, except
using \textit{supercommutators}
rather than commutators.  This will be the convention for the
remainder of this section, as well as for Section \ref{bvs}.

The super versions will be called the double and wheeled
\textit{Schouten-Nijenhuis} brackets (the double S-N bracket was first
defined in \cite{VdB}). We omit the details of this definition, since
everything is essentially the same as before.

\subsection{The formula $\tau(\nabla)(\xi, \eta) = \nabla_\xi \eta - 
\nabla_\eta \xi - \ldb \xi, \eta \rdb$}
Let us use the notation
\begin{equation}
\nabla:=\nabla_\ell + \nabla_r. \label{ndfn}
\end{equation}
The goal of this section is to prove
\begin{prop} \label{torprop}
\begin{gather} \label{torfla1}
\tau(\nabla) (\xi, \eta) = (\nabla_r)_\xi(\eta) - \tau_{(21)}(\nabla_\ell)_\eta(\xi) - 
\ldb \xi, \eta \rdb_L; \\
\label{torfla} \tau(\nabla)(\xi, \eta) - \tau_{(21)} \tau(\nabla)(\eta, \xi) = \nabla_\xi(\eta) - \tau_{(21)} \nabla_\eta(\xi) - \ldb \xi, \eta \rdb.
\end{gather}
\end{prop}
The notation $(\nabla_\ell)_\xi, (\nabla_r)_\xi, \nabla_\xi$ is explained below.
Note that \eqref{torfla} includes both \eqref{torfla1} and the result
of swapping $\xi, \eta$ in \eqref{torfla1}, in the distinct components
$\DDA \o A$ and $A \o \DDA$.
\begin{ntn} \label{pairntn} Define the pairing
  $\ctr: \DDA \o \Omega^1 A \rightarrow A^e$ by the composition (cf.~\eqref{dstardfn})
  $$ 
\xymatrix{\DDA \o \Omega^1 A \ar[rr]_<>(0.5){\sim}^{(d^*)^{-1} \o \Id}&&
 (\Omega^1
  A)^{\vee} \o \Omega^1 A \ar[rr]^<>(0.5){\phi \o \omega \mapsto \phi(\omega)}&& A
  \o A,}$$
  and similarly $\Omega^1 A \o \DDA \rightarrow A^e$ by identifying
  first $\DDA$ with $(\Omega^1 A)^\vee$ using $d^*$. That is, one has
  $\xi \ctr \omega = ((d^*)^{-1} \xi) (\omega)$.
\end{ntn}
\begin{ntn} We denote
\begin{equation} \label{nsub}
(\nabla_\ell)_\xi \eta := (\xi \o 1) \ctr \nabla_\ell(\eta); \quad
(\nabla_r)_\xi \eta := (1 \o \xi) \ctr \nabla_r(\eta); \quad
\nabla_\xi := (\nabla_\ell)_\xi + (\nabla_r)_\xi.
\end{equation}
\end{ntn}
  One may use the pairing $\ctr$ to fix an isomorphism
\begin{equation}
\Omega^1 A \iso (\DDA)^\vee: 
  \omega \mapsto \bigl( \xi \mapsto \tau_{(21)} (\xi \ctr \omega)
  \bigr).
\end{equation}
 The $\tau_{(21)}$ is needed because $(\DDA)^\vee = \Hom_{A-\text{bimod}}(\DDA, (A \o A)_\out)$ carries the $A$-bimodule
structure on $\DDA$ to the \textbf{outer} bimodule structure on $A \o A$.
We always take $A$ to have degree $0$ with respect to the super structure, 
so that $\tau_{(21)}$ is 
the unsigned flip.  This will henceforth be the assumed isomorphism $\Omega^1 A \iso (\DDA)^\vee$.  See the following Caution.
\begin{caut}\label{doc} Note that the tautological contraction
  $\DDA^\vee \o_{A^e} \DDA \rightarrow A \o A$ differs from the
  contraction $\Omega^1 A \o_{A^e} (\Omega^1 A)^\vee$ via the flip: 
  One has the commutative diagram
\begin{equation}
\xymatrix{
\DDA^\vee \o_{A^e} \DDA \ar[r] \ar[d]^\sim 
& A \o A \\
\Omega^1 A \o (\Omega^1 A)^\vee \ar[r] & A \o A \ar[u]^{\tau_{(21)}}
}
\end{equation}
By Notation \eqref{pairntn}, the contraction $\Omega^1 A \o_{A^e} \DDA \rightarrow
A \o A$ uses the bottom arrow (identifying $\DDA \cong (\Omega^1 A)^{\vee}$);
when we actually mention $(\DDA)^\vee$, then one will use the top
arrow.
\end{caut}

  If we have multiple components, we will use notation of the sort
  $(\xi_1 \o \xi_2) \ctr (\omega_1 \o \omega_2) = (\xi_1 \ctr
  \omega_1) \o (\xi_2 \ctr \omega_2)$
  (we can replace the $\o$ between $\omega_1$ and $\omega_2$ by
  $\o_A$; then the whole element can be
  replaced by an element of $\Omega^2 A$, by Definition \ref{ondfn}).

Let's make our choice of the pairing $\ctr$ more explicit.  
\begin{ntn} \label{swntn}
  For any element $f \in M \o N$, we use the Sweedler summation notation
  $f = f' \o f''$, shorthand for $f = \sum_i f_i' \o f_i''$.  So, for
  example, $(g \ctr f'') \o f'$ means $\sum_i (g \ctr f_i'') \o f_i'$.
  Similarly, $f \in M_1 \o \cdots \o M_n$ can be described as
  $f = f^{(1)} \o \cdots \o f^{(n)}$, with superscripts of $(i)$
  equivalent to $i$ primes.
\end{ntn}
\begin{lemma}
By the choice in Notation \ref{pairntn}, one has
\begin{equation} \label{paireqn}
  (a \xi b) \ctr (c \omega d) = c (\xi \ctr \omega)' b \o a 
(\xi \ctr \omega)'' d; \quad \xi \in \DDA, \omega \in \Omega^1 A.
\end{equation}
\end{lemma}
The proof is immediate.  
\begin{proof}[Proof of Proposition \ref{torprop}]
For convenience, we use the notation $(1 \o x) \ctr (a \o b) := a \o (x \ctr b)$ and similarly with $1$ appearing in other components, e.g., 
$(x \o 1) \ctr (a \o b) := (x \ctr a) \o b$.

We can pair $\tau(\nabla) \in \DDA \o_{A^e} \Omega^2 A$ with elements
$\xi, \eta \in \DDA$ using the isomorphism $\Omega^2 A \cong \Omega^1
A \o_A \Omega^1 A$. Precisely, define $\tau(\nabla)(\xi, \eta)$ in
$\DDA \o A$ as the image under the composition
\begin{multline}
\tau(\nabla) \in \DDA \o_{A^e} \Omega^2 A \cong \DDA \o_{A^e} (\Omega^1 A \o_A \Omega^1 A) \\ \traa^{(1 \o \xi \o \eta) \ctr} \DDA \o_{A^e} ((A \o A) \o_A (A \o A)) \cong \DDA \o A.
\end{multline}

Recall that $\iota \in \DDA \o_{A^e} \Omega^1 A$ corresponds to
$\Id \in \End_{A^e}(\DDA)$ (and $\Id \in \End_{A^e}(\Omega^1 A)$).  Let us suppose
that $\iota = \sum_s \xi_s \o_{A^e} \omega_s$ for some
$\xi_s \in \DDA$ and $\omega_s \in \Omega^1 A$.  We will need to use
the resulting identities
\begin{gather}
\sum_s (\xi \ctr \omega_s)''\cdot  \xi_s \cdot (\xi \ctr \omega_s)' = \xi, \quad
\sum_s (\xi_s \ctr \omega)'\cdot \omega_s\cdot (\xi_s \ctr \omega)'' = \omega,
\end{gather}
which follow from definitions (see \eqref{paireqn}).  
For any element $\xi \in \DDA$, let us use the notation $\theta_\xi: A \rightarrow A \o A$ for the associated map (to avoid confusion with multiplying by $\xi$ in $T_A \DDA$).
We then have, using the natural
identifications $(\DDA) \o_{A} (A \o A) \cong \DDA \o A$, and $(A \o A) \o_A \DDA \cong A \o \DDA$:
\begin{multline}\label{toreq0}
  \tau(\nabla)(\xi, \eta) = -\sum_s \tau_{21} \circ \bigl( (\eta \o 1) \ctr
  \nabla_\ell(\xi_s) \bigr) \o_{A^e} \bigl(\xi \ctr
  \omega_s\bigr) \\ + \bigl((1 \o \xi) \ctr
  \nabla_r(\xi_s)\bigr) \o_{A^e} \bigl(\eta \ctr
  \omega_s\bigr) + \xi_s \o_{A^e} \bigl((\xi \o \eta) \ctr d
  \omega_s\bigr),
\end{multline}
\begin{equation}\label{toreq1} -\sum_s \bigl( (\eta \o 1) \ctr
  \nabla_\ell(\xi_s) \bigr) \o_{A^e} \bigl(\xi \ctr
  \omega_s\bigr)  = -(\nabla_\ell)_{\eta} \xi + \sum_s
  \theta_\eta((\xi \ctr \omega_s)'') \o_{A} \xi_s \cdot (\xi \ctr
  \omega_s)',
\end{equation}
and
\begin{equation}\label{toreq2}
  \sum_s \bigl((1 \o \xi) \ctr
  \nabla_r(\xi_s)\bigr) \o_{A^e} \bigl(\eta \ctr
  \omega_s\bigr) =  (\nabla_r)_{\xi} \eta - \sum_s (\eta \ctr \omega_s)'' \xi_s \o_{A} \theta_\xi((\eta \ctr \omega_s)').
\end{equation}
Let $i_1: (A \o A)_{\inn} \o_{A^e} (A \o A \o A) \iso A \o A \o A$ be given by $i_1((a \o b)_{\inn} \o_{A^e} (c \o e \o f)) = ac \o fb \o e$, which is what is needed to have the commutative diagram
\begin{equation}
\xymatrix{
\DDA \o_{A^e} (A \o A \o A) \ar[rd]^{\ctr (\omega \o 1)} \ar[r]^{\sim}
& \DDA \o A \ar[r]^{\ctr (\omega \o 1)} & A \o A \o A \\
 & (A \o A)_{\inn} \o_{A^e} (A \o A \o A) \ar[ru]^{i_1}
}
\end{equation}
Then one has:
\begin{multline}\label{toreq3}
\sum_s i_1 \bigl((\xi_s \ctr \omega)_\inn \o_{A^e} \bigl((\xi \o \eta) \ctr d
  \omega_s\bigr)\bigr) = \tau_{(32)} \Bigl((\xi \o \eta) \ctr d\omega  \\
-\sum_s \bigl(\theta_\xi((\xi_s \ctr \omega)') \o_{A} \eta \ctr \omega_s 
(\xi_s \ctr \omega)'' + (\xi_s \ctr \omega)' \xi \ctr \omega_s
\o_{A} \theta_\eta((\xi_s \ctr \omega)'')\bigr)\Bigr).
\end{multline}
Contracting \eqref{toreq1}, \eqref{toreq2} with $\omega$ and adding
\eqref{toreq3}, one obtains
\begin{multline}\label{bigtoreq}
(-\tau(\nabla)(\xi, \eta) + (\nabla_r)_\xi(\eta) -\tau_{(21)} (\nabla_\ell)_\eta(\xi))\ctr
\omega \\= \sum_s \tau_{(32)} \Bigl(-(\xi \o \eta) \ctr d \omega + \bigl( \theta_\xi\bigl((\xi_s \ctr \omega)' (\eta \ctr \omega_s)'\bigr) \o 
(\eta \ctr \omega_s)'' (\xi_s \ctr \omega)''\bigr) \\ - \bigl( 
(\xi_s \ctr \omega)'(\xi \ctr \omega_s)' \o \theta_\eta \bigl((\xi \ctr \omega_s)'' (\xi_s \ctr \omega)''\bigr) \bigr) \Bigr)
\\ = \tau_{(32)} \bigl(-(\xi \o \eta) \ctr d \omega + (\theta_\xi \o 1)(\eta \ctr \omega) - (1 \o \theta_\eta) (\xi \ctr \omega)\bigr) = \ldb \xi, \eta \rdb_L \ctr (\omega \o 1).
\end{multline}
For the last equality, we note that $\ldb \xi, \eta \rdb_L$ is defined by
$\ldb \xi, \eta \rdb_L \ctr (da \o 1) = \tau_{(32)} \bigl((\theta_\xi \o 1)\circ \theta_\eta(a) - (1 \o
\theta_\eta) \circ \theta_\xi(a) \bigr)$;
there is then a unique extension of $\ldb \xi, \eta \rdb_L$ to a map
$\Omega^1 A \o A \rightarrow A \o A \o A$ as indicated in the last
line.
\end{proof}

\subsection{Torsion of connections on $\Omega^1$}
Note: This subsection and the next will not be used elsewhere in this paper.

It turns out that torsion of connections on $\Omega^1$ is easier
to define. In the classical setting, a connection $\nabla$ on
$\Omega^1_X$ for a manifold $X$ is a map
$\nabla: \Omega^1_X \rightarrow \Omega^1_X \o \Omega^1_X$ which is a
derivation in the sense that
$\nabla(a \omega) = a \nabla(\omega) + \omega \o da$, for local sections $\omega$ and $a$
of $\Omega^1_X$ and $\mathcal{O}_X$, respectively.  Then, to
compare with
$d: \Omega^1_X \rightarrow \Omega^2_X \cong \La ^2 \Omega^1_X$,
let us define $q: \Omega^1_X \o \Omega^1_X \onto \Omega^2_X$ to be the
quotient. Then we may consider
\begin{equation}
\tau(\nabla):=q \circ \nabla + d: \Omega^1_X \rightarrow \Omega^2_X,
\end{equation}
which we can call the \textbf{torsion} of $\nabla$.

The definition of torsion in \cite{CQ} for connections on $\Omega^1 A$ is then the noncommutative analogue of the above.  Namely, in the noncommutative case $\Omega^2 A = \Omega^1 A \o_A \Omega^1 A$, so $\nabla_r, \nabla_\ell,$ and $d$ are all maps $\Omega^1 \rightarrow \Omega^2 A$.  To get an $A^e$-module map, it is clear that one considers the combination
\begin{equation} \label{o1tordfn}
\tau(\nabla) := -\nabla_\ell + \nabla_r + d,
\end{equation}
which is called the \textbf{torsion} of $\nabla$. 
\begin{rem}
As is
mentioned in \cite{CQ}, given any choice of torsion and left or right
connection, \eqref{o1tordfn} defines a right or left connection: so
there are one-to-one correspondences
\begin{equation}
\xymatrix{
\text{left connections on $\Omega^1 A$} \ar@{<->}[rr]^{\nabla_\ell
  \mapsto \nabla_\ell - d} & &
\text{right connections on $\Omega^1 A$}
}
\end{equation}
and
\begin{equation}
\ds \text{left connections on $\Omega^1 A$ with choice of torsion} \longleftrightarrow
\text{connections on $\Omega^1 A$}.
\end{equation}
One may easily see the corresponding statement
for connections on $\DDA$.
\end{rem}

\subsection{Dual connections and equivalence of torsion with dual torsion}
Classically, given a connection
$\nabla: E \rightarrow E \o \Omega^1_X$ on a vector bundle $E$, one
may define a natural dual connection $\nabla^\vee$ on $E^\vee$, such
that, if $E$ is of finite-rank, $(\nabla^\vee)^\vee \cong \nabla$
under the canonical isomorphism $E \cong (E^\vee)^\vee$.  Namely, one
uses the formula (for $f \in \Gamma(U,E^\vee), e \in \Gamma(U,E)$):
\begin{equation}
\nabla^\vee(f) \ctr e := d(f \ctr e) - (f \o 1) \ctr \nabla(e).
\end{equation}

For a (projective) bimodule $M$ over a ring $A$ with connection $\nabla=(\nabla_\ell, \nabla_r)$, one may define a dual connection $\nabla^\vee$ on $M^\vee
:= \Hom_{A^e}(M, A^e)$.  Let us define the pairing $\ctr: M^\vee \o M \rightarrow A$ by applying $M^\vee$ to $M$, which means that one has \eqref{paireqn}, considering $\xi \in M^\vee$ and $\omega \in M$.

Let us first try to dualize $\nabla_\ell$. It is natural to consider
the two possible compositions
$\ds M^\vee \o M \tra^{\ctr} A \o A \tra^{d \o 1} \Omega^1 A \o A$ and
$\ds M^\vee \o M \tra^{\ctr} A \o A \tra^{1 \o d} A \o \Omega^1 A$.  
As for using the connection $\nabla_\ell$, it is natural to consider
$\ds M^\vee \o M \tra^{1 \o \nabla_\ell} M^\vee \o \Omega^1 \o_A  M \tra^\ctr 
(\Omega^1 \o_A A) \o A \cong \Omega^1 \o A$.

In the end, we will need something that is $A^e$-linear in $M$ (so as
to get a map $M^\vee \rightarrow M^\vee \o \Omega^1 A$ or
$M^\vee \rightarrow \Omega^1 A \o M^\vee$). As the latter map is
right $A$-linear in $M$ and a derivation on the left, we need
to consider
\begin{equation}\bigl((d \o 1) \circ \ctr\bigr) - \bigl(\ctr \circ (1 \o \nabla_\ell)\bigr): M^\vee \o M
\rightarrow \Omega^1 A \o A.
\end{equation}
In $M^\vee$, this is left $A$-linear and a derivation on the right.
Hence the resulting map should be considered as a \textbf{right} connection: dualizing left connections results in right connections.  Finally, we define
the operations $\ctr m: \Omega^1 A \o_A M^\vee \rightarrow A \o \Omega^1 A, M^\vee \o_A \Omega^1 A \rightarrow \Omega^1 A \o A$ by $(\omega \o_A f) \ctr m = (f \ctr m)' \o \omega (f \ctr m)'', (f \o_A \omega) \ctr m = (f \ctr m)' \omega \o (f \ctr m)''$. We then have the
\begin{lemdfn}
For any left connection $\nabla_\ell:M \rightarrow \Omega^1 A \o_A M$ on an $A^e$-module $M$, the map $\nabla_r^\vee: M^\vee \rightarrow M^\vee \o_A \Omega^1 A$ given by
\begin{equation} \label{rdualc}
\nabla_r^\vee(f) \ctr m =  \bigl( (d \o 1) (f \ctr m) - \tau_{(21)}
(1 \o f) \ctr \nabla_\ell(m) \bigr), \quad(\text{given }\nabla_\ell, 
\text{ a left conn.}) 
\end{equation}
is a right connection.  Similarly, if $\nabla_r$ is right connection,
\begin{equation} \label{ldualc}
\nabla_\ell^\vee(f) \ctr m = \bigl( (1 \o d) (f \ctr m) - 
\tau_{(21)} (f \o 1) \ctr \nabla_r(m) \bigr), 
\quad(\text{given }\nabla_r, \text{ a right conn.})
\end{equation}
defines a left connection.  So, for any bimodule connection
$\nabla=(\nabla_\ell, \nabla_r)$, one may define a dual connection
$\nabla^\vee=(\nabla^\vee_\ell, \nabla^\vee_r)$ by
\eqref{rdualc}, \eqref{ldualc}.
\end{lemdfn}
Note that, if $m \mapsto m^{\vee \vee}$ under $M \iso M^{\vee \vee}$, and $f \in M^\vee$, then $f \ctr m = \tau_{(21)} (m^{\vee \vee} \ctr f)$.  
It is then immediate that
\begin{lemma}
If $M$ is finitely-generated (projective), then $(\nabla^\vee)^\vee \cong
\nabla$ under the natural isomorphism $(M^\vee)^\vee \cong M$.
\end{lemma}

We now prove that, in both the classical and the noncommutative cases,
torsion of a connection on $T_X$ ($\DDA$) and of a connection on
$\Omega^1_X$ ($\Omega^1 A$) are identical (under the appropriate
natural identifications of spaces).
\begin{prop}
  Let $\nabla$ be a connection on $\Omega^1_X$ and $\nabla^\vee$ its
  dual connection on $T_X$.  Then
  $\tau(\nabla)=q \circ \nabla + d: \Omega^1_X \rightarrow \Omega^2_X$
  is naturally identified with
  $\tau(\nabla^\vee) = (\nabla^\vee \o 1)(\iota)$ in $\Gamma(X, T_X \o \Omega^2_X)$
  (or with $\tau(\nabla^\vee): T^2_X \rightarrow T_X$ given by
  $\tau(\nabla^\vee)(\xi, \eta) = \nabla^\vee_\xi \eta -
  \nabla^\vee_\eta \xi - \{\xi, \eta\}$).
\end{prop}
\begin{proof}
  We show equivalence of $q \circ \nabla + d$ with
  $\tau(\nabla^\vee)(\xi, \eta) = \nabla^\vee_\xi \eta -
  \nabla^\vee_\eta \xi - \{\xi, \eta\}$. Let us consider
\begin{equation}
  (q \circ \nabla(\omega) + d \omega)(\xi, \eta) = 
  \xi \ctr (\nabla_\eta \omega) - \eta \ctr (\nabla_\xi \omega) 
  + \xi (\eta \ctr \omega) - \eta (\xi \ctr \omega) - 
  \{\xi, \eta\} \ctr \omega,
\end{equation}
while 
\begin{equation}
(\nabla^\vee_\xi \eta - \nabla^\vee_\eta \xi) \ctr (\omega \o 1)
= \xi(\eta \ctr \omega) - \eta(\xi \ctr \omega) - \eta \ctr \nabla_\xi \omega + \xi \ctr \nabla_\eta \omega,
\end{equation}
proving the desired result.
\end{proof}

\begin{prop}
  Let $\nabla$ be a bimodule connection on $\Omega^1 A$ and
  $\nabla^\vee$ the dual connection on $\DDA$.  Then
  $\tau(\nabla): \Omega^1 A \rightarrow \Omega^2 A$ is naturally
  identified with
  $\tau(\nabla^\vee)\in \DDA \o_{A^e} \Omega^2 A
  \cong  \Hom_{A^e \o A^e}((\DDA \o \DDA)_{1,2}, (\DDA \o A)_{\inn, \out})$, 
  where the last 
 space is where the formula
\eqref{torfla1} lives (see the end of Section \ref{tcs}).
\end{prop}
\begin{proof}
  As before, we show equivalence of $-\nabla_\ell + \nabla_r + d$ with
  $\xi \o \eta \mapsto (\nabla^\vee_r)_\xi \eta - \tau_{(21)}
  (\nabla^\vee_\ell)_\eta \xi - \{\xi, \eta\}_\ell$.
  We first note that
  $f \in \Hom_{A^e}(\Omega^1 A, \Omega^2 A)$ corresponds to
  $f' \in \Hom_{A^e \o A^e}((\DDA \o \DDA)_{1,2}, (\DDA \o
  A)_{\inn, \out})$
  if and only if we have
  $(\xi \o \eta) \ctr (f(\omega)) = \tau_{(32)} f'(\xi \o \eta) \ctr
  (\omega \o 1)$,
  which can be checked (for example) by looking at how each side of the
  formula changes when $\xi, \eta,$ and $\omega$ are acted on by
  $A^e$.

Now, we compute
\begin{multline} \label{dp1}
(\xi \o \eta) \ctr 
\bigl(-\nabla_\ell(\omega) + \nabla_r(\omega) + d \omega\bigr) =
-(1 \o \eta) \ctr (\nabla_\ell)_\xi(\omega) + (\xi \o 1) \ctr (\nabla_r)_\eta(\omega) \\+ \theta_{\xi \o 1}(\eta \ctr \omega) - \theta_{1 \o \eta} (\xi \ctr \omega) - \tau_{(32)} \bigl(\{\xi, \eta\}_\ell \ctr (\omega \o 1)\bigr),
\end{multline} 
where we used the last line of \eqref{bigtoreq} to expand $(\xi \o \eta) \ctr d \omega$.

Now, let us expand $\tau(\nabla^\vee)(\xi, \eta) +\{\xi, \eta\}_\ell =
(\nabla^\vee_r)_\xi \eta - \tau_{(21)}
  (\nabla^\vee_\ell)_\eta \xi$ applied to $\omega \o 1$:
\begin{multline}
\bigl((\nabla^\vee_r)_\xi \eta - \tau_{(21)}(\nabla^\vee_\ell)_\eta \xi\bigr) \ctr (\omega \o 1) = \tau_{(32)} \bigl(\xi \ctr 
(\nabla^\vee_r \eta \ctr \omega) - \eta \ctr 
(\nabla^\vee_\ell \xi \ctr \omega)\bigr)\\=
\tau_{(32)} \bigl(\theta_{\xi \o 1}(\eta \ctr \omega) - (1 \o \eta) \ctr (\nabla_\ell)_\xi(\omega) - \theta_{1 \o \eta} (\xi \ctr \omega) + (\xi \o 1) \ctr (\nabla_r)_{\eta}(\omega) \bigr),
\end{multline}
which proves that
$\tau_{(32)} (\tau(\nabla^\vee)(\xi, \eta) \ctr (\omega \o 1)) =$ the
RHS of \eqref{dp1}, as desired.
\end{proof}
\section{The BV operator $D_\nabla$} \label{bvs}
\subsection{The classical story} \label{cdns} In this subsection we
briefly recall a classical construction, following \cite{Kz}, of BV
structures on $\La  T_X$ for a finite-dimensional smooth manifold
$X$, which generate the Schouten-Nijenhuis bracket.

Let $\nabla$ be a connection on $T_X$.  The connection extends to a
connection $\nabla: \La ^n T_X \rightarrow \La ^n T_X \o 
\Omega^1_X$ satisfying
$\nabla(\xi \wedge \eta) = \nabla(\xi) \wedge \eta + \xi \wedge
\nabla(\eta)$.
(Note that the derivation property of $\nabla$ does not allow one to
put a sign such as $(-1)^{|\xi|}$ in front of the second term.)

Let $\iota \in \Gamma(X, T_X \o \Omega^1_X)$ be the canonical section
which is the identity on fibers.  Given a section of
$\La ^n T_X \o \La ^m \Omega^1_X$, one may consider the contraction
$i_\iota$ with $\iota$: this means
\begin{multline}
i_\iota(\xi_1 \wedge \cdots \wedge \xi_n \wedge \omega_1 \wedge \cdots \wedge \omega_m) \\= \sum_{j=1}^n \sum_{\ell=1}^m (-1)^{i+j} i_\iota(\xi_i \o \omega_j)
\xi_1 \wedge \cdots \wedge \hat \xi_i \wedge \cdots \wedge \xi_n \wedge \omega_1 \wedge \cdots \wedge \hat \omega_j \wedge \cdots \wedge \omega_m.
\end{multline}
One may then consider
\begin{equation}
  D_\nabla := 
  i_\iota \circ \nabla: \La ^n T_X \rightarrow \La ^{n-1} T_X,
\end{equation}
which can also be written as
\begin{equation}
D_\nabla = \sum_s i_{\omega_s} \nabla_{\xi_s}, \quad \iota = \sum_s \xi_s \o \omega_s.
\end{equation}
It is easy to see that $D_\nabla$ is a differential operator of order
$\leq 2$.
Since $\nabla$ is torsion-free if and only if
$\nabla_\xi \eta - \nabla_\eta \xi = \{\xi, \eta\}$ for all
$\xi, \eta \in \Gamma(U,T_X)$ (where $U \subset X$ denotes any
open subset), one may easily show
\begin{prop}\cite{Kz} The connection $\nabla$ is torsion-free if and only if
the principal symbol of $D_\nabla$ (as an
operator of order $2$) is $\pm$ the Schouten-Nijenhuis bracket; precisely,
\begin{equation} \label{cdnsb}
D_\nabla(\xi \wedge \eta) - \xi D_\nabla(\eta) - D_\nabla(\xi) \eta
  = (-1)^{|\xi|+1} \{\xi, \eta\} \text{ for all $\xi, \eta \in \Gamma(U, \La  T_X)$}.
\end{equation}
\end{prop}

Furthermore, one may get a formula for $D_\nabla^2$. First, \cite{Kz}
remarks that $D_\nabla^2 = \frac{1}{2} [D_\nabla, D_\nabla]$ where
$[,]$ is the supercommutator (using degree of differential operator,
{\em not} order: $D_\nabla$ lowers degree by one), which shows that
$D_\nabla^2$ has order $\leq 3 = 2+2-1$ as a differential operator
(because $D_\nabla$ itself has order $\leq 2$). Then, \cite{Kz}
computes the principal symbol $\Gamma_3(D_\nabla^2)$ as an operator of
degree $\leq 3$, in terms of the Jacobiator of the principal symbol
$\Gamma_2(D_\nabla)$ (as an operator of degree $\leq 2$).  Namely,
$\Gamma_3(D_\nabla^2) = \Gamma_3(\frac{1}{2} [D_\nabla, D_\nabla]) =
\frac{1}{2} [\Gamma_2[D_\nabla], \Gamma_2[D_\nabla]]$.
This means that, if $\Gamma_2(D_\nabla)$ satisfies the Jacobi
identity, then $\Gamma_3(D_\nabla^2)=0$, so $D_\nabla^2$ is actually
an operator of degree $\leq 2$.

Now, $\nabla$ is torsion-free if and only if $\Gamma_2(D_\nabla)$ is the
Schouten-Nijenhuis bracket. In this case, $D_\nabla^2$ must have order
$\leq 2$.  Since it is also an operator of degree $-2$, it must be
contraction with a two-form.

It remains to compute this two-form. Koszul gives the following formula:
\begin{prop}\label{cdnsfp}\cite{Kz} One has
\begin{equation}\label{cdnsf}
D_\nabla^2 = i_{\tr(\nabla^2)},
\end{equation}
where $\nabla^2: T_X \rightarrow T_X \o \Omega^2_X$ is the curvature and
$\tr(\nabla^2) \in \Gamma(X, \Omega^2_X)$ is its trace, which can also
be written as $\ds \tr(\nabla^2): \
\La ^{\text{dim } X} T_X \stackrel{\nabla^2}\too
\La ^{\text{dim } X} T_X \o \Omega^2_X.$
\end{prop}
In particular, $D_\nabla$ gives a BV structure generating the
Schouten-Nijenhuis bracket in the case that $\tr(\nabla^2) = 0$ and
$\nabla$ is torsion-free.  

In Koszul's paper, the verification of \eqref{cdnsf} is omitted; this
seems to be the hardest technical part of the proof.  In \S \ref{pfcdnsf}, 
we will give a new proof
of \eqref{cdnsf}, since we will need to apply the same proof to the
noncommutative setting. Our proof actually works in the smooth
algebraic setting, and is based on purely global arguments on any
affine variety, replacing the tangent bundle by its global sections
($\Der(A)$, where $X = \Spec A$), viewed as a projective module.

Let us return to Koszul's setting.  Note that the $D_\nabla$ are all
the BV structures possible which generate the Schouten-Nijenhuis
bracket: any two differential operators of order $2$ with the same
principal symbol must differ by an operator of order $\leq 1$ (a
derivation).  If, furthermore, the operators have degree $-1$ (as is
the case here), then the resulting difference must be given by an
$\mathcal{O}_X$-linear map $T_X \rightarrow \mathcal{O}_X$, i.e., contraction with
a global one-form $\omega$.  Similarly, two torsion-free
connections differ by a linear map $T_X \rightarrow T_X \o \Omega^1_X$,
i.e., a global $\End(T_X)$-valued one-form $\beta$, such that $(\eta \ctr \beta)(\xi) = (\xi \ctr
\beta)(\eta)$ for all $\xi, \eta \in \Gamma(X, T_X)$. 
It remains only to
show that one can produce, for every global section $\omega \in \Gamma(X, \Omega^1_X)$,
such a one-form $\beta$, so that $\tr(\beta) = \omega$. This is clear.
Thus, one can always go from a differential operator $D$ with principal
symbol $[-,-]$ to a connection $\nabla$ such that $D = D_\nabla$.

Finally, it is evident by \eqref{cdnsb} that $D_\nabla$ is determined
by its restriction to vector fields, where one sees that
\begin{equation}
  D_\nabla\bigl|_{T_X} = \div \nabla := \tr \circ \nabla: T_X \tra^{\nabla} T_X \o \Omega^1_X \tra^{\ctr} \mathcal{O}_X,
\end{equation}
where $\ctr: T_X \o \Omega^1_X \rightarrow \mathcal{O}_X$ is the
contraction.  Hence, two connections $\nabla, \nabla'$ induce the same
differential operator if and only if they have the same divergence:
$\div \nabla = \div \nabla'$. (One may also express $\tr \circ \nabla$
as the trace of the $\mathcal{O}_X$-linear map
$T_X \rightarrow \End_k(T_X) \cong T_X \o \Omega^1_X, \xi \mapsto
\nabla_\xi(\text{\textbf{---}})$.)  Summarizing:
\begin{prop}\cite{Kz}  The map $\nabla \rightarrow D_\nabla$ gives
  a one-to-one correspondence between torsion-free connections up to
  equivalence, and differential operators of order $2$ on $T_X$ whose
  principal symbol is $\pm$ the Schouten-Nijenhuis bracket. Here
  $\nabla, \nabla'$ are ``equivalent'' if and only if they have the
  same divergence.  Furthermore, torsion-free connections which are
  trace-flat ($\tr(\nabla^2)= 0$) correspond to BV structures
  generating the Schouten-Nijenhuis bracket.
\end{prop}

\subsection{BV and Gerstenhaber wheelgebras}\label{tbvps}
Now, we pursue a wheeled analogue of the preceding
section.  The first step is to define BV and Gerstenhaber wheelgebras:
these are just BV and Gerstenhaber algebras in the category of
wheelspaces.

To interpret BV wheelgebras in terms of differential operators on a
Gerstenhaber wheelgebra (which is, in particular, a supercommutative wheelgebra), we
need to apply the categorical definition of differential operators to
the category of super wheelspaces: these are wheelspaces
equipped with an additional $\Z$-grading, with the super (i.e.,
signed) braidings. We will use ``degree'' to refer to this new degree,
and  ``wheeled degree'' to refer to the wheelspace
degree.

We thus obtain an almost-supercommutative wheelgebra of differential
operators.  
Explicitly, $D$ has degree $m$ if
$|D(x)|=|x|+m$ for any homogeneous $x$, and then
$[D,D'] := D \circ D' - (-1)^{|D| |D'|} D' \circ D$.

\subsection{Overview of wheeled version of $D_\nabla$}\label{ovwsec} We would like
to mimic Section \ref{cdns} in the wheeled (noncommutative geometry) context, by
replacing classical notions with the twisted versions of Section
\ref{tbvps}.  To do this, we first need to establish some
preliminaries concerning multilinearity and connections in the
noncommutative case.  To demonstrate how this is important, let us
begin with a bimodule connection $\nabla=(\nabla_\ell, \nabla_r)$ on
$\DDA$.  We note that $\nabla$ induces a ``connection'' on $T_A \DDA$
as follows (recall \eqref{ndfn}):
\begin{gather}
\nabla: T_A^n \DDA \rightarrow \bigoplus_{0 \leq i \leq n} T_A^i \DDA \o_A \Omega^1 A \o_A T_A^{n-i} \DDA; \\
\nabla: \xi_1 \o_A \cdots \o_A \xi_n \mapsto \sum_i \xi_1 \o_A \cdots
\o_A \nabla(\xi_i) \o_A \cdots \o_A \xi_n.
\end{gather}
This is well defined because $\nabla_r$ in the $i$-th component and
$\nabla_\ell$ in the $i+1$-th component are compatible in the sense
that 
\begin{equation}
\nabla_r(\xi a) \o_A \eta +(\xi a) \o_A \nabla_\ell(\eta) = \nabla_r(\xi)
\o_A (a \eta) + \xi \o_A \nabla_\ell(a \eta).
\end{equation}
Then, we would like to define $D_\nabla$ by $i_\iota \circ \nabla$.
However, when we contract with $\iota$, the result does not live in
$T_A \DDA$ anymore: one may see (by Corollary \ref{fbvc}) that the
result makes sense in a space of the form
$T_A \DDA \o (T_A \DDA / [T_A \DDA, T_A \DDA])$.

To iterate, we can proceed by using the space $\Ff(T_A \DDA)$. For
example,
the connection $\nabla$ descends to a well defined map
\begin{equation}
\nabla: (T_A \DDA)_\cyc \rightarrow  
\bigoplus_{0 \leq i \leq n}
T_A \DDA \o_{A^e} \Omega^1 A.
\end{equation} 
On this latter sum, we may still define the contraction $i_\iota$, which
lands back in $\Ff(T_A \DDA)$.

\subsection{The wheeled Schouten-Nijenhuis bracket}
The wheeled Schouten-Nijenhuis bracket is the super version of
the wheeled Poisson bracket from Definition \ref{wpbd}, with
$|\DDA|=1, |A|=0$.  (This gives the structure of a Gerstenhaber wheelgebra).

In particular, one obtains an ordinary Gerstenhaber bracket on
$\Ff_0(\DDA)$.  This identifies with the supersymmetric algebra on
 $(T_A \DDA)_\cyc$. The latter Lie
algebra, noticed in \cite{VdB}, is a generalization of the ``necklace
Lie algebra,'' defined in \cite{BLB, G}.

\subsection{Connections and bimodule contractions on spaces $\Ff(T_A M)$}
In this section, we will use \textit{bimodule contractions} to refer
to operations of the sort $\xi \ctr \omega$, for
$\xi \in \DDA, \omega \in \Omega^1 A$, or more generally,
$\xi \in M, \omega \in M^\vee$, for any finitely-generated projective
bimodule $M$. In contrast, \textit{wheeled contractions} refer to the
contractions $\mu_{i,j}$  that are part of the structure of wheelspaces.

Suppose we are given an $A$-bimodule, $M$, and an $A$-bimodule
connection $\nabla=(\nabla_\ell, \nabla_r)$.  We would like to say
that $\nabla$ extends to an operator
$\Ff_m(T_A M) \rightarrow \Ff_m(T_A (M \oplus \Omega^1 A))$, which
commutes with wheeled contractions and permutations, and lands in the
subspace of degree one in $\Omega^1 A$.  It suffices to define the
operator $\nabla$ on $T_A M, T_A M/[T_A M, T_A M]$ and extend by
\begin{equation}
\nabla(a \o b) = \nabla(a) \o b + a \o \nabla(b), \quad \nabla(\sigma \o a) = \sigma \o \nabla(a).
\end{equation}
We have
\begin{equation}
\nabla(m_1 \o_A \o \cdots \o_A m_n) = \sum_{i=1}^n (m_1 \o_A \cdots
\o_A (\nabla_\ell + \nabla_r)(m_i) \o_A \cdots \o_A m_n).
\end{equation}
and we obtain the formula on $\Ff_0(T_A M)$ by taking a wheeled contraction of this.

It is not difficult to check that $\nabla$ extends to an operator on
$\Ff(T_A (M \oplus \Omega^1 A))$, which commutes with wheeled contractions and
permutations, as follows. Let us use superbraidings with respect to a bigrading
with $|A| = (0,0)$, $|\Omega^1 A|=(0,-1)$, and $|M| \subset \Z \times \{0\}$.  
\begin{lemdfn}
  The operator $\nabla: T_\k(M \oplus A \oplus \Omega^1 A) \rightarrow \Ff(M \oplus \Omega^1 A)$
  given by
  $\nabla: f_1 \o \cdots \o f_m \rightarrow \sum_{i=1}^m
  (-1)^{|\{j < i: f_j \in \Omega^1 A\}|}\nabla^{(i)}(f_1 \o \cdots \o f_m)$,
  where 
\begin{equation}
\nabla(f) := \begin{cases} \nabla_\ell(f) + \nabla_r(f), & a \in M, \\
                           df, & f \in A \text{ or } \Omega^1 A,
\end{cases}
\end{equation}
extends uniquely  to an endomorphism of wheelspaces of
$\Ff(M \oplus \Omega^1 A)$. Hence, we define $\nabla$ in this way.
\end{lemdfn}

Next, we need to define bimodule contractions.  Suppose that
$M, M^\vee$ are dual finitely-generated projective $A$-bimodules
(i.e., $M^\vee \iso \Hom_{A^e}(M, A^e)$).  Then we have a map
$i_{pre}: M^{\o n} \o M^{\vee} \rightarrow M^{\o (n-1)} \o (A \o A)$,
by contracting the last two factors.  This map is $(A^e)^{\o
  (n+1)}$-linear,
sending the $A^e$-action on the $m$-th component to the $A^e$-action
on the $m$-th component on the right hand side for $m \leq n$, if one
considers the $(A \o A)$ term as a single component with outer action;
then, the $A^e$ action on $M^\vee$ gets sent to the inner $A^e$-action
on the $A \o A$ term.

From this, we may consider a contraction
\begin{equation}
(\Ind_{S_m}^{S_m \times S_m} M^{\o m}) \o (\Ind_{S_n}^{S_n \times
  S_n}
(M^{\vee})^{\o n}) \rightarrow (\Ind_{S_{m+n}}^{S_{m+n} \times
  S_{m+n}} (M^{\o (n-m)} \o A^{\o 2m})),
\end{equation}
if $m \leq n$, by total contraction, as follows.  By
$\Ind_{S_m}^{S_m \times S_m}$, we mean induced with respect to the
diagonal embedding, and we will write the action of $S_m \times S_m$
on an element $x$ with the notation $\sigma \cdot x \cdot \tau$ for
$\sigma, \tau \in S_m$.
Now, for any
$\sigma_L  (a_1 \o \cdots \o a_m) \sigma_R$ and $\sigma_L' (a_1' \o \cdots \o a_n') \sigma_R'$, 
\begin{enumerate}
\item Tensor the two elements, giving $(\sigma_L \times \sigma_L') (a_1 \o \cdots \o a_m \o a_1' \o
\cdots \o a_n') (\sigma_R \times \sigma_R')$.
\item Summing over
  all ways to rewrite the above expression by applying a 
cyclic permutation to $\sigma_L$ and cyclically permuting the
$a_i$, contract the two adjacent
  elements $M \o M^{\vee}$ that appear, sending the $A^e$-bimodule structure on
  $M$ to the outer structure (and the $A$-bimodule structure on
  $M^{\vee}$ to the inner structure) on $A \o A$. This yields
\begin{multline}
\sum_{i=1}^{m} 
(12 \cdots m)^{-i} (\sigma_L \times \sigma_L') \o (a_{i+1} \o a_{i+2} \o \cdots \o a_{i-1} \o (a_1'(a_i)) \\ \o a_2' \o \cdots \o a_n')
\o (m,m+1)(12 \cdots m)^{-i}(\sigma_R \times \sigma_R').
\end{multline}
\item Rewrite the above by moving the ($A \otimes A$)-components all the way to the right, and adding appropriate permutations to the left and right of the above expression.
Continue iterating
until one has an element of 
$\Ind_{S_{n+m}}^{S_{n+m} \times S_{n+m}} (M^{\o n-m} \o A^{\o 2m})$.
\end{enumerate}
\begin{prop} \label{tcprop}
\begin{enumerate}
\item[(i)]The above procedure is well defined, and using wheeled contractions,
yields \textbf{total contraction maps}
\begin{multline}
  i_{tot}: \Ff_m(T_A M) \o \Ff_n(T_A M^{\vee}) \rightarrow \Ff_{n+m}(T_A M) \oplus_{\Ff_{n+m}(A)} \Ff_{n+m}(T_A M^{\vee}) \\ \subset \Ff_{n+m}(T_A(M \oplus M^{\vee})),
\end{multline}
which respect the grading,
$|M|=1, |A|=0, |M^{\vee}| = -1$, and the $(A^e)^{\o (m+n)}$-module structure.
\item[(ii)] The above map extends, for any $N, N'$, to a map 
$$\Ff_m(T_A(M\oplus N)) \o \Ff_n(T_A(M^\vee \oplus N'))
\rightarrow \Ff_{m+n}(T_A (M \oplus M^\vee \oplus N \oplus N')),$$
that respects the above grading and
has image contained in the sum of tensors which have only occurrences of $M$ or $M^{\vee}$, but not both.
\end{enumerate}
\end{prop}

Given $f \in \Ff_n(T_A M^{\vee})$, let $i_f: \Ff_m(T_A M) \rightarrow
\Ff_{n+m}(T_A M) +
 \Ff_{n+m}(T_A M^\vee)$ be the
map $i_{tot}(- \o f)$, and denote the induced map $\Ff_m(T_A(M \oplus
N)) 
\rightarrow \Ff_{n+m}(T_A(M \oplus M^\vee \oplus N))$ by $i_f$ as well.

Finally, we define trace. Let $M, M_1, M_2$ be any $A$-bimodules.
\begin{defn} \label{trdefn} Given an  $A$-bimodule map $
\phi: M \rightarrow M_1 \o_A M \o_A M_2,
$ where $M$ is projective,  let  $\phi' \in M^\vee \o_{A^e} (M_1 \o_A M
\o_A M_2)$
be the element corresponding to $\phi$. We put
 $\tr(\phi) := i_{\iota}(\phi')$, where
$\iota \in M^\vee \o_{A^e} M$ is the canonical element.
\end{defn}

\subsection{The differential operator $D_\nabla$}
We are now prepared to define the map $D_\nabla$.  Consider
the canonical element $\iota$ from Definition \ref{iodfn}. Let us now
consider $\iota$ as an element of $\DDA^\vee \otimes_{A^e} (\Omega^1 A)^\vee
 \subset \Ff_0(T_A (\DDA \oplus \Omega^1 A)^{\vee})$.  We then have the obtained
contraction $i_\iota$ as in Proposition \ref{tcprop}.
Gradings are as in Section \ref{signntns}.


\begin{defn} \label{dndfn}
We define an operator $D_\nabla: \Ff_m(T_A \DDA) \rightarrow \Ff_m(T_A \DDA)$ by 
$
D_\nabla = i_{\iota} \circ \nabla.
$
\end{defn}
The map $D_\nabla$ satisfies BV-like identities involving the double
Poisson bracket, analogous to \eqref{cdnsb}. We have the following main
result (some notation will be explained after the statement).
\begin{thm} \label{bvthm} \begin{enumerate}
Let $\nabla$ be any torsion-free bimodule connection.
\item[(i)] $D_\nabla$ is a differential operator of order $\leq 2$ and
  degree $-1$ (for $|A|=1, |\DDA|=1$) on $\Ff(T_A \DDA)$, commuting with
  wheeled contractions, 
whose principal symbol $\Gamma_2(D_\nabla)$ is $\pm$ the
  Schouten-Nijenhuis bracket.  That is, one has the BV identity, for
homogeneous $\xi, \eta \in \Ff(T_A \DDA)$,
\begin{equation}\label{ncbvid}
(-1)^{|\xi|+1} \{\xi, \eta\} = D_\nabla(\xi \o \eta) - D_\nabla(\xi)
\o \eta - (-1)^{|\xi|} \xi \o D_\nabla(\eta).
\end{equation}
\item[(ii)] The operator $\nabla^2: \DDA \rightarrow \Ff_1(T_A(\DDA \oplus \Omega^1 A))$
is $A^e$-linear, and one has
\begin{equation} \label{ncdnsf}
D_\nabla^2 = i_{\tr(\nabla^2)}.
\end{equation}
\item[(iii)] More generally, if we adjoin formally the element
  $(1 + \rk(\DDA))^{-1}$ to $\Ff(A)$ (for $1 \in \Ff_0(A)$ the
  wheelgebra unit, not the unit element of $A$), then the above yields
  a one-to-one correspondence between generalized torsion-free
  bimodule connections $\nabla = (\nabla_\ell, \nabla_r)$ on $\DDA$
  and differential operators $D \in \Df_0(A)_{\leq 2}[-2]$ whose
  principal symbol is the Schouten-Nijenhuis bracket.  Under this
  correspondence, the trace-flat connections are the ones which map to
  wheeled BV operators.  Two torsion-free generalized connections
  $\nabla, \nabla'$ map to the same differential operator
  $D_\nabla = D_{\nabla'}$ if and only if $\div \nabla = \div
  \nabla'$.
\end{enumerate}
\end{thm}
Here, $\tr(\nabla^2)$ is defined, as in Definition \ref{trdefn}, by
contracting the input with the $\DDA$ in the output, but viewing the
output $\Omega^1 A$'s as separate from $\DDA$ (like $M_1, M_2$ in
Definition \ref{trdefn}).  The notation $\Df_0(A)_{\leq 2}[-2]$ means
sending $\Ff_m(A)$ to $\Ff_m(A)$, having order $\leq 2$, and degree
$-2$ using $|\DDA|=1$, $A$=0.

A \textit{generalized connection} $\nabla = (\nabla_\ell, \nabla_r)$ is
a pair of maps
$\ds
\nabla_\ell, \nabla_r: \DDA \rightarrow \Omega^1 A \o_{A^e} \WDDA_2,
$
such that 
\begin{gather}
\nabla_\ell(a \xi) = a \o_A \nabla_\ell(\xi) + da \o_A \xi, \quad \nabla_\ell(\xi a) = \nabla_\ell(\xi) \o_A a, \\
\nabla_r(a \xi) = a \o_A \nabla_r(\xi), \quad \nabla_r(\xi a) = \nabla_r(\xi) \o_A a + \xi \o_A da.
\end{gather}
One may define the torsion exactly as in Section \ref{tcddas}, and
hence the torsion-free condition. 

By formally adjoining the element $(1 + \rk(\DDA))^{-1}$, we mean to
consider polynomials in this element (considered to have wheeled
degree zero) with coefficients in $\Ff(A), \Df(A)$; the element
$\rk(\DDA) \in \Syns^2 (A/[A,A]) \subset \Ff_0(A)$ is
$\rk(\DDA) := \pi(i_{\iota}(\iota))$, where
$\pi: (A/[A,A])^{\otimes 2} \onto \Syns^2 (A/[A,A])$ is the
projection.\footnote{Note that the element
  $i_{\iota}(\iota) \in (A/[A,A])^{\otimes 2}$ itself is already
  canonical, and under the representation functor, it is taken to the
  dimension of the representation varieties.}
Finally, $\div \nabla := i_{\iota} \circ \nabla$.

\begin{cor} If $\nabla$ is torsion-free, then $D_\nabla$ endows
  $\Ff(T_A \DDA)$ with the structure of a BV wheelgebra if and only if
  $\tr(\nabla^2) = 0$.  The induced wheeled Gerstenhaber structure is the
  Schouten-Nijenhuis one.
\end{cor}
In particular, this motivates the definition
\begin{defn} An associative algebra $A$ is wheeled Calabi-Yau if $A$ has
a trace-flat torsion-free bimodule connection.
\end{defn}
By the above, if $A$ is wheeled Calabi-Yau, then $\Ff(T_A \DDA)$ is
equipped with a wheeled BV structure extending the Schouten-Nijenhuis
wheeled Gerstenhaber structure.  (Also, by the above remarks, the two
notions are not that different).
\begin{cor} \label{fbvc} Suppose $\nabla$ is a torsion-free
  connection on $\DDA$. Then, $D_\nabla$ induces maps
\begin{align*}
D_\nabla:\
&T_A \DDA \o T_A \DDA \too T_A \DDA \o T_A \DDA \o \Syns^{\leq 1} A_\cyc, \quad\text{and}\\
&T_A \DDA \too T_A \DDA \o (T_A \DDA)_\cyc,
\end{align*}
which satisfy the identities, for $\xi, \eta \in T_A \DDA$,
\begin{gather} \label{fbvet}
D_\nabla(\xi \o \eta) - D_\nabla(\xi) \o \eta - (-1)^{|\xi|} \xi \o D_\nabla(\eta) =
(-1)^{|\xi| + 1} \{ \xi, \eta \}, \\
\label{fbvep}
D_\nabla(\xi\eta) - D_\nabla(\xi) \eta - (-1)^{|\xi|} \xi D_\nabla(\eta) = (-1)^{|\xi| + 1} (\pr \o 1)\,
\{ \xi, \eta \}, 
\end{gather}
where $\pr:
T_A \DDA \onto (T_A \DDA)_\cyc$ is the projection. 
\end{cor}

\begin{proof}[Proof of Corollary \ref{fbvc}.]
  We will prove the corollary independently of the theorem, to help
  explain what is going on in a simpler setting. Recall that
  $D_\nabla = i_\iota \circ \nabla$, and the contraction $i_\iota$
  with the canonical element is, essentially by definition, a
  signed sum over ways of contracting a copy of $\DDA$ with
  $\Omega^1 A$.  Hence, if $\xi = \xi_1 \o_A \o \cdots \o_A \xi_m$ and
  $g = \eta_1 \o_A \eta_2 \o_A \cdots \o_A \eta_n$, the LHS of
  \eqref{fbvet} can be expanded as
\begin{multline}
\sum_{i=1}^m \sum_{j=1}^n \pm \Bigl(\bigl(\eta_1 \o_A \o \cdots \o_A \eta_{j-1}
\o_A \bigl(( (\nabla_\ell)_{\xi_i} \eta_j)' - ((\nabla_r)_{\eta_j} \xi_i)''\bigr) 
\o_A \xi_{i+1} \o_A \cdots \o_A \xi_m \bigr) 
\\ \o \bigl(\xi_1 \o_A \cdots \o_A \xi_{i-1} \o_A \bigl(
( (\nabla_\ell)_{\xi_i} \eta_j)'' - ((\nabla_r)_{\eta_j} \xi_i )'\bigr) \o_A \eta_{j+1} \o_A \cdots \o_A \eta_n\bigr) \\
+ \pm \bigl( \eta_1 \o_A \o \cdots \o_A \eta_{j-1}
\o_A \bigl( ((\nabla_r)_{\xi_i} \eta_j)' - ((\nabla_\ell)_{\eta_j} \xi_i)''\bigr) 
\o_A \xi_{i+1} \o_A \cdots \o_A \xi_m \bigr)
\\ \o \bigl(\xi_1 \o_A \cdots \o_A \xi_{i-1} \o_A \bigl(
((\nabla_r)_{\xi_i} \eta_j)'' - ((\nabla_\ell)_{\eta_j} \xi_i )' \bigr)
 \o_A \eta_{j+1} \o_A \cdots \o_A \eta_n\bigr) \Bigr),
\end{multline}
where above we use the restriction/quotient described in the statement
of the proposition, and where signs are determined in the appropriate
way from the superbraidings. By \eqref{torfla1} (or \eqref{torfla} in
the case $m$ is odd), since $\nabla$ is torsion-free, \eqref{fbvet}
follows from the graded double Poisson condition for Van den Bergh's
Schouten-Nijenhuis bracket $\ldb -,-\rdb$ (or, equivalently, the
wheeled Poisson condition for $\{-,-\}$), from which we deduce
\eqref{fbvep} using appropriate wheeled contractions.
\end{proof}
\subsection{The quiver case}
Before proving the theorem, we explain what it says in the case of
quivers.  Let $A = P_Q$, the path algebra of a quiver $Q$.  In this
case, it is well known (\cite{BLB, G}) that $(P_{\dq})_\cyc$ is a Lie
algebra, so one may consider its universal enveloping algebra, which
quantizes $\Sym P_{\dq} \cong \Ff_0(T_A \DDA)$.  However, this is
\emph{not} the quantization $\Df_0(\Ff(T_A \DDA))$: this is a
quantization of the \emph{Lie bialgebra} structure defined in
\cite{S}, called the \emph{quantized necklace algebra}.  So, in a
sense, taking differential operators performs two quantizations at
once: the universal enveloping, followed by the quantized universal
enveloping.
\begin{thm}\label{quivthm}
\begin{enumerate}
\item[(i)] $A$ is wheeled Calabi-Yau, equipped with the trivial connection on
  $\DDA$, whose associated wheeled BV structure incorporates the
  double/necklace Lie (co)bracket as in \eqref{neckbv}.
\item[(ii)] The degree-zero part $\Df_0(\Ff(P_Q))$ of
wheeled differential operators on $P_Q$ is isomorphic to the quantized
necklace algebra from \cite{S} (for $\hbar = 1$).
\item[(iii)] The action of $\Df_0(\Ff(P_Q))$ on $\Ff_0(P_Q)$ is the
  limit, as $\mathbf{d} \rightarrow \infty$, of the representations of
  the quantized necklace algebra as differential operators on $\Rep_{\mathbf{d}}(Q)$ (for $\hbar = 1$).
\item[(iv)] One has a wheeled analogue of the infinitesimal Weil representation: the
quadratic-in-path-length subspace of $(P_{\dq})_\cyc$ is a Lie subalgebra
of $\Df(\Ff(P_Q))$ isomorphic to $\mathfrak{sp}(Q)$, which acts via a
\textit{wheeled infinitesimal Weil representation} on $\Ff(P_Q)$.
\end{enumerate}
\end{thm}
Note that, in particular, the degree-zero part of a BV wheelgebra is
an ordinary BV algebra (indeed, the degree-zero part of any type of
wheelgebra is an ordinary algebra of that type), and so the algebra
$\Ff_0(P_{\dq}) = \SuperSym (P_\dq)_\cyc$, in particular obtains an
ordinary BV algebra structure. This coincides with the one obtained
from the construction of \S \ref{nlbasec}, viewing $(P_\dq)_\cyc$ as
an involutive Lie bialgebra.
\begin{exam}
  In the case that the quiver has one vertex, $A = TV$ is merely a
  tensor algebra on a vector space $V$.  In this case, the BV algebra
  $\SuperSym (TV)_\cyc$ coincides with the BV algebra $F$ considered
  in \cite{Barmobvg}, \cite[\S 1]{Barncbvg} (cf.~\S
  \ref{futdirsec}), although there $V$ is allowed to be a supervector
  space rather than merely a vector space (we could also work in this
  generality, or in the quiver setting, we could allow edges to be
  even or odd).
\end{exam}
\begin{proof}
  (i) The algebra $T_A \DDA$, as was noticed in \cite{VdB}, is
  $P_{\dq}$, equipped with the grading $|Q| = 0, |Q^*| = 1$.  There is
  a canonical bimodule connection $\nabla$ on $\DDA$ given by
  $\nabla(e^*) = 0$ for all $e^* \in Q^*$, called the \emph{trivial
    connection}.  Using this connection, the operator $D_\nabla$ acts
  as in \eqref{neckbv}.  It is immediate that $\nabla^2 = 0$ (since
  $\nabla(e^*) = 0$ for all $e^* \in Q^*$, and $\nabla^2$ is
  $A^e$-linear).  Hence, it follows that $P_Q$ is wheeled Calabi-Yau,
  and that the operator $D_\nabla$ is wheeled BV.

(ii) This follows because $\DDA$ is a free bimodule generated by $\partial_e$ for $e \in Q$, using Theorem \ref{dathm} and the description of the quantized
necklace algebra from \cite{S}.

(iii) This follows from asymptotic bijectivity (\cite{S,GS,EG}) of 
the representation functor $\Ff_0(P_{\dq}) \rightarrow \Rep_{\mathbf{d}}(\dq)^{\mathsf{GL}_{\mathbf{d}}}$, together with Theorem \ref{dathm} and the fact that
$\Ff(P_\dq) \cong \Ff(T_{P_Q} \DDer(P_Q))$ (using
that $\DDA$ is free as in (ii)).

(iv) We claim that the map $[P_{\dq}[2]]_{\cyc} \into \Df_0(\Ff(A))$ given
by symmetrization,
\begin{equation}
[e e^*] \mapsto \frac{1}{2} ([e \circ \partial_e]_\cyc + [\partial_e \circ e]_\cyc), \quad e \in Q,
\end{equation}
is an embedding of Lie algebras (using the Lie structure on
$[P_{\dq}[2]]_\cyc$ obtained from the symplectic pairing of $\langle
\dq \rangle$), where $e \circ \partial_e$ says to first perform
$\partial_e$ and then multiply by $e$ (on the input and output
corresponding to $\partial_e$), and similarly for $\partial_e \circ
e$.  All other elements $[e f], [e^* f^*]$ for all $e,f \in Q$, and
$[e f^*]$ for $f \neq e$, can be mapped by the same symmetrization
procedure, but we don't need to symmetrize since $e, f$
twisted-commute (and similarly $\partial_e, \partial_f$ and
$e, \partial_f$). This claim is easy to verify explicitly, and it is
all we need.
\end{proof}

\subsection{Proof of Proposition \ref{cdnsfp}}
\label{pfcdnsf} In order to prove the theorem, we first need to
complete the proof of Koszul's Proposition \ref{cdnsfp}, as promised.

We prove a more general result than Proposition \ref{cdnsfp}.  Namely,
we give a proof that works in the algebraic setting, as well as the
smooth or complex settings.  It is enough to assume that our variety $X$
is affine, so that $\O_X$, $T_X$, and $\Omega_X$ are generated by local
sections, since the statements (particularly \eqref{cdnsf}) are local.

Now, assuming that $X$ is smooth affine and finite-dimensional, we
give a global argument, that does not rely on further localization (we
will be very careful about this at every step).  Let
$A = \Gamma(X, \O_X)$, $\Der A = \Gamma(X, T_X)$, and
$\Omega^1 A = \Gamma(X, \Omega_X)$. Since $X$ is smooth, $\Der A$ and
$\Omega^1 A$ are finitely-generated projective. Also, $\Omega^1 A$ and
$\Der A$ are projectively dual (i.e.,
$\Omega^1 A \cong \Hom_A(\Der A, A)$ and
$\Der A \cong \Hom_A(\Omega^1 A, A)$).

Let us explicitly
write the canonical element $\iota \in \Der A \o_A \Omega^1 A$
corresponding to the identity as
\begin{equation} \label{ioxo}
\iota = \sum_{i=1}^n \xi_i \o_A \omega^i,
\end{equation}
for some $\xi_i \in \Der A, \omega^i \in \Omega^1 A$. (NOTE: the ``$i$'' appearing
in $\omega^i$ is a
superscript, \emph{not} an exponent.)

Let $\nabla$ be a connection on $T_X$ which is torsion-free.  Explicitly,
\begin{equation}
  \nabla_{\xi_i} \xi_j - \nabla_{\xi_j} \xi_i = [\xi_i, \xi_j], \quad \forall i,j.
\end{equation}
\begin{ntn} Let $\circ$ denote the action of vector fields on
  functions: $\xi \circ f := \xi \ctr df$ for all functions $f \in A$ and
  vector fields $\xi \in \Der A$.
\end{ntn}

We now wish to verify explicitly the identity
\begin{equation}
D_\nabla^2 = -i_{\tr(\nabla^2)}.
\end{equation}
Since we already know (as in Section \ref{cdns}) that the LHS is a
differential operator of order $\leq 2$ and degree $-2$, it must be
given by contraction with some two-form.  To verify that the above
formula holds, it suffices to compute both sides applied to an
arbitrary element of the form $\xi \wedge \eta$, for
$\xi, \eta \in \Der A$.

\begin{ntn}
  For any vector $\eta \in \Der A$, let us define coefficients
  $\eta^i \in A$ by $\eta^i := \eta \ctr \omega^i$ (when there is no
  possible confusion). In particular, this implies that
$
\eta = \sum_i \eta^i \xi_i.
$
\end{ntn}
Whenever we have a tensor $a \in V \o W$, we write $a = a' \o
a''$, using Sweedler's notation (\ref{swntn}).

Finally, note that one has the well known formula for $\nabla^2$, which
we prove as a warmup (and because we will need to prove this in the wheeled setting):
\begin{lemma}
\begin{equation}
(\xi \wedge \eta) \ctr \nabla^2 = \nabla_\xi \nabla_\eta - \nabla_\eta \nabla_\xi - \nabla_{[\xi, \eta]}.
\end{equation}
\end{lemma}
\begin{proof}
We have
\begin{equation}
\nabla(\theta) = \sum_i \nabla_{\xi_i} \theta \o_A \omega^i.
\end{equation}
Then, we have
\begin{equation} \label{n21}
\nabla^2(\theta) = \sum_{i,j} \nabla_{\xi_i} \nabla_{\xi_j} \theta \o_A (\omega^i \wedge
\omega^j) + \sum_i \nabla_{\xi_i} \theta \o_A d \omega^i.
\end{equation}
We need to use the following formula (which we note can be proved by writing $\omega = a\, db$ without localizing):
\begin{equation}
d\omega \ctr (\xi \wedge \eta) = \xi \circ (\omega \ctr \eta) - \eta \circ (\omega \ctr \xi) - [\xi, \eta] \ctr \omega.
\end{equation}
Thus, we may expand \eqref{n21} applied to $\xi \wedge \eta$, for each of the summations
on the RHS separately, as
\begin{multline}\label{n22}
(\xi \wedge \eta)\ctr \sum_{i,j} \nabla_{\xi_i} \nabla_{\xi_j} \theta (\omega^i \wedge \omega^j) 
= \sum_{i,j} (\xi^i \eta^j - \xi^j \eta^i) \nabla_{\xi_i} \nabla_{\xi_j} \theta
\\ = \sum_{j} \eta^j \nabla_{\xi} \nabla_{\xi_j} \theta - \xi^j \nabla_{\eta} \nabla_{\xi_j} \theta \\
= \bigl(\nabla_{\xi} \nabla_{\eta} - \nabla_{\eta} \nabla_{\xi}\bigr) \theta - 
\bigl( (\xi \circ \eta^j) - (\eta \circ \xi^j) \bigr) \nabla_{\xi_j} \theta.
\end{multline}
\begin{multline}\label{n23}
\sum_{i} \nabla_{\xi_i} \theta \o_A d\omega_i \ctr (\xi \wedge \eta) = \sum_i \bigl( (\xi \circ \eta^i) - (\eta \circ \xi^i) - [\xi, \eta]^i\bigr) \nabla_{\xi_i} \theta \\
= -\nabla_{[\xi, \eta]} \theta + \sum_i \bigl( (\xi \circ \eta^i) - (\eta \circ \xi^i) \bigr) \nabla_{\xi_i} \theta.
\end{multline}
Summing \eqref{n22} and \eqref{n23}, we deduce the desired identity.
\end{proof}

We will need the formula
\begin{equation}
\tr(F) = \sum_i F(\xi_i)^i,
\end{equation}
for any endomorphism $F: \Der A \rightarrow \Der A$.  It is immediate that this is the
same as applying the contraction $i_{\iota}$ to $F$ considered as an element of $\Der A \o_A \Omega^1 A$.

Now, we proceed to the main
\begin{lemma}
\begin{equation}
D_\nabla^2(\xi \wedge \eta) = -\tr(\nabla_\xi \nabla_\eta - \nabla_\eta \nabla_\xi
- \nabla_{[\xi, \eta]}).
\end{equation}
\end{lemma}
\begin{proof}
First, by applying $D_\nabla$ twice, and canceling the $\pm D_\nabla(\xi) D_\nabla(\eta)$
terms, we obtain
\begin{equation}
D_\nabla^2(\xi \wedge \eta) = \eta \circ D_\nabla(\xi) - \xi \circ D_\nabla(\eta) +
D_\nabla [\xi, \eta].
\end{equation}
By definition, the RHS expresses in terms of generators as
\begin{equation}
\sum_i \eta \circ (\nabla_{\xi_i}(\xi)^i) - \xi \circ(\nabla_{\xi_i}(\eta)^i)
+ \nabla_{\xi_i}([\xi, \eta])^i. 
\end{equation}
Next, we apply torsion-freeness to rewrite the RHS as
\begin{equation} \label{dn21}
\sum_i \eta \circ (\nabla_{\xi}(\xi_i)^i - [\xi, \xi_i]) - \xi \circ
(\nabla_{\eta}(\xi_i)^i - [\eta, \xi_i]) + \nabla_{[\xi, \eta]}(\xi_i)^i -
[[\xi, \eta], \xi_i]^i.
\end{equation}
Furthermore, we note that
\begin{multline}\label{dn22}
\nabla_{\xi} \nabla_\eta(\xi_i)^i = \sum_{\ell} \nabla_\xi (\nabla_\eta(\xi_i)^\ell \xi_\ell)^i
= \sum_\ell \xi \circ (\nabla_\eta(\xi_i)^\ell) (\xi_\ell \ctr \omega^i) + \nabla_\xi(\xi_\ell)^i \nabla_\eta(\xi_i)^\ell 
\\ = \sum_\ell \xi \circ \bigl((\xi_\ell \ctr \omega^i) \nabla_\eta(\xi_i)^\ell \bigr) - \xi \circ (\xi_\ell \ctr \omega^i) \nabla_\eta(\xi_i)^\ell + \nabla_\xi(\xi_\ell)^i \nabla_\eta(\xi_i)^\ell \\
= \xi \circ (\nabla_\eta(\xi_i)^i) + \sum_\ell \bigl(\nabla_\xi(\xi_\ell)^i -\xi \circ (\xi_\ell \ctr \omega^i) \bigr) \nabla_\eta(\xi_i)^\ell.
\end{multline}
Substituting this formula (and the same with $\xi, \eta$ swapped) into \eqref{dn21},
 we obtain
\begin{multline}\label{dn23}
D_\nabla^2(\xi \wedge \eta) + \tr(\nabla_\xi \nabla_\eta - \nabla_\eta \nabla_\xi - \nabla_{[\xi, \eta]}) \\ = \sum_i -\eta \circ [\xi, \xi_i] + \xi \circ [\eta, \xi_i] - [[\xi, \eta], \xi_i]^i
\\ + \sum_{i, \ell} (\eta \circ (\xi_\ell\ctr\omega^i)) \nabla_\xi(\xi_i)^\ell
- \xi \circ (\xi_\ell\ctr\omega^i)\nabla_\eta(\xi_i)^\ell.
\end{multline}
It remains to prove that the RHS is zero.  Using the Jacobi identity, we have
\begin{multline} \label{dn24}
[[\xi, \eta], \xi_i]^i = [\xi,[\eta, \xi_i]]^i - [\eta,[\xi, \xi_i]]^i
 = \sum_\ell [\xi, [\eta, \xi_i]^\ell \xi_\ell]^i - [\eta, [\xi, \xi_i]^\ell \xi_\ell]^i
\\ = \xi \circ ([\eta, \xi_i]^\ell) (\xi_\ell \ctr \omega^i)  - \eta \circ ([\xi, \xi_i]^\ell)(\xi_\ell \ctr \omega^i) + \sum_\ell [\xi, \xi_\ell]^i[\eta, \xi_i]^\ell - [\eta, \xi_\ell]^i [\xi, \xi_i]^\ell \\
= \xi \circ [\eta, \xi_i]^i - \eta \circ[\xi, \xi_i]^i + \sum_\ell - \xi \circ (\xi_\ell \ctr \omega^i) [\eta, \xi_i]^\ell + \eta \circ (\xi_\ell \ctr \omega^i) [\xi, \xi_i]^\ell.
\end{multline}
Thus, combining \eqref{dn23} and \eqref{dn24}, and applying
torsion-freeness, we need only show that
\begin{equation}
\sum_{i, \ell} (\eta \circ (\xi_\ell\ctr\omega^i)) \nabla_{\xi_i}(\xi)^\ell
- \xi \circ (\xi_\ell\ctr\omega^i)\nabla_{\xi_i}(\eta)^\ell = 0.
\end{equation}
More generally, we claim that
\begin{equation} \label{iellid}
\sum_{i, \ell} d(\xi_\ell\ctr\omega^i) \o_A \xi_i \o_A \omega^\ell = 0.
\end{equation}
To show this, let us define
\begin{equation}
M_{ij} := \xi_i \ctr \omega^j,
\end{equation}
and contract the second and third components of the LHS \eqref{iellid} with
$\omega^j \o_A \xi_k$. We obtain
\begin{equation} \label{mmm}
\sum_{i, \ell} d(M_{\ell i}) M_{ij} M_{k\ell}.
\end{equation}
To show that this is zero (for all $j, k$), first note that
$\xi_i = \sum_j (\xi_i \ctr \omega^j) \xi_j$ implies
\begin{equation} \label{mmp}
\sum_j M_{ij} M_{jk} = M_{ik}.
\end{equation}
Applying \eqref{mmp} together with the Leibniz rule (i.e., $d(M_{\ell i}) M_{ij} =
d(M_{\ell i} M_{ij}) - M_{\ell i} d(M_{ij})$),
we obtain
\begin{equation}
\sum_{i, \ell} d(M_{\ell i}) M_{ij} M_{k \ell} = \sum_\ell d(M_{\ell j}) M_{k \ell}
- \sum_i M_{ki} d(M_{ij}) = 0.
\end{equation}
This concludes the proof.
\end{proof}

\subsection{Proof of Theorem \ref{bvthm}}
The difficulty is in proving part (ii), which we do first, using part
(i).  For this, the proof in the previous subsection applies to our
setting as well, provided we make the following interpretations and
conventions: First, $\xi, \eta$ are now double derivations. We use
\eqref{ndfn}: $\nabla_\xi := (\nabla_\ell)_\xi + (\nabla_r)_\xi$.
Secondly, in computing $D_\nabla^2(\xi \o \eta)$, we can afford to
lose track of (signed) permutations of tensor components, if we apply,
at the end, whatever signed wheeled permutation is necessary so as to
have the substitution $\xi \mapsto a \xi b$ apply left-multiplication
by $a$ under the first $A$-module structure, and similarly apply
right-multiplication by $b$ to the first $A^{\op}$-module structure.
Here we use that each $\Ff_m(T_A M)$ is not merely an $S_m$-module,
but an $S_m \times S_m$-module, given by permuting the left and right
components separately.

Furthermore, the tensor products over $A$ must become tensor products
over $A^e$, and $\wedge$ becomes $\o_A$. Then $\ctr$ becomes the total
contraction.  Finally, in the arguments about $M_{ij}$ to prove
\eqref{iellid}, the Leibniz rule still applies, but now
$M_{ij} \in A \o A$, and the multiplication is $A^e$-multiplication:
if $\xi_i \ctr \omega^j$ takes the bimodule action on $\omega^j$ to outer
action, as in \eqref{paireqn}, then $M_{ij} M_{kl}$ glues the outer
action on $M_{ij}$ to the inner on $M_{kl}$ (i.e., puts $M_{ij}$
inside $M_{kl}$ and multiplies).  

Making all of these changes, the proof of (ii) goes through.

(i) Since $D_\nabla = i_\iota \circ \nabla$, the definitions show that
we are summing over applying operations to two terms in a tensor
product. This shows that $D_\nabla$ must have order $\leq 2$.  Also,
the number of $\DDA$'s (the degree) goes down by one. Then, we apply
Theorem \ref{dathm} and Proposition
\ref{torprop}.

(iii) It is easy to see that a generalized connection is torsion-free
if and only if it satisfies the BV identity \eqref{ncbvid} (we say it ``generates
the S-N bracket'').

We first show that any differential operator $\phi$ that generates the
S-N bracket is of the form $\phi = D_\nabla$ for an appropriate
torsion-free generalized bimodule connection
$\nabla=(\nabla_\ell, \nabla_r)$.  To see this, first let $\nabla'$ be
\textbf{any} torsion-free bimodule connection (which exists because we
can let $\nabla'_\ell$ be arbitrary, and then $\nabla'_r$ is
determined by \eqref{torfla1})).  Then,
$\Gamma_2(\phi - D_{\nabla'}) = 0$, so $\phi - D_{\nabla'}$ is a
differential operator of order one and degree $-1$.  Such a map is the
same as an $A^e$-linear map $\DDA \rightarrow \Ff_1(A)$, i.e.,
contraction with an element $\alpha \in \Ff_0(\Omega^1 A)[-1]$.  It
suffices to show that any such element can be realized as
$D_{\nabla} - D_{\nabla'}$ for some connection $\nabla$.  Finding
such a $\nabla$ is the same as finding $\nabla_\ell - \nabla'_{\ell} \in \WDDA_{2} \o_{A^e \o A^e} (\Omega^1 A \o (\Omega^1 A)')$ (where $(\Omega^1 A)'$ is
a distinctly-labeled copy of $\Omega^1 A$),
which uniquely determines $\nabla_r$ by the torsion-free condition. Precisely, $\nabla_r - \nabla_{r'}$ swaps $\Omega^1 A$ with $(\Omega^1 A)'$.

Next, let us write $\alpha = \sum_i [\alpha_i]_{\cyc} \o X_i$ for some $\alpha_i \in \Omega^1 A$, and some $X_i \in \Ff_0(A)$.  It suffices to find an element
$Y \in \WDDA_{2} \o_{A^e \o A^e} (\Omega^1 A \o \Omega^1 A)$ such that 
\begin{equation}
i_{\iota} (Y) = [\alpha_i]_\cyc.
\end{equation}
There is a natural choice, $Y = [\alpha_i]_\cyc \o \iota$. We compute
\begin{equation}
i_\iota(Y) = [\alpha_i]_\cyc + [\alpha_i]_\cyc \o \rk(\DDer).
\end{equation}
Thus, if we are allowed to divide by $(1 + \rk(\DDer))$ (where
$1 \in \Ff_0(T_A \DDA)$ is the unit of the wheelgebra, not of $A$), then
we can produce such a $\nabla$.

Now, a straightforward generalization of (ii) shows that the operator
$D_\nabla$ associated to a torsion-free generalized bimodule
connection $\nabla$ satisfies $D_\nabla^2 = 0$ if and only if $\tr(\nabla^2) = 0$:
in fact, \eqref{ncdnsf} generalizes to this setting.

\section{The representation functor}\label{repsec}
\subsection{Main constructions}
If $A$ is a finitely-generated associative algebra and $V$ is a
finite-dimensional vector space, one can consider the affine variety
$\Rep_V A$ of representations of $A$ in $V$.  In this section we
sketch how wheeled constructions on $A$ correspond to the usual
notions on $\Rep_V A$: for instance, a wheeled differential operator
on $A$ gives rise to a differential operator on $\k[\Rep_V A]$ for
every $V$.

Throughout, $A$ will denote an associative algebra over $\k$ and $V$ a
finite-dimensional vector space over $\k$. Let
$\Hom_{\k-\text{alg}}(X, Y)$ denote the space of $\k$-algebra
homomorphisms from $X$ to $Y$. Recall from, e.g., \cite[\S 12]{G2} the
\begin{defn}\label{repvadfn}
  The affine $\k$-scheme $\Rep_V(A)$ is defined as the scheme
  representing the functor $B \mapsto \Hom_{\k-\text{alg}}(A, B
  \otimes \End(V))$.
\end{defn}
Next, we define from this a certain commutative wheelgebra.  First,
we define the \emph{endomorphism wheelgebra} of $V$ (which is analogous
to its endomorphism operad, (wheeled) PROP, etc.):
\begin{defn}
  Let $\WEnd(V) := \bigoplus_m \End(V)^{\otimes m}$ be the endomorphism wheelgebra
  of $V$, where the $S_m \times S_m$-module structure is the natural
  one on $(V^*)^{\otimes m} \otimes V^{\otimes m}$ by permuting
  components separately in $V$ and $V^*$, and the contraction
  operation is given by the trace pairings $\tr_{i,j}$ of the $i$-th
  component of $V^*$ with the $j$-th component of $V$.
\end{defn}
\begin{defn}
  Given a (commutative, Lie, Poisson, etc.) wheelgebra $\cW$, and a
  commutative algebra $B$, let $\cW \otimes B$ denote the (commutative,
  Lie, Poisson, etc.) wheelgebra $(\cW \otimes B)(m) := \cW(m)
  \otimes B$, with all structure maps given by tensoring by $\Id_B$.
\end{defn}
The following will be our main object of study:
\begin{defn}
The \emph{representation wheelgebra} $\WRep_V(A)$ is defined as
\begin{equation}
\WRep_V(A) := \k[\Rep_V(A)] \otimes \WEnd(V).
\end{equation}
\end{defn}
This is equipped with canonical evaluation map $\ev: \Ff(A)
\rightarrow \WRep_V(A)$, defined as the unique extension of the
standard evaluation map $\ev: A \rightarrow \k[\Rep_V(A)]
\otimes \End(V)$ to a morphism of wheelgebras.  Recall from, e.g.,
\cite[\S 12]{G2}, that the standard evaluation map is the universal
map such that any algebra map $\rho: A \rightarrow B \otimes \End(V)$
factors through $\ev$ by a morphism of algebras $B \rightarrow
\k[\Rep_V(A)]$. Likewise, its extension has the universal property
that any morphism of wheelgebras $\Ff(A) \rightarrow B \otimes
\WEnd(V)$ factors through $\WRep_V(A)$ via a morphism of the form
$(\phi \otimes \Id): B \otimes \WEnd(V) \rightarrow \k[\Rep_V(A)]
\otimes \WEnd(V)$, where $\phi: B \rightarrow \k[\Rep_V(A)]$ is an
algebra morphism.  The image of the evaluation maps is the space of
$\GL(V)$-invariants: it is clear that the image is contained in the
space of invariants. For surjectivity, we argue similarly to
\cite{LBP}: as $\GL(V)$-representations, $\k[\Rep_V(A)]$ is generated
by $(A \otimes \End(V))$ under the map $A \otimes T \mapsto \tr(\ev(a)
\cdot T)$, which is a direct sum of copies of $V \otimes V^*$. On the
other hand, $\End(V)^{\otimes m} = V^{\otimes m} \otimes
(V^*)^{\otimes m}$, and the fundamental theorem of invariant theory
implies that all $\GL(V)$-invariants are spanned by complete
contractions (matchings of the $V$ components with the $V^*$
components using the canonical pairing).

In fact, the evaluation morphism of wheelgebras gives a natural
presentation of $\k[\Rep_V(A)]$: see \eqref{repvaeq} below.

This motivates the following more general definition:
\begin{defn}\label{repwdfn}
  Given any wheelgebra $\cW$ and vector space $V$, let $(\RA_V(\cW),
  \ev)$ denote the universal associative algebra together with
  evaluation morphism $\ev: \cW \rightarrow \RA_V(\cW) \otimes
  \WEnd(V)$, such that any morphism of wheelgebras $\cW \rightarrow B
  \otimes \WEnd(V)$ factors through $\ev$ via an algebra morphism $B
  \rightarrow \RA_V(\cW)$.  It is called the \emph{representation
    algebra} of $\cW$. Similarly, let $\WRep(\cW) := \RA_V(\cW)
  \otimes \WEnd(V)$.
\end{defn}
In the case of ordinary algebras, instead of wheelgebras, this
definition is equivalent to Definition \ref{repvadfn} for
$\k[\Rep_V(A)]$: the latter is the universal associative algebra
equipped with a morphism $A \rightarrow \k[\Rep_V(A)] \otimes \End(V)$
such that any morphism $A \rightarrow B \otimes \End(V)$ factors
through $\ev$ via an algebra morphism $B \rightarrow \k[\Rep_V(A)]$.

Definition \ref{repwdfn} relies on 
\begin{thm}\label{whrepthm}
  For all wheelgebras, $\RA_V(\cW)$ exists.  It is equipped
  with an action of $\GL(V)$, and the image of the evaluation morphism
  is $\WRep_V(\cW)^{\GL(V)}$. Finally, if $\cW$ is (almost)
  commutative, Poisson, or BV, so is $\RA_V(\cW)$.
\end{thm}
The proof is postponed to the end of the section, since the considerations
are a generalization of the following theorem in the case of $\cW = \Ff(A)$
(which we prove independently of the theorem):

\begin{thm} \label{repthm}
Let $A$ be an associative algebra and $V$ a vector space.
\begin{enumerate}
\item[(i)] The representation algebra $\WRep_V(\Ff(A))$ coincides with
  $\WRep_V(A)$.
\item[(ii)] There is a canonical evaluation morphism of filtered
  wheelgebras, $\ev: \Df(\Ff(A)) \rightarrow (D(\Rep_V(A)) \otimes
  \WEnd(V))$.
\item[(iii)] If $A$ is smooth, then $\RA_V(\Df(\Ff(A))) =
  D(\Rep_V(A))$, the algebra of differential operators on the smooth
  affine scheme $\Rep_V(A)$, and the morphism of (ii) is its
  evaluation morphism.
\end{enumerate}
\end{thm}
\begin{proof}[Sketch of proof.] 
  (i) Given any morphism $\Ff(A) \rightarrow (B \otimes \WEnd(V))$ of
  wheelgebras, the restriction to $A \rightarrow (B \otimes \End(V))$
  factors uniquely through the standard evaluation morphism $A
  \rightarrow (\k[\Rep_V(A)] \otimes \End(V))$ via an algebra map $B
  \rightarrow \k[\Rep_V(A)]$.  Since $B \otimes \End(V)$ generates $B
  \otimes \WEnd(V)$ by wheelgebra operations, the morphism $B
  \rightarrow \k[\Rep_V(A)]$ must in fact factor the original
  wheelgebra morphism through the wheelgebra evaluation morphism
  $\Ff(A) \rightarrow (\k[\Rep_V(A)] \otimes \WEnd(V))$.

  (ii) We have a tautological map of wheelgebras  (where $\otimes_{\text{wh}}$ is the tensor product in the category of wheelspaces, as in \eqref{owheqn})
\begin{equation}\label{dfwhff}
\Df(\Ff(A))
  \otimes_{\text{wh}} \Ff(A) \rightarrow \Ff(A),
\end{equation}
and can compose with the evaluation morphism to obtain a wheelgebra
map $\Df(\Ff(A)) \otimes_{\text{wh}} \Ff(A) \rightarrow \WRep_V(A)$.
We claim that this descends to a morphism $\Df(\Ff(A))
\otimes_{\text{wh}} \WRep_V(A)^{\GL(V)} \rightarrow \WRep_V(A)$.  To
see this, we note that $\k[\Rep_V(A)]$ can be presented as the free
commutative algebra generated by the linear space $A \otimes \End(V)$
subject to the condition that the diagram commutes:
  \begin{equation} \label{repvaeq} \xymatrix{
      A \otimes A \ar[rr]^-{\ev' \otimes \ev'} \ar[d]^{\text{mult}} && (\k[\Rep_V(A)] 
\otimes \End(V))^{\otimes 2} \ar[d]^{\text{mult}} \\
      A \ar[rr]^-{\ev'} && \k[\Rep_V(A)] \otimes \End(V). }
\end{equation}
Here, the morphism $\ev'$ is obtained from the original
morphism $\ev: A \otimes \End(V) \rightarrow \k[\Rep_V(A)]$ using
the trace pairing: $\ev'(a) = a \otimes c$, where $c \in \End(V)
\otimes \End(V)$ is the canonical element inverse to the trace
pairing.  In other words, in a basis of $V$, we can write $\k[\Rep_V(A)]$
as generated by the matrix coefficients $a_{ij}$ of elements $a \in A$,
subject to the linearity relation $\lambda a_{ij} + \mu b_{ij} = (\lambda a+\mu b)_{ij}$ and the matrix multiplication relation $\sum_k a_{ik} a_{kj} = a_{ij}$; the map $\ev'$ then sends $a \in A$ to the canonical matrix whose $i,j$-th coefficient is $a_{ij}$.

As a result, the fact that our original map \eqref{dfwhff} is a morphism of
wheelgebras says that the morphism descends as desired.  Moreover, we can
extend the map to all of $\Df(\Ff(A)) \otimes_{\text{wh}} \WRep_V(A)$, using
that this is all obtainable from the $\GL(V)$-invariants by a linear operation of
the form 
\begin{multline}
\tr_E: (\k[\Rep_V(A)] \otimes \End(V)^{\otimes m})^{\GL(V)} \rightarrow \k[\Rep_V(A)] \otimes \End(V)^{\otimes (m-1)}, \\ T_1 \otimes \cdots \otimes T_m \mapsto \tr(ET_1) T_2 \otimes \cdots \otimes T_m,
\quad E \in \End(V).
\end{multline}
We have to show this is well defined. First, if $E$ is a multiple of
the identity, this is a wheelgebra contraction, so the action commutes
with our morphism.  Next, if $E$ is traceless, then the target can
only be $\GL(V)$-invariant if it is zero.  So, it remains to show that
if $E$ is traceless, then the kernel of $\tr_E$ is preserved by any
wheeled differential operator on $\Ff(A)$.  By the fundamental theorem of
invariant theory, this kernel is spanned by elements of the form
$\ev(f \otimes 1)$, where $f \in \Ff_{m-1}(A)$ and $1 \in A$ is the
identity.  It remains to show that any wheeled differential operator
sends an element of the form $f \otimes 1$ to another element of the
same form.  By induction on order, we can reduce to the case of the
element $1 \in A$ itself (setting $f$ to be the identity of $\Ff(A)$,
of degree zero).  Suppose such a differential operator $\phi$ sends
$1$ to $g \in \Ff(A)$. Then again by induction on order, $\phi(1
\otimes 1) = 1 \otimes g + g \otimes 1 + h \otimes 1 \otimes 1$ for
some $h \in \Ff(A)$.  By the wheelgebra contraction, this implies that
$\phi(1) = 2g + h \otimes 1$, and hence that $g = -h \otimes 1$, as
desired.

It remains to show that the resulting operations on $\WRep_V(A)$ are indeed
wheeled differential operators. This follows by definition.



(iii) By Theorem \ref{dathm}.(ii) and induction, one can reduce to the
statement that $\RA_V(\Ff(T_A \DDA))$ coincides with $\k[T^*
\Rep_V(A)]$. Thanks to part (i), this is the same as showing that
$\Rep_V(T_A \DDA) = T^* \Rep_V(A)$.  This last equality amounts to
showing that $\DDA \otimes_{A^e} (\k[\Rep_V(A)] \otimes \End(V)) \cong
\Der(\k[\Rep_V(A)])$, which follows from the fact that $\Der(A, A
\otimes A) \otimes_{A^e} A = \Der(A, A)$, since $A$ is projective as
an $A$-bimodule.
\end{proof}
As before, the above considerations work equally well in the quiver
setting, i.e., over the semisimple ground ring $\k^I$ rather than
$\k$. We then deduce
\begin{cor} When $A$ is the path algebra of a quiver, the evaluation map
$\Df(\Ff(A)) \rightarrow D(\Rep_V(A)) \otimes \WEnd(V)$ restricts, in degree
zero, to the map $\Df_0(\Ff(A)) \rightarrow D(\Rep_V(A))^{\GL(V)}$ of
\cite{S}, cf.~Theorem \ref{quivthm}. Moreover, applying this to the
odd cotangent bundle, the wheeled BV structure on $\Ff(T_A \DDA)$ (in this case, viewing $\DDA$ as \emph{odd}), induces the standard one on 
$\RA_V(\Ff(T_A \DDA)) = \La T_{\Rep_V(A)}$.  Equivalently, the wheeled
Calabi-Yau structure on $A$ induces the Calabi-Yau structure on $\Rep_V(A)$.
\end{cor}
\begin{proof}[Sketch of proof] The first statement follows by
  comparing explicitly the evaluation morphism above with the
  construction of \cite{S}: the two closely resemble each other.
  There is an easy dictionary between terminology here, such as the
  order of composing partial derivatives, with that of \cite{S}, in
  this case, heights of links over the quiver.  The second statement
  then follows by explicitly comparing the formulas for the BV
  derivative. In some sense, there is little choice for the BV
  derivative other than the one we made, in order for it to induce the
  necklace bracket. More precisely, the choice is only one of which
  connection to use, and the trivial (tautological) one is the choice
  made here and in \cite{S}.
\end{proof}
According to Theorem \ref{quivthm}, in this case we also deduce that
the action of $\Df_0(\Ff(A))$ on $\Ff(A)$ is recovered from the action
of $D(\Rep_V(A))^{\GL(V)}$ on $\Rep_V(A)$ in the limit as $\dim V
\rightarrow \infty$.  In the next section, we will generalize this
arbitrary commutative wheelgebras.

\begin{proof}[Sketch of proof of Theorem \ref{whrepthm}]
  To construct $\RA_V(\cW)$, first form the free associative algebra
  generated by the vector spaces $\cW_m \otimes \End(V)^{\otimes m}$
  for all $m \geq 0$. Along with this, as in \eqref{repvaeq}, form the
  evaluation morphisms $\ev': \cW_m \rightarrow (\cW_m
  \otimes \End(V)^{\otimes m}) \otimes \End(V)^{\otimes m}$ using the
  canonical element in $\End(V)^{\otimes m} \otimes \End(V)^{\otimes
    m}$ inverse to the trace pairing.  Then, as we did after
  \eqref{repvaeq}, we quotient by the condition that $\ev'$
  morphism is a map of wheelgebras.  For the next statement, the
  action of $\GL(V)$ is the one induced by the tautological diagonal
  action on the second factor of the generating spaces $\cW_m
  \otimes \End(V)^{\otimes m}$.  To see that the morphism is
  surjective onto $\GL(V)$-invariants is an application of the
  fundamental theorem of invariant theory, similarly to the case
  $\cW=\Ff(A)$. For the final statement, in the case that $\cW$ is
  (almost) commutative, we deduce that $\RA_V(\cW) \otimes \WEnd(V)$
  is as well, and hence so is $\RA_V(\cW)$. For BV or Poisson
  structures, we have to show that the BV differential or Poisson
  bracket descends (and extends) to $\RA_V(\cW)$.  This follows in the
  same manner as in the proof of Theorem \ref{repthm}.(ii), or
  alternatively, one can apply the result using that the BV
  differential is a differential operator of order two, and a Poisson
  bracket is a differential operator of two inputs, of order one in
  each input.
\end{proof}
\begin{rem}
  The above considerations can also be applied to wheeled PROPs (a
  generalization of commutative wheelgebras), or even wheeled PROs
  (the noncommutative analogue of wheeled PROPs).  For example, let
  $\cW$ be the wheeled associative PROP, whose representations, over
  finite-dimensional vector spaces, are the same as finite-dimensional
  algebras. Then, $\RA_V(\cW)$ is the structure-constant algebra: in
  coordinates $V \cong \k^n$, this is the free commutative algebra
  generated by elements $c_{ij}^k$ for $i,j,k \in \{1, 2, \ldots,
  n\}$, modulo the associativity relation $\sum_{k=1}^n c_{ij}^k
  c_{kl}^m = \sum_{k=1}^n c_{ik}^m c_{jl}^k$.  We interpret the
  elements $c_{ij}^k$ as structure constants for an associative
  algebra on $\k^n$.  Then, an associative algebra structure on $V$
  over a base ring $B$ is the same as an algebra morphism $\RA_V(\cW)
  \rightarrow B$, which is in turn the same as a representation of the
  original wheeled PROP in $B \otimes \WEnd(V)$.  The same applies to
  the wheeled commutative, Lie, etc., PROPs.
\end{rem}

\subsection{The limit $\dim V \rightarrow \infty$}
By a \emph{finitely-generated wheelgebra} we mean a wheelgebra which
is generated, using wheelgebra operations, by a finite-dimensional
$\SS$-module, or equivalently, by a finite-dimensional wheelspace.
\begin{prop}
  Let $\cW$ be a finitely-generated wheelgebra.  Then, as $d
  \rightarrow \infty$, the maps $\ev_d: \cW \rightarrow
  \WRep_{\k^d}(\cW)$ are asymptotically injective: for every
  finite-dimensional vector subspace $U$ of $\cW$, there exists $N \gg
  0$ such that $\ev_d$ is injective restricted to $U$.  Hence, $\cW
  \cong \lim_{d \rightarrow \infty} \WRep_{\k^d}(\cW)^{\GL(\k^d)}$.
\end{prop}
\begin{proof}This follows, as in \cite{G}, from the second fundamental
  theorem of invariant theory, in the same way that the fact that
  $\ev$ is surjective onto the $\GL(V)$-invariant part is a
  consequence of the first fundamental theorem. In more detail, recall
  from the proof of Theorem \ref{whrepthm} that $\RA_V(W)$ is
  generated by elements of the form $\cW_p \otimes \End(V)^{\otimes
    p}$.  Next, restrict $\ev_d$ to $M \otimes \End(V)^{\otimes m}$,
  where $M \subseteq \RA_{V}(\cW)$ is a subspace generated by elements
  of $\cW_p \otimes \End(V)^{\otimes p}$.  It follows from the second
  fundamental theorem of invariant theory that, when $d \geq m+p$,
  $\ev_d$ is injective on $M \otimes \End(V)^{\otimes m}$.  On the
  other hand, $\cW$ is linearly spanned by such subspaces, which
  proves the result.
\end{proof}

{\footnotesize
\bibliographystyle{amsalpha}
\bibliography{nc-bv}

}
\end{document}

%% file: chip.eepic
\setlength{\unitlength}{0.00050000in}
\begingroup\makeatletter\ifx\SetFigFontNFSS\undefined%
\gdef\SetFigFontNFSS#1#2#3#4#5{%
  \reset@font\fontsize{#1}{#2pt}%
  \fontfamily{#3}\fontseries{#4}\fontshape{#5}%
  \selectfont}%
\fi\endgroup%
{\renewcommand{\dashlinestretch}{30}
\begin{picture}(3699,3459)(0,-10)
\path(12,2289)(3687,2289)(3687,1314)
	(12,1314)(12,2289)
\path(537,2289)(537,2964)
\path(912,2289)(912,2964)
\path(3012,2289)(3012,2964)
\path(537,639)(537,1314)
\path(912,639)(912,1314)
\path(3012,639)(3012,1314)
\put(1587,789){\makebox(0,0)[lb]{\smash{{\SetFigFontNFSS{17}{20.4}{\rmdefault}{\mddefault}{\updefault}. . .}}}}
\put(1662,2589){\makebox(0,0)[lb]{\smash{{\SetFigFontNFSS{17}{20.4}{\rmdefault}{\mddefault}{\updefault}. . .}}}}
\put(1287,3189){\makebox(0,0)[lb]{\smash{{\SetFigFontNFSS{12}{14.4}{\rmdefault}{\mddefault}{\updefault}outputs}}}}
\put(1437,114){\makebox(0,0)[lb]{\smash{{\SetFigFontNFSS{12}{14.4}{\rmdefault}{\mddefault}{\updefault}inputs}}}}
\put(1737,1689){\makebox(0,0)[lb]{\smash{{\SetFigFontNFSS{12}{14.4}{\rmdefault}{\mddefault}{\updefault}$w$}}}}
\end{picture}
}

%% file: chipmuij.eepic
\setlength{\unitlength}{0.00050000in}
\begingroup\makeatletter\ifx\SetFigFontNFSS\undefined%
\gdef\SetFigFontNFSS#1#2#3#4#5{%
  \reset@font\fontsize{#1}{#2pt}%
  \fontfamily{#3}\fontseries{#4}\fontshape{#5}%
  \selectfont}%
\fi\endgroup%
{\renewcommand{\dashlinestretch}{30}
\begin{picture}(9237,4014)(0,-10)
\path(5550,2652)(9225,2652)(9225,1677)
	(5550,1677)(5550,2652)
\path(6075,2652)(6075,3327)
\path(6450,2652)(6450,3327)
\path(8550,2652)(8550,3327)
\path(6075,1002)(6075,1677)
\path(6450,1002)(6450,1677)
\path(8550,1002)(8550,1677)
\path(525,2652)(4200,2652)(4200,1677)
	(525,1677)(525,2652)
\path(1050,2652)(1050,3327)
\path(1050,1002)(1050,1677)
\path(3525,1002)(3525,1677)
\path(2625,2652)(2625,3327)
\path(3675,2652)(3675,3327)
\path(2175,1002)(2175,1677)
\path(2625,3327)(2623,3328)(2620,3330)
	(2613,3333)(2603,3339)(2588,3347)
	(2568,3357)(2544,3370)(2514,3386)
	(2479,3404)(2441,3425)(2398,3447)
	(2351,3471)(2302,3497)(2251,3523)
	(2199,3550)(2146,3577)(2093,3604)
	(2041,3630)(1989,3656)(1939,3681)
	(1890,3704)(1843,3727)(1797,3749)
	(1754,3769)(1712,3788)(1671,3806)
	(1632,3823)(1595,3838)(1559,3853)
	(1524,3867)(1490,3879)(1457,3891)
	(1425,3902)(1393,3913)(1362,3922)
	(1331,3931)(1300,3939)(1271,3947)
	(1242,3954)(1212,3960)(1183,3966)
	(1154,3972)(1124,3976)(1094,3980)
	(1064,3983)(1034,3986)(1004,3987)
	(974,3987)(943,3987)(913,3985)
	(882,3982)(851,3978)(821,3972)
	(791,3966)(760,3958)(731,3948)
	(701,3937)(672,3925)(643,3911)
	(615,3896)(587,3879)(560,3861)
	(533,3841)(507,3819)(482,3796)
	(458,3772)(434,3746)(411,3719)
	(389,3690)(367,3660)(347,3628)
	(327,3594)(308,3559)(289,3522)
	(271,3484)(254,3444)(238,3402)
	(225,3367)(212,3331)(200,3293)
	(188,3254)(176,3214)(165,3172)
	(154,3128)(143,3084)(133,3038)
	(123,2991)(113,2942)(104,2892)
	(95,2841)(86,2788)(78,2735)
	(71,2680)(63,2624)(56,2567)
	(50,2510)(44,2452)(38,2393)
	(33,2333)(29,2273)(25,2213)
	(21,2152)(18,2092)(16,2031)
	(14,1971)(13,1910)(12,1850)
	(12,1791)(13,1732)(14,1674)
	(15,1616)(17,1559)(20,1503)
	(23,1449)(27,1395)(31,1342)
	(36,1290)(41,1240)(47,1190)
	(53,1142)(60,1095)(68,1050)
	(76,1005)(84,962)(93,920)
	(102,879)(112,840)(113,839)
	(125,794)(139,749)(153,706)
	(168,664)(184,624)(201,584)
	(219,546)(238,509)(257,473)
	(278,438)(299,405)(321,373)
	(344,342)(368,312)(393,283)
	(418,256)(444,230)(470,206)
	(497,183)(525,161)(553,141)
	(581,123)(610,106)(638,90)
	(667,76)(697,63)(726,52)
	(755,42)(784,34)(813,27)
	(841,22)(870,17)(898,15)
	(926,13)(953,12)(981,13)
	(1007,15)(1034,18)(1060,22)
	(1086,27)(1112,33)(1138,39)
	(1167,49)(1196,59)(1225,72)
	(1255,85)(1284,101)(1314,118)
	(1344,138)(1375,159)(1406,183)
	(1439,209)(1472,237)(1507,268)
	(1543,301)(1580,336)(1618,374)
	(1658,414)(1698,456)(1739,500)
	(1781,546)(1823,592)(1865,639)
	(1906,685)(1946,731)(1984,775)
	(2020,816)(2052,854)(2081,889)
	(2106,918)(2127,943)(2143,963)
	(2156,979)(2165,990)(2171,997)
	(2174,1000)(2175,1002)
\put(7125,1152){\makebox(0,0)[lb]{\smash{{\SetFigFontNFSS{17}{20.4}{\rmdefault}{\mddefault}{\updefault}. . .}}}}
\put(7200,2952){\makebox(0,0)[lb]{\smash{{\SetFigFontNFSS{17}{20.4}{\rmdefault}{\mddefault}{\updefault}. . .}}}}
\put(4650,2052){\makebox(0,0)[lb]{\smash{{\SetFigFontNFSS{12}{14.4}{\rmdefault}{\mddefault}{\updefault}$=$}}}}
\put(1725,2952){\makebox(0,0)[lb]{\smash{{\SetFigFontNFSS{17}{20.4}{\rmdefault}{\mddefault}{\updefault}. . .}}}}
\put(2925,2952){\makebox(0,0)[lb]{\smash{{\SetFigFontNFSS{17}{20.4}{\rmdefault}{\mddefault}{\updefault}. . .}}}}
\put(1350,1227){\makebox(0,0)[lb]{\smash{{\SetFigFontNFSS{17}{20.4}{\rmdefault}{\mddefault}{\updefault}. . .}}}}
\put(2625,1227){\makebox(0,0)[lb]{\smash{{\SetFigFontNFSS{17}{20.4}{\rmdefault}{\mddefault}{\updefault}. . .}}}}
\put(2625,3552){\makebox(0,0)[lb]{\smash{{\SetFigFontNFSS{12}{14.4}{\rmdefault}{\mddefault}{\updefault}i}}}}
\put(2175,552){\makebox(0,0)[lb]{\smash{{\SetFigFontNFSS{12}{14.4}{\rmdefault}{\mddefault}{\updefault}j}}}}
\put(6975,2127){\makebox(0,0)[lb]{\smash{{\SetFigFontNFSS{12}{14.4}{\rmdefault}{\mddefault}{\updefault}$\mu_{i,j} w$}}}}
\put(2175,2127){\makebox(0,0)[lb]{\smash{{\SetFigFontNFSS{12}{14.4}{\rmdefault}{\mddefault}{\updefault}$w$}}}}
\end{picture}
}

%% file: chipperm.eepic
\setlength{\unitlength}{0.00050000in}
\begingroup\makeatletter\ifx\SetFigFontNFSS\undefined%
\gdef\SetFigFontNFSS#1#2#3#4#5{%
  \reset@font\fontsize{#1}{#2pt}%
  \fontfamily{#3}\fontseries{#4}\fontshape{#5}%
  \selectfont}%
\fi\endgroup%
{\renewcommand{\dashlinestretch}{30}
\begin{picture}(8724,3564)(0,-10)
\path(5037,2262)(8712,2262)(8712,1287)
	(5037,1287)(5037,2262)
\path(5562,2262)(5562,2937)
\path(5937,2262)(5937,2937)
\path(8037,2262)(8037,2937)
\path(5562,612)(5562,1287)
\path(5937,612)(5937,1287)
\path(8037,612)(8037,1287)
\path(537,2262)(537,2937)
\path(537,612)(537,1287)
\path(3012,612)(3012,1287)
\path(3162,2262)(3162,2937)
\path(912,612)(912,1287)
\path(12,2262)(3687,2262)(3687,1287)
	(12,1287)(12,2262)
\path(2787,2262)(2787,2937)
\path(2787,2937)(2790,2939)(2796,2944)
	(2807,2953)(2823,2966)(2842,2982)
	(2865,3000)(2888,3020)(2911,3039)
	(2933,3059)(2954,3077)(2973,3095)
	(2989,3111)(3004,3127)(3017,3142)
	(3029,3157)(3040,3172)(3050,3187)
	(3059,3203)(3067,3219)(3076,3236)
	(3084,3255)(3091,3275)(3099,3297)
	(3107,3322)(3115,3349)(3123,3378)
	(3131,3408)(3138,3438)(3145,3466)
	(3151,3491)(3156,3511)(3159,3525)
	(3161,3533)(3162,3536)(3162,3537)
\path(3162,2937)(3159,2939)(3153,2944)
	(3142,2953)(3126,2966)(3107,2982)
	(3084,3000)(3061,3020)(3038,3039)
	(3016,3059)(2995,3077)(2976,3095)
	(2960,3111)(2945,3127)(2932,3142)
	(2920,3157)(2909,3172)(2900,3187)
	(2890,3203)(2882,3219)(2873,3236)
	(2865,3255)(2858,3275)(2850,3297)
	(2842,3322)(2834,3349)(2826,3378)
	(2818,3408)(2811,3438)(2804,3466)
	(2798,3491)(2793,3511)(2790,3525)
	(2788,3533)(2787,3536)(2787,3537)
\path(912,612)(909,610)(903,605)
	(892,596)(876,583)(857,567)
	(834,549)(811,529)(788,510)
	(766,490)(745,472)(726,454)
	(710,438)(695,422)(682,407)
	(670,392)(659,377)(650,362)
	(640,346)(632,330)(623,313)
	(615,294)(608,274)(600,252)
	(592,227)(584,200)(576,171)
	(568,141)(561,111)(554,83)
	(548,58)(543,38)(540,24)
	(538,16)(537,13)(537,12)
\path(537,612)(540,610)(546,605)
	(557,596)(573,583)(592,567)
	(615,549)(638,529)(661,510)
	(683,490)(704,472)(723,454)
	(739,438)(754,422)(767,407)
	(779,392)(790,377)(800,362)
	(809,346)(817,330)(826,313)
	(834,294)(841,274)(849,252)
	(857,227)(865,200)(873,171)
	(881,141)(888,111)(895,83)
	(901,58)(906,38)(909,24)
	(911,16)(912,13)(912,12)
\put(6612,762){\makebox(0,0)[lb]{\smash{{\SetFigFontNFSS{17}{20.4}{\rmdefault}{\mddefault}{\updefault}. . .}}}}
\put(6687,2562){\makebox(0,0)[lb]{\smash{{\SetFigFontNFSS{17}{20.4}{\rmdefault}{\mddefault}{\updefault}. . .}}}}
\put(4137,1662){\makebox(0,0)[lb]{\smash{{\SetFigFontNFSS{12}{14.4}{\rmdefault}{\mddefault}{\updefault}$=$}}}}
\put(1737,837){\makebox(0,0)[lb]{\smash{{\SetFigFontNFSS{17}{20.4}{\rmdefault}{\mddefault}{\updefault}. . .}}}}
\put(912,2562){\makebox(0,0)[lb]{\smash{{\SetFigFontNFSS{17}{20.4}{\rmdefault}{\mddefault}{\updefault}. . .}}}}
\put(5337,1737){\makebox(0,0)[lb]{\smash{{\SetFigFontNFSS{12}{14.4}{\rmdefault}{\mddefault}{\updefault}$((1,2) \times (m,m-1)) w$}}}}
\put(1812,1737){\makebox(0,0)[lb]{\smash{{\SetFigFontNFSS{12}{14.4}{\rmdefault}{\mddefault}{\updefault}$w$}}}}
\end{picture}
}